\documentclass[A4]{amsart}
\usepackage{amssymb,stmaryrd}
\usepackage{amsfonts}
\usepackage{amstext}
\usepackage{algorithmic}
\usepackage{algorithm}
\usepackage{graphicx}
\usepackage{epstopdf}
\usepackage[all]{xy}
\usepackage{MnSymbol}
\usepackage{tikz}

\numberwithin{equation}{section}
\newtheorem{theorem}{Theorem}[section]

\newtheorem{proposition}[theorem]{Proposition}

\theoremstyle{definition}

\theoremstyle{remark}

\newcommand{\C}{{\mathbb{C}}}
\newcommand{\Z}{{\mathbb{Z}}}

\newcommand{\F}{{\mathbb{F}}}

\newcommand{\<}{{\langle}}

\renewcommand{\>}{{\rangle}}

\newcommand{\CF}{{\mathcal{F}}}
\newcommand{\CR}{{\mathcal{R}}}
\newcommand{\CS}{{\mathcal{S}}}

\newcommand{\CQ}{{\mathcal{Q}}}
\newcommand{\wedgeq}{{\wedge\kern-5pt\cdot}}

\newcommand{\tens}{\otimes}

\newcommand{\id}{{\rm id}}

\renewcommand{\o}{{}_{(1)}}
\renewcommand{\t}{{}_{(2)}}

\newcommand{\eps}{\epsilon}

\newcommand{\la}{{\triangleright}}

\def\lcross{{>\!\!\!\triangleleft}}

\begin{document}

\title{Digital quantum groups}
\keywords{Classification, finite-dimensional,  bialgebra, Hopf algebra, quantum group, Fourier, quasitriangular, involutive, factorisable, electronics, finite field.  Ver 1.0}

\subjclass[2010]{Primary 16T05, 16T10, 17B37, 16Q12}

\author{S. Majid and A. Pacho{\l}}
\address{Queen Mary, University of London\\
School of Mathematics \& School of Biological and Chemical Sciences\\
Mile End Rd, London E1 4NS, UK}

\email{s.majid@qmul.ac.uk, a.pachol@qmul.ac.uk}

\date{}

\begin{abstract} We find and classify all bialgebras and Hopf algebras or `quantum groups' of dimension $\le 4$ over the field $\F_2=\{0,1\}$. We summarise our results as a quiver, where the vertices are the inequivalent algebras and there is an arrow for each inequivalent bialgebra or Hopf algebra built from the algebra at the source of the arrow and the dual of the algebra at the target of the arrow. There are 314 distinct bialgebras, and among them 25 Hopf algebras with at most one of these from one vertex to another. We find a unique smallest noncommutative and noncocommutative one, which is moreover self-dual and resembles a digital version of  $u_q(sl_2)$. We also find a unique self-dual Hopf algebra in one anyonic variable $x^4=0$. For all our Hopf algebras we determine the integral and associated Fourier transform operator, viewed as a representation of the quiver. We also find all quasitriangular or `universal R-matrix' structures on our Hopf algebras. These induce solutions of the Yang-Baxter or braid relations in any representation. 
 \end{abstract}
\maketitle 

\section{Introduction}

Quantum groups have been around in modern form since the 1980s, originally arising in quantum inverse scattering but also at the heart of TQFT's leading to knot and 3-manifold invariants on the one hand (the `quasitriangular' Drinfeld-Jimbo quantum groups \cite{Dri}) and key to the first modern models of quantum spacetime as quantum isometry groups (the `bicrossproduct' quantum groups) as in \cite{MaRue}, where they could model quantum gravity effects. The latter quantum groups relate to ideas of Born-reciprocity and observable-state symmetry \cite{Ma:pla}  for quantum gravity. TQFTs and their quantum groups are also behind quantum gravity in 2+1 dimensions with point sources. Even finite-dimensional quantum groups are potentially useful, such as $u_q(g)$ at roots of unity, for example in the Kitaev model of topological quantum computing \cite{Kit} among others. They also provide `re-write rules' in ZX-calculus for quantum computing more generally \cite{CoDun}. Mathematicians first considered Hopf algebras starting in the 1940s, including the development of a theory of characters and Hopf algebra Fourier transform \cite{Ma:book} (but without many examples at that time truly beyond those associated to groups and Lie algebras). This generalises usual Fourier theory and transforms functions on a nonAbelian group to the group algebra viewed as a noncommutative Fourier dual space (for example curved momentum space Fourier transforms to noncommutative spacetime). By now, they also provide sources of quantum geometries \cite{BegMa} with quantum differential structures on them well studied following \cite{Wor}.   

Aside from very special classes, Hopf algebras or quantum groups in general have, however, defied classification, although there are some partial results over algebraically closed fields of characteristic zero such as $\C$, see e.g.\cite{And,Bea}. In this paper, we try a new approach which {\em does} achieve a complete classification of all Hopf algebras and bialgebras (these are like `quantum semigroups') but only of dimension $n\le 4$ and by working over the field $\F_2=\{0,1\}$ of two elements. The latter is neither  algebraically closed nor characteristic zero and we find many even to this low dimension. If one similarly classified them over $\F_3, \F_5$ etc. and looked for common families, one would get a sense of what the moduli of quantum groups looks like over a generic field. Although there is no guarantee and one would need to go to higher dimensions to be interesting, natural constructions over $\C$ (even if that is the case of interest) indeed tend to have versions over finite fields. Note that this is necessarily a hard problem as the classification of all quasitriangular Hopf algebras and their representations typically implies the classification of R-matrices or solutions of the Yang-Baxter equations which itself remains open after some decades due to its cubic-matrix form (there are partial results for example in the upper triangular case). The classification of factorisable quasitriangular Hopf algebras is in some sense dual to the classification of knots, another unsolved problem with only partial results. Therefore even the $\F_2$ case here could be seen as a significant step.

It is possible that the $\F_2$ or `digital' case could also be interesting in its own right. Indeed, this is the third in our series on `digital geometry' where previous works \cite{MaPac1,MaPac2} as well as \cite{BasMa,Ma:dem} looked at noncommutative geometry over $\F_2$. The work  \cite{MaPac2} classified `digital' quantum geometries to dimension 3 with few results for dimension $n=4$. Our complete classification of digital quantum groups for $n\le 4$ feeds into this in the same way that Lie groups and their homogeneous spaces are key examples of classical geometry. Since the axioms of a quantum group are well-known by now in mathematical physics, their digital versions could also be the first to have applications in other contexts. Whether or not they are actually useful in signal processing or electrical engineering remains to be seen but digital quantum groups allow the transfer of ideas from group theory, Fourier theory, topological quantum computing \cite{Kit} as well as ZX-calculus in general quantum computing \cite{CoDun} over to digital algorithms and digital electronics. Digital quantum computers, for example, while lacking the exponential increase of actual quantum computers, could be built now and provide training wheels and experience towards quantum computing over $\C$. 

Our results at the most rudimentary level can be summarised as
\[\begin{tabular}{|c|c|c|c|c|} 
\hline
$n$ & algebras & bialgebras & Hopf algebras & nontrivial quasitri. Hopf pairs \\ 
\hline 2 & 3 & 4 & 3 & 1 \\ 
 \hline 3 & 7 & 24 & 2 & 0 \\
 \hline 4 & 25 & 286 &  20 & 28\\ \hline\end{tabular}\]
The starting point here is that algebras over $\F_2$ in low dimension were already classified in our previous work \cite{MaPac1,MaPac2} but only in the commutative case for $n=4$ where there are 16 of them denoted A-P; our first step will be to include noncommutative ones and for $n=4$ we find 9 of these, NA-NI. For $n=3$, only the last of the algebras A-G is noncommutative and for $n=2$, none of the algebras is noncommutative. Next, a bialgebra is a vector space that is both an algebra and a {\em coalgebra} in a compatible way. The latter is the arrow-reversed concept of an algebra and in the finite-dimensional case corresponds to an algebra on the dual space. We will say, for example, that a {\em coalgebra is of type} B* if it is isomorphic to the dual of B from our list, and  that a {\em bialgebra is of type} (A,B*) if the bialgebra is isomorphic as an algebra to A from our list and isomorphic as a coalgebra to the dual of B from our list.  For a Hopf algebra, we also need an `antipode' map (a kind of linearised inverse) but this is uniquely determined if it exists, so not additional data. 

In fact, some of the possible algebras do not feature in any bialgebra, for example if they are simple as algebras (such as a field extension of $\F_2$ or a matrix algebra) as they won't then be able to admit a counit. Meanwhile, our results about the allowed bialgebras and Hopf algebras will be shown  by a {\em quiver diagram} on the set of algebras where each arrow A$\to$B, say, indicates a distinct bialgebra or Hopf algebra of type (A,B*). In this form, all possible bialgebras and Hopf algebras for $n=2,3,4$ are shown in (\ref{n2quiver}) in Section~\ref{secn2},  (\ref{n3quiver}) in Section~\ref{secn3} and (\ref{n4quiver}) (Hopf algebras only)  in Section~\ref{secn4} respectively.  Our results for $n=4$ bialgebras are given as an extended graph in Figure~\ref{n4big} with the number of a given type shown on the relevant edge. 

Finding these results will be a two-step process. First, we find all possible bialgebras where the algebra is one of our fixed algebras. Thus, the appendix explicitly lists all the coproducts that we find for $n=3$. The second step, in Section~\ref{secn3} for $n=3$, is to analyse which of them are isomorphic as bialgebras so as to find the distinct ones. For $n=4$, there are too many to list here but we have made the data available online \cite{github}. This time we depend entirely on computer analysis except for some simpler cases which are analysed by hand as a check. The $n=2,3$ cases provide the template as well as checks on the coding since these were analysed by hand. 

Most of the bialgebras and some of the Hopf algebras we find are new, after we identify the low-dimension version of known general constructions. None of the bialgebras for $n\le 3$ are `strictly quantum' in the sense of both noncommutative and noncocommutative  (the latter in our case means the algebra dual to the coalgebra is noncommutative). For  $n=4$ there are many, but of these a {\em unique} one which is a Hopf algebra (it is self-dual of type (NF,NF*)), see Proposition~\ref{nfnf}. As $\F_2$ is itself the smallest finite field, this is in some sense the {\em absolute smallest possible} strict quantum group. It appears to be new as a Hopf algebra and we call it $d_{sl_2}$ as it resembles a `digital analogue' of a quotient of $u_q(sl_2)$ at $q=-1$ in the conventions in \cite{AziMa}. The latter is 8-dimensional and commutative over $\F_2$,  so we have more as well as different relations. Also of note is a unique self-dual Hopf algebra in one anyonic variable $x^4=0$, which is not the obvious additive one, see Proposition~\ref{propG}.

Hopf algebra Fourier transform \cite{Ma:book} is also of interest.  In the finite-dimensional case over $\F_2$, there is always a unique `Haar integral' $\int:H\to \F_2$  and we compute this for every $n\le 4$ Hopf algebra in Section~\ref{secfou} and Appendix~\ref{AppB}. This is then used to compute the canonical Fourier transform $H\to H^*$. An innovation is to view Fourier transform as a representation of the quiver of Hopf algebras (i.e. a linear `Fourier transport' map associated to each arrow). Geometrically, one can thinks of this as a connection on the quiver, with curvature expressed as holonomies. 
 
Finally, with applications in mind, we look in Section~\ref{secqt} at which of our Hopf algebras admit quasitriangular structures and find  all of them for our Hopf algebras (the numbers above count the pairs $(H,\CR)$ where $H$ is a Hopf algebra and $\CR$ is a  quasitriangular structure or `universal R-matrix' \cite{Dri,Ma:book}). For example, the Grassmann plane admits 16 quasitriangular structures. We say which of the ones we find are involutive (`triangular') and which are factorisable \cite{Dri,Ma:book}. For example, $d_{sl_2}$ is triangular. In fact, we find relatively few that are non-involutive (strictly quasitriangular), but they include 8 on the Grassmann plane Hopf algebra as well as the expected factorisable ones for the Drinfeld double of $\F_2(\Z_2)$.  We conclude in Section~\ref{secconc} with some directions for further work. 
 
 The paper starts with preliminary definitions in Section~\ref{secpre}. There are many books on quantum groups but our starting point is to write everything out in terms of structure constant tensors, as these will be solved for by computer. We used Mathematica to find all the coproducts and R to look for isomorphisms. 
 
 \section*{Acknowledgements}
A.P. would like to thank S. Zonetti for helpful discussions and suggestions on the development of the coding in R.

\section{Preliminaries} \label{secpre}

In this section, we work over any field $k$, but we will be interested in the rest of the paper in  $k=\F_2$. 
Let $x^{\mu }$ be a basis of our algebra $A$ with $x^{0}=1$ the unit
and $\mu =0,\cdots ,n-1$. We write structure constants by 
\begin{equation}
x^{\mu } x^{\nu }=V^{\mu \nu }{}_{\rho }x^{\rho },\quad V^{\mu \nu
}{}_\rho\in k. 
\end{equation}%
For a unital associative algebra we of course need 
\begin{equation}\label{Veqn}
V^{0\mu }{}_{\nu }=\delta _{\nu }^{\mu }=V^{\mu 0 }{}_{\nu },\quad V^{\rho \nu }{}_{\lambda
}V^{\lambda \mu }{}_{\gamma }=V^{\nu \mu }{}_{\lambda }V^{\rho \lambda
}{}_{\gamma }. 
\end{equation}
If we do not assume that $x^0=1$ then more generally we assume $1=\eta_\mu x^\mu$ and then the unity axiom is
\[ \eta_\nu V^{\nu\mu}{}_\rho = \delta^\mu_\rho=\eta_\nu V^{\mu\nu}{}_\rho. \]

If $A$ admits the bialgebra structure, then we express the coproduct in terms of structure constants as,
\begin{equation}\label{C_coef}
\Delta x^{\mu }=C^{\mu
}{}{}_{\nu \rho }x^{\nu }\otimes x^{\rho },\quad C^\mu{}_{\nu
\rho}\in k,\quad \epsilon \left( x^{\mu
}\right) =\epsilon ^{\mu }\in k.
\end{equation}
 
For the Hopf algebra one also has the antipode
\[ S x^{\mu} =s^{\mu }{}_{\nu }x^{\nu },\quad s^{\mu }{}_{\nu }\in k.\]

On the unit of the algebra we require: $\Delta 1=1\otimes 1$ and the counit $\epsilon 1 =1\,$.
If the Hopf algebra structure exists then the antipode obeys $S1 =1$ on the unit of $A$.

The co-associativity and counity axioms for the coalgebra structure are
\begin{equation}
(\Delta \otimes \id)\circ \Delta =(\id\otimes \Delta )\circ \Delta, \qquad
(\epsilon \otimes \id)\circ \Delta =\id=(\id\otimes \epsilon )\circ \Delta, 
\end{equation}
which in tensor terms become
\begin{equation}\label{CoAs}  C^{\mu }{}{}_{\nu \gamma }C^{\nu }{}{}_{\alpha \beta } =C^{\mu
}{}{}_{\alpha \rho }C^{\rho }{}{}_{\beta \gamma },  \quad   
C^{\mu }{}{}_{\nu \rho }\epsilon ^{\nu } =\delta^{\mu }{}_{\rho }, \quad C^{\mu }{}{}_{\nu \rho }\epsilon ^{\rho }=\delta^{\mu }{}_{\nu }. 
\end{equation}
This is because we need equality of 
\[ (\Delta \otimes \id)\circ \Delta x^{\mu }=C^{\mu }{}{}_{\nu \rho
}\Delta x^{\nu }\otimes x^{\rho }=C^{\mu }{}{}_{\nu \rho }C^{\nu
}{}{}_{\lambda \gamma }x^{\lambda }\otimes x^{\gamma }\otimes x^{\rho }=
C^{\mu }{}{}_{\nu \gamma }C^{\nu }{}{}_{\alpha \beta }x^{\alpha }\otimes
x^{\beta }\otimes x^{\gamma }\]
\[ (\id\otimes \Delta )\circ \Delta x^{\mu }=C^{\mu }{}{}_{\nu \rho }x^{\nu
}\otimes \Delta x^{\rho }=C^{\mu }{}{}_{\nu \rho }C^{\rho }{}{}_{\alpha
_{1}\beta _{1}}x^{\nu }\otimes x^{\alpha _{1}}\otimes x^{\beta _{1}}=
C^{\mu }{}{}_{\alpha \rho }C^{\rho }{}{}_{\beta \gamma }x^{\alpha }\otimes
x^{\beta }\otimes x^{\gamma }\]
 and similarly for the counit. Henceforth, we leave similar such derivations 
below to the reader.

Additionally, $\Delta $ and $\epsilon $ are algebra homomorphisms which in terms of the structure constants is: 
\begin{gather} \label{CoHo}
V^{\mu \nu }{}_{\rho }C^{\rho }{}{}_{\lambda \gamma }=C^{\mu }{}{}_{\alpha
\beta }C^{\nu }{}{}_{\rho \delta }V^{\alpha \rho }{}_{\lambda }V^{\beta
\delta }{}_{\gamma }, \quad
V^{\mu \nu }{}_{\rho }\epsilon ^{\rho }=\epsilon ^{\mu }\epsilon ^{\nu }. 
\end{gather}
The antipode map obeys $
m\circ (S\otimes \id)\circ \Delta =\eta \circ \epsilon =m\circ \(id\otimes
S)\circ \Delta$ which in tensor terms is 
\begin{equation}Sx^\mu=s^\mu{}_\nu x^\nu,\quad C^{\mu }{}{}_{\nu \rho }s^{\nu }{}_{\alpha }V^{\alpha \rho }{}_{\beta }=\epsilon^\mu\delta_{\beta 0}
=C^{\mu }{}{}_{\nu \rho }s^{\rho }{}_{\lambda }V^{\nu \lambda }{}_{\beta }.
\label{CopS1}
\end{equation}

The antipode $S$ of a Hopf algebra $H$ is necessarily an algebra anti-homomorphism and a
coalgebra anti-homomorphism of $H$, hence will obey
\begin{equation}  \label{SHo}
V^{\mu \nu }{}_{\rho }s^{\rho }{}_{\lambda } =s^{\nu }{}_{\alpha }s^{\mu
}{}_{\beta }V^{\alpha \beta }{}_{\lambda },\quad s^0{}_\mu=\delta^0{}_\mu,\quad 
s^\mu{}_\nu C^\nu{}_{\alpha\beta}= C^\mu{}_{\tau\eta}s^\tau{}_\beta s^\eta{}_\alpha,\quad  s^\mu{}_\nu \eps^\nu=\eps^\mu.  \end{equation}

An algebra homomorphism $\phi(x^\mu)=\phi^\mu{}_\nu x{}^\nu$ from an algebra with product $V$ to one with product $V'$ means
\begin{equation}\label{AHo} V\phi=(\phi\tens\phi) V',\quad V^{\mu\nu}{}_\rho \phi^\rho{}_\tau=\phi^\mu{}_\alpha\phi^\nu{}_\beta V'{}^{\alpha\beta}{}_\tau.\end{equation}
and we also demand that $\eta_\mu\phi^\mu{}_\nu=\eta'_\nu$ for the units (if both algebras are in standard form then this is  $\phi^0{}_\nu=\delta^0{}_\nu$). If $\phi$ is surjective  (such as an isomorphism) then this unit condition is automatic.  
Similarly, a coalgebra homomorphism  $\psi(x^\mu)=\psi^\mu{}_\nu x^\nu$ from a coalgebra with coproduct $C'$ to one with coproduct $C$ means 
\begin{equation}\label{CHo}
C' (\psi\tens\psi)= \psi C,\quad  C'{}^\tau{}_{\alpha\beta} \psi^\alpha{}_\mu\psi^\beta{}_\nu=\psi^\tau{}_\rho C^\rho{}_{\mu\nu}.
\end{equation}
and we also demand that $\psi^\mu{}_\nu\eps{}^\nu=\eps'{}^\mu$ for the counit, which is automatic in the case of $\psi$ injective (such as an isomorphism).  Finally,  we note that the data for a coalgebra is exactly the same as the data for an algebra by
\[ C^\mu{}_{\nu\rho}\leftrightarrow V^{\nu\rho}{}_\mu,\quad \eps^\mu \leftrightarrow \eta_\mu\]
and an algebra homomorphism $\phi$ as above is equivalent to a corresponding coalgebra homomorphism  $\psi^\mu{}_\nu\leftrightarrow \phi^\nu{}_\mu$. Under this switch in interpretation,  a bialgebra is swapped  to the dual bialgebra on the dual basis and likewise in the Hopf algebra case. Another fact is that a bialgebra map between Hopf algebras is automatically a Hopf algebra map, in that it connects the antipodes, so we do not have to consider this additionally.    

We also mention some Hopf algebras notation \cite{Ma:book}. An element $x$ is called {\em primitive} if it has an `additive' coproduct $\Delta x=x\tens 1+1\tens x$ and {\em grouplike} if $\Delta x=x\tens x$. Also, if $H$ is a bialgebra then $H^{op}$ (with reversed product) and $H^{cop}$ (with reversed coproduct) are again bialgebras. The same applies for Hopf algebras if the antipode is invertible, which is always true for finite-dimensional Hopf algebras. We leave the details for Fourier theory and quasitriangular structures to their later sections, but note that the above duality swaps a quasitriangular structure with a coquasitriangular one, so we do not have to classify these separately.   

\section{Algebras and bialgebras of dimension $n\le 3$}\label{secle3}

Here, $n=2$ will be done entirely by hand, while $n=3$ will be done by manual analysis of raw data in the Appendix  (which lists all the possible coproducts as generated) obtained with the help of Mathematica. This provides the methodology for the much harder $n=4$ case in Section~\ref{secn4}. 

\subsection{Quiver for $n=2$}\label{secn2}

Here, there are just three distinct unital algebras with basis $1,x$ as in \cite{MaPac1} and looking for coproducts, one finds: 

A:\quad $x^2=0$ i.e., $\F_2\Z_2=\F_2[x]/\<x^2\>$. This is also the unital algebra with all other products zero.

{\em 2 Hopf algebras and no further bialgebras}
\begin{itemize}
\item Hopf algebra of type (A,B*) with a coalgebra A.1:  $\Delta x=x\tens 1+ 1\tens x+ x\tens x$ with $\eps x=0$ and $S x=x$ (this is $\F_2\Z_2$, $s=1+x$ is grouplike and obeys $s^2=1$). 
\item Self-dual Hopf algebra of type (A,A*) with coalgebra A.2: $\Delta x=x\tens 1+ 1\tens x$ with $\eps x=0$ and $S x=x$ (we call this the `Grassmann line').
\end{itemize}

B:\quad $x^2=x$ i.e., $\F_2(\Z_2)=\F_2[x]/\<x^2+x\>$. 

{\em 1 Hopf algebra and 2 further bialgebras}
\begin{itemize}
\item Hopf algebra of type (B,A*) with coalgebra B.1: $\Delta x=x\tens 1+1\tens x$ with $\eps x=0$ and $S x=x$ (this is $\F_2(\Z_2)$, $x=\delta_1$ the delta-function at 1).
\item Bialgebra of type (B,B*) with coalgebra B.2:  $\Delta x=x\tens x$ with $\eps x=1$  (we call this the `projector bialgebra').
\item Bialgebra  of type (B,B*) with coalgebra B.3: $\Delta x=x\tens 1+1\tens x+x\tens x$ with $\eps x=0$. \end{itemize}

In fact B has an algebra isomorphism $\phi(x)=1+x$ and $(\phi\tens\phi)\Delta_{B.3}=\Delta_{B.2}\phi$ (while B.1 is invariant) with the result that B.2$\cong$B.3 as bialgebras. Hence, up to equivalence, we have only one distinct Hopf algebra and one distinct self dual bialgebra. 

C:\quad $x^2=x+1$ (this is $\F_4=\F_2[x]/\<x^2+x+1\>$ as an algebra over $\F_2$).\\
{\em No bialgebras}

Altogether, we find 3 distinct Hopf algebras and 1 further bialgebra. The Hopf algebras are the group algebra and function Hopf algebra on $\Z_2$, which are clearly dual to each other, and the Grassmann line which arises from the observation that $\F_2[x]/\<a(x)\>$, where $a$ contains only terms of degrees that are powers of $2$, is always a Hopf algebra over $\F_2$ with primitive $\Delta x=x\tens 1+1\tens x$. The group function algebra $\F_2(\Z_2)=A_1$ is also in this family, as the smallest member of a series of Hopf algebras $A_d$ introduced in \cite{BasMa} over any characteristic $p>0$. 

We can represent our results as a quiver by drawing an arrow for each bialgebra of Hopf algebra according to its type. For example,  ${\rm A}\to {\rm B}$ since we have a Hopf algebra $\F_2\Z_2$ with algebra A and coalgebra isomorphic to the dual of B, i.e. of type (A,B*) in our notation. The dual bialgebra then implies an arrow the other way, in our case the Hopf algebra $\F_2(\Z_2)$  as the arrow B$\to$A. Altogether, we have four arrows for $n=2$:
\begin{equation}\label{n2quiver}\includegraphics[scale=.76]{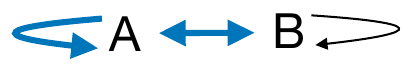}  \end{equation}
with the 3 thicker ones signalling Hopf algebras. The self-arrow on the left is the Grassmann line of type (A,A*). This is self-dual with the pairing
\[ \<1,1\>=\<x,x\>=1,\quad \<x,1\>=\<1,x\>=0,\]
which one can check obeys the Hopf algebra pairing axioms with itself. The self-arrow on the right is the projector bialgebra of type (B,B*). If we take it in the form B.2, say, then this is dually paired as a bialgebra with itself by
\[ \<1,1\>=\<1,x\>=\<x,1\>=1,\quad \<x,x\>=0.\]

\subsection{Quiver for $n=3$}\label{secn3}

In 3 dimensions we have found 7 unital associative algebras over $\mathbb{F}_2$ in \cite{MaPac2}, six of them are commutative (A - F) and one is noncommutative (G). These can all be written with basis $1,x,y$ and products

${\rm A}:\quad  x^2=y^2=xy=0$ (the unital algebra with all other products zero.)

${\rm B}:\quad x^2=x,\ y^2=y,\ xy=0$ (this is the algebra of ${\mathbb{F}}_{2}({\mathbb{Z}}_{3})$ or functions on a triangle.)

${\rm C}:\quad x^2=x,\ y^2=xy=0$ (this is $\mathbb{F}_{2}[z]/\<z^{3}+z\>$ with $z=1+x+y$ or conversely $x=1+z^2$ and $y=z+z^2$.) 

${\rm D}:\quad x^2=y,\ y^2=x, xy=x+y$  (this is the group algebra $\F_2\Z_3=\F_2[z]/\<z^3+1\>$ with $z=1+x$.)

${\rm E}:\quad x^2=y,\ y^2=xy=0$  (this is $\F_2[x]/\<x^3\>$, the anyonic line.)  

${\rm F}:\quad x^2=y,\ xy=1+y,\  y^2=1+x+y$ (this is the field $\F_8=\F_2[x]/\<x^3+x^2+1\>$.) 

${\rm G}:\quad x^2=x,\ y^2=0,\ xy=y,\ yx=0$ (this is noncommutative but ${\rm G}\cong {\rm G}^{op}$ by $x\mapsto 1+x$,  $y\mapsto y$.) 

Of these we will find that B and D are the only algebras admitting a Hopf algebra structure (namely the unique one  indicated by the notation as group algebra or function algebra on a group). The algebras  B, C, D and G 
admit many bialgebras (but no further Hopf algebras) and the algebras A, E, F admit no bialgebra structures. If we make a graph as we did for $n=2$  then we can graph our results up to equivalence as  the quiver 
\begin{equation}\label{n3quiver} \includegraphics[scale=0.75]{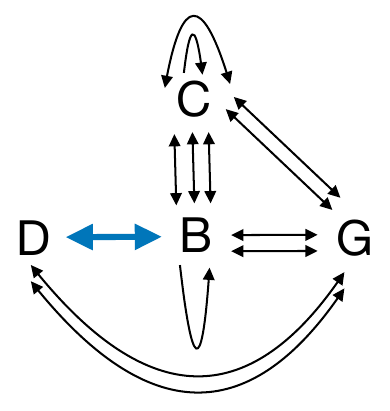} \end{equation}
where the thick arrows are Hopf algebras. Here D$\to$B is the expected group Hopf algebra $\F\Z_3$ and B$\to$D is its dual $\F(\Z_3)$ of functions on the group $\Z_3$. All the rest are strictly bialgebras and we see that two of them  are self-dual and 12 of them (those connecting to G) are noncommutative or noncocommutative (but not both). 

We will obtain this as follows. For each type of algebra, first step is to solve (\ref{Veqn}), (\ref{CoAs}), (\ref{CoHo}) for all possible bialgebra structure constants in a basis including $1$. For each solution we see if it is a Hopf algebra and we compute the dual algebra to the coalgebra to see what type it is. This `raw data' was generated using Mathematica and is collected for $n=3$ in the Appendix. We now proceed below to analyse this raw data for each algebra to determine which of the bialgebras are isomorphic, so as to have the moduli of inequivalent or distinct ones up to isomorphism. 

\subsection{Analysis for algebra B} The code produces 33 bialgebras as detailed in the Appendix:
\begin{itemize}
\item 6 bialgebras with coalgebras of type B* 
\item 18 bialgebras with coalgebras of type C*
\item 3 Hopf algebras with coalgebras of type D* 
\item 6 bialgebras with coalgebras of type G*. These are grouped as ${\rm B.8}= {\rm B.9}{}^{cop}$, ${\rm B.14}= {\rm B.15}{}^{cop}$, ${\rm B.26}= {\rm B.27}{}^{cop}$ as bialgebras with opposite coproduct by inspection. 
\end{itemize}
To narrow down to the isomorphism classes we note that B has 6 algebra automorphisms (including the identity) forming the group  $S_3$, with order 2 generators providing bialgebra isomorphisms
\[ x\mapsto 1+x+y:\quad{\rm  (i)\ B.4\cong B.21\ ;\quad (ii)\ B.1\cong B.12,\quad B.5\cong B.20,\quad B.24\cong B.30\quad B.28\cong B.31}\]
\[ {\rm B.{23}\cong B.{32},\quad B.2\cong B.17,\quad B.3\cong B.16,\quad B.7\cong B.18,\quad B.11\cong B.13} \]
\[ {\rm (iii)\ B.8\cong B.14,\quad B.9\cong B.15,\quad 
 (iv)\ B.6\cong B.19,\quad B.10\cong B.22,\quad B.25\cong B.29}
 \]
\[ x\leftrightarrow y:\quad {\rm (i)\ B.21\cong B.33;\quad (ii)\ B.1\cong B.11,\quad B.2\cong B.5,\quad B.12\cong B.24,\quad B.17\cong B.32}\]
\[ {\rm B.3\cong B.7,\quad B.16\cong B.28,\quad B.18\cong B.31,\quad B.13\cong B.30,\quad B.20\cong B.23}\]
\[{\rm  (iii) \ B.14\cong B.27,\quad B.15\cong B.26,\quad 
 (iv)\ B.6\cong B.10,\quad B.19\cong B.25,\quad B.22\cong B.29}
 \]
Now looking at the orbits under the action of the automorphism group generated by these, we see from the four cases that:

(i) The 3 bialgebras are one orbit so there is one distinct bialgebra of type (B,D*), a Hopf algebra (the group function Hopf algebra  $\F_2(\Z_3)$). 

(ii) The 18 bialgebras form three orbits and hence there are  three distinct bialgebras of type (B,C*).

(iii) The 6 bialgebras form two orbits hence there are two distinct bialgebras of type (B,G*).  By the above, one is the co-opposite of the other.

(iv) The 6 bialgebras are one orbit so there is one distinct bialgebra of type (B,B*).  This implies that B is {\em self-dual}, ${\rm B}\cong {\rm B}^*$ as bialgebras. 

\begin{proposition}\label{propBself} The coproduct {\rm B.19},  say, i.e., 
\[ \Delta x=x\tens x,\quad\Delta y=y\tens 1+1\tens y+y\tens y,\quad \eps x=1,\ \eps y=0\]
makes {\rm B} into a self-dual bialgebra with self-pairing
\[ \<1,1\>=\<1,x\>=\<x,1\>=\<x,y\>=\<y,x\>=1,\quad \<1,y\>=\<y,1\>=\<y,y\>=0. \]
\end{proposition}
\proof This is obtained as follows. The Appendix gives the algebra on the dual bases $y_0,y_1,y_2$ to $x_0=1, x_1=x, x_2=y$. Here $y_0=1+y_1$ and $y_1,y_2$ obey the relations for B. The dual of the product of B induces a coproduct on these, 
\begin{align*}\Delta y_\rho=\<x_\mu x_\nu, y_\rho\>y_\mu \tens y_\nu&=\<1,y_\rho\>y_0\tens y_0+\<x,y_\rho\>(y_0\tens y_1+y_1\tens y_0+y_1\tens y_1)\\
&\quad +\<y,y_\rho\>(y_0\tens y_2+y_2\tens y_0+y_2\tens y_2)
\end{align*}
in view of the relations $x^2=x, y^2=y$ and $xy=yx=0$. This gives 
\[ \Delta y_1=y_1\tens 1+1\tens y_1+y_1\tens y_1,\quad \Delta y_2=y_2\tens 1+1\tens y_2+ y_1\tens y_2+y_2\tens y_1+y_2\tens y_2,\]
which from the Appendix we recognise as B.10. But this is isomorphic to B.6 by the 2nd automorphism above (swapping $y_1,y_2$) and then B.6 is isomorphic to B.19 by the first isomorphism. Hence we can replace the $y_i$ by $x=1+y_1+y_2$ and $y=y_1$ and have a duality pairing of B.19 with itself. The $y_i$ being dual bases gives the pairing among the $x,y$ as shown. One can check directly that it is indeed a bialgebra pairing as it has to be by construction. \endproof

\subsection{Analysis for the algebra C} The code produces the 8 bialgebras in the Appendix:
\begin{itemize}
\item 3 bialgebras  C.2, C.7, C.8 with coalgebras of type B*
\item 3 bialgebras C.1, C.3, C.6 with coalgebras of type C* 
\item 2 bialgebras  ${\rm C4}={\rm C5}^{cop}$ by inspection, with coalgebras of type G*.
\end{itemize}

Algebra C has only the identity as an algebra automorphism, so these are all distinct. The 3 bialgebras of type (C,C*) consist of one which is self-dual and one pair related by bialgebra duality. By the same methods as in Proposition~\ref{propBself} we find that
C.6 is the self-dual coproduct on C
\[ \Delta x=x\tens x+ x\tens y+y\tens x,\quad\eps x=1,\ \eps y=0\]
with self-pairing 
\[ \< 1,1\>=\<1,x\>=\<x,1\>=\<y,y\>=1,\quad \<x,x\>=\<1,y\>=\<y,1\>=\<x,y\>=\<y,x\>=0. \]
(Use $1=y_0+y_2$, $x=y_0$, $y=y_2$ in terms of dual basis elements  for C.6 in the Appendix.) A similar calculation for the coproduct C.1 in the Appendix with $1=y_0$, $x=y_0+y_1$ and $y=y_2$ in terms the dual basis there now gives the coproduct C.3 on these, thus   (C,C.1)*=(C,C.3) as bialgebras (viewing a bialgebra as a pair consisting of an algebra and a compatible coalgebra) for the remaining two bialgebras of this type.

\subsection{Analysis for the algebra D}\label{secn3D} The code produces the 3 bialgebras in the Appendix:
\begin{itemize}
\item 1 Hopf algebra  D.1 with coalgebra of type B* (the group Hopf algebra $\F_2\Z_3$) 
\item 2 bialgebras  ${\rm D.2}= {\rm D.3}^{cop}$ by inspection with coalgebras of type G*.
\end{itemize}

The algebra D has only one nontrivial algebra automorphism $ x\leftrightarrow y $ which is not an isomorphism between any of the coalgebras. Hence these are all distinct and we have one Hopf algebra of type (D,B*) and two bialgebras of type (D,G*).

\subsection{Analysis for the noncommutative algebra G} The code produces the 8 bialgebras in the Appendix:
\begin{itemize}
\item 2 bialgebras  G.3, G.7 with coalgebras of type B*
\item 4 bialgebras  G.1, G.2, G.5, G.6 with coalgebras  of type C*
\item 2 bialgebras  G.4, G.8 with coalgebras of D*.
\end{itemize}

One can see by hand that there are {\em coalgebra} isomorphisms
 \[ x\mapsto 1+x:\quad {\rm G.1}\cong {\rm G.5},\quad {\rm G.2}\cong {\rm G.6},\quad {\rm G.3\cong G.7,\quad G.4\cong G.8}\]
which, however, reverses the product. So ${\rm G.1}\cong {\rm G.5}^{op}$, ${\rm G.2}\cong {\rm G.6}^{op}$, ${\rm G.3}\cong {\rm G.7}^{op}$, ${\rm G.4}\cong {\rm G.8}^{op}$ as bialgebras. 

The algebra G has only one nontrivial algebra automorphism which provides bialgebra isomorphisms
\[ x\mapsto x+y: \quad \ {\rm G.1\cong G.2,\quad G.5\cong G.6}\ \]
(and bialgebra automorphisms on G.3, G.7, G.4, G.8). 
As a result, there are  two distinct bialgebras each of type (G,B*), (G,C*), (G,D*), with one the opposite algebra of the other in each pair.

\section{Algebras and bialgebras of dimension $n=4$}\label{secn4}

Here it is known from \cite{MaPac1} that there are 16 unital commutative algebras and the same computer classification of algebras up to isomorphism now finds a further 9 noncommutative ones, listed in Section~\ref{secn4noncomm}. Several are known to have at least one or two commutative and cocommutative Hopf algebra structures, so part of our work will be to identify known Hopf algebras and check that all of them turn up.  Writing the basis elements $x^\mu$ explicitly as  $1,x,y,z$, Section~\ref{secn4comm} and Section~\ref{secn4noncomm} summarise all possible coproducts or `raw data' in the commutative and noncommutative cases -- they are too many to list explicitly and we refer to \cite{github} for the actual lists in machine readable form. This data plays the role of the Appendix \ref{AppA} for $n=3$.  Then, in Section~\ref{secisom}, we identify the equivalence classes, mostly by computer but with some smaller cases analysed by hand in different subsections, as a check on the coding. 

Our final result for $n=4$ in our previous quiver notation 
(where A$\to$B means a bialgebra or Hopf algebra of type (A,B*)) is at the Hopf algebra level 
\begin{equation}\label{n4quiver} \includegraphics[scale=0.73]{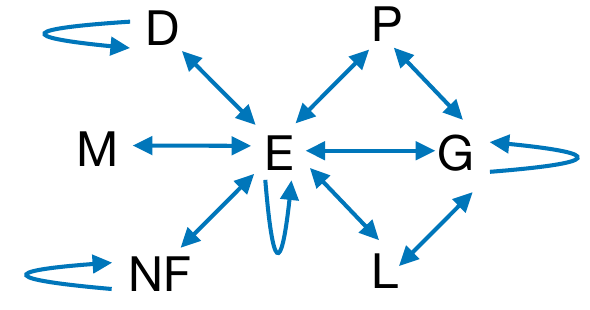}\end{equation}
The vertices here are the $n=4$ algebras below, with NF the only noncommutative one. We will come back to this diagram in Section~\ref{secisom} in a fully decorated form where we identify all the arrows as either known Hopf algebras or new ones. Also note that there is just one Hopf algebra which is both noncommutative and noncocommutative, namely the self-arrow on NF, studied in Proposition~\ref{nfnf}. The full picture for all bialgebras is also found but has too many arrows to draw as a quiver, so this is presented instead as an extended weighted graph Figure~\ref{n4big} in  Section~\ref{secisom}.

\subsection{Commutative algebras for $n=4$ and all their coproducts}\label{secn4comm} We list the commutative unital algebras, where possible, in a tensor form and/or a quartic form with relations $x^4=a x^3+ b x^2+ c x+ d$ for $a,b,c,d\in\{0,1\}$. This description is more systematic than in \cite{MaPac1} but the names of the algebras are the same.  Where helpful, we will specify the original $x,y,z$ with primes where needed, for the original description in \cite{MaPac1} for the same algebra.

A:\quad The unital algebra with all other products of $x,y,z$ zero. \\
\textit{No bialgebras}

B:\quad  All products of $x,y,z$ zero except $x^2=z$.\\
\textit{No bialgebras}

C:\quad All products of $x,y,z$ zero except $x^2=x$.\\
\textit{90 bialgebras and no Hopf algebras}:
\begin{itemize}
\item 1 bialgebra with coalgebra of type C*
\item 6 bialgebras with coalgebras of type D*
\item 3 bialgebras with coalgebras of type J*
\item 24 bialgebras with coalgebras of type K*
\item 3 bialgebras with coalgebras of type L*
\item 9 bialgebras with coalgebras of type P*
\item 6 bialgebras with coalgebras of type NC*
\item 6 bialgebras with coalgebras of type ND*
\item 2 bialgebras with coalgebras of type NE*
\item 30 bialgebras with coalgebras of type NG*
\end{itemize}

D:\quad $\F_2(\Z_2)\otimes \F_2 \Z_2\cong \F_2[w]/ \<w^4+w^2\>$ with $x^2=x$,  $y^2=0$  for the two commuting subalgebras (with $z=xy$ 
and the inherited relations $zx =z$ and $z^2=yz =0$). The quartic description is related via  $w=x+y$ and conversely $x=w^2, y=w+w^2$. One has  $x'=1+x,\ y'=y+z,\ z$  for the original description in \cite{MaPac1}.  There is a canonical Hopf algebra structure on each tensor factor as indicated by the notation.\\
\textit{4 Hopf algebras and 48 further bialgebras}:


\begin{itemize}
\item 2 bialgebras with coalgebras of type C* 
\item 2 Hopf algebras with coalgebras of type D* (includes selfdual  double $D(\F_2\Z_2)$)
\item 4 bialgebras with coalgebras of type D*  (includes selfdual proj. $\tens$ Grass.)
\item 2 Hopf algebras with coalgebras of type E* (includes $\F_2(\Z_2)\tens$ Grass. line)
\item 10 bialgebras with coalgebras of type K*
\item 8 bialgebras with coalgebras of type P* (includes projector bialgebra$\tens \F_2 \Z_2$)
\item 4 bialgebras with coalgebras of type NC*
\item 4 bialgebras with coalgebras of type ND*
\item 16 bialgebras with coalgebras of type NG*
\end{itemize}

E:\quad $\F_2\Z_2\tens\F_2\Z_2$ with all products of $x,y,z$ zero except $z=xy$ as in \cite{MaPac1}. Setting $s=x+1, t=y+1$ we have relations $s^2=t^2=1$ and $z=1+s+t+st$ and at least two Hopf algebra structures, namely either with $s,t$ grouplike and another  \cite{BasMa} is the dual of the algebra $A_2=\F_2[x]/\<x^4+x\>=L$ below.\\
\textit{76 Hopf algebras and no further bialgebras:}
\begin{itemize}
\item 24 Hopf algebras with coalgebras of type D* (includes $\F_2\Z_2\tens$ Grass. line)
\item 4 Hopf algebras with coalgebras of type E* (includes Grass. plane) 
\item  12 Hopf algebras with coalgebras of type G* (included dual of the anyonic line)
\item  12 Hopf algebras with coalgebras of type L* (includes $A_2^*$ from \cite{BasMa})
\item  8 Hopf algebras with coalgebras of type M*
\item 4 Hopf algebras with coalgebras of type P* (includes $\F_2\Z_2\tens \F_2\Z_2=\F_2\Z_2^2$)
\item  12 Hopf algebras with coalgebras of type NF*
\end{itemize}

F:\quad $\F_2[x,y]/\<y^2,x(x+y)\>$. Here all products of $x,y,z$ are zero except $x^2=z$, $z=xy$, as in \cite{MaPac1}.\\
\textit{No bialgebras}

G:\quad $\F_2\Z_4=\F_2[x]/ \<x^4\> $ is a group algebra if we take $s=1+x$ grouplike or the `anyonic line' if we take $x$ primitive. The latter is like the Grassmann line but higher order. Setting  $y=x^2$ and $z=xy=yx$ gives all other products of $x,y,z$ zero as in \cite{MaPac1}. 
\\
\textit{8 Hopf algebras and no further bialgebras:}
\begin{itemize}
\item 2 Hopf algebras with coalgebras of type E*  (includes the anyonic line)
\item 2 Hopf algebras with coalgebras of type G* 
\item 2 Hopf algebras with coalgebras of type L*
\item 2 Hopf algebras with coalgebras of type P* (includes $\F_2\Z_4$) 
\end{itemize}

H:\quad $\F_4\tens\F_2\Z_2\cong \F_2[w]/ \<w^4+w^2+1\>$ where $x^2=1+x$ for $\F_4$ and $y^2=0$ with $z=xy$.  The quartic version is related by $w=x+z$ and conversely $x=1+w^2,\ y=1+w^3$. Here $x,y,z'=y+xy$ for the original description in \cite{MaPac1}.\\ 
\textit{No bialgebras}


I:\quad $\F_2[w]/ \<w^4+w^3+w^2\>\cong\F_2[a]/\<a^4+a+1\>$ where the second version is related $a=w+1$. One also has that $x=w^2+w^3$ obeys $x^3=x^2+x$ (which implies that $x^4=x$) and together with $y=w^2$, $z=w+w^4$ fits the original description in \cite{MaPac1}.\\
\textit{4 bialgebras and no Hopf algebras:} 


\begin{itemize}
\item 1 bialgebra with coalgebra of type NC*
\item 1 bialgebras with coalgebra of type ND*
\item 2 bialgebras with coalgebras of type NG* 
\end{itemize}

J:\quad  $\F_2[y]/\<y^4+y^3\>\cong\F_2[w]/\<w^4+w^3+w^2\>$ where $x=y^2$ and $z=x(y+1)$ and the relation $yz=0$ is equivalent to the quartic for $y$,  and one can check that  $z^2=0$ as in \cite{MaPac1}.  The second quartic version is related by $w=1+y$.  \\
\textit{6 bialgebras and no Hopf algebras:}
\begin{itemize}
\item 1 bialgebra with coalgebra of type C*
\item 1 bialgebra with coalgebra of type J*
\item 2 bialgebras with coalgebras of type P* 
\item 2 bialgebras with coalgebras of type NE*
\end{itemize}

K:\quad Includes algebra $\F_2(\Z_3)$ with basis $1,x,y$ where $x^2=x$, $y^2=y$, $xy=0$, {\em plus} an additional $z$ with all products zero, as in \cite{MaPac1}.\\
\textit{96 bialgebras and no Hopf algebras:} 
\begin{itemize}
\item 8 bialgebras with coalgebras of type C*
\item 10 bialgebras with coalgebras of type D*
\item  26 bialgebras with coalgebras of type K*
\item  12 bialgebras with coalgebras of type P*
\item  7 bialgebras with coalgebras of type NC*
\item  7 bialgebras with coalgebras of type ND*
\item  4 bialgebras with coalgebras of type NE*
\item  22 bialgebras with coalgebras of type NG*
\end{itemize}

L:\quad $A_2=\F_2[x]/\<x^4+x\>$ as in \cite{BasMa} and $z=x^2$, $y=1+xz$ gives $y^2=y, z^2=x, xy=yz=0$ as in \cite{MaPac1}. There is a canonical Hopf algebra structure with $x$ primitive. \\
\textit{4 Hopf algebras and further 28 bialgebras:} 
\begin{itemize}
\item 2 bialgebras with coalgebras of type C*
\item 2 Hopf algebras with coalgebras of type E*  (includes $A_2$ from \cite{BasMa})
\item 2 Hopf algebras with coalgebras of type G*
\item 2 bialgebras with coalgebras of type L*
\item 4 bialgebras with coalgebras of type P*
\item 4 bialgebras with coalgebras of type NC*
\item 4 bialgebras with coalgebras of type ND*
\item 4 bialgebras with coalgebras of type NE*
\item 8 bialgebras with coalgebras of type NG*
\end{itemize}

M:\quad  $\F_2[z]/\<z^4+z^3+z\>\cong \F_2[w]/ \<w^4+w^2+w\>\cong \F_2[a]/\<a^4+a^3+a^2+1\>\cong \F_2[i]/ \<i^4+i^2+i+1\>$ with $x^2=1+x+y+z, \ y^2=y, \ z^2=x, \ xy=0,\ xz=1+x+y, \ yz=0$ as in \cite{MaPac1}. The isomorphism is given by $w=z^3$ and $z=w^2+w^3$. The other two quartics are just shifts to $a=z+1$ and $i=w+1$. There is a canonical Hopf algebra structure with $w$ primitive.\\
\textit{1 Hopf algebra and 2 further bialgebras:} 
\begin{itemize}
\item 1  Hopf algebra with coalgebra of type E*
\item 2 bialgebras with coalgebras of type NE*
\end{itemize}


N:\quad This is $\F_4\otimes\F_4\cong \F_2(\Z_2)\otimes \F_4$ with $x^2=1+x,\ y^2=y+1,\ z=xy$. 
 The second version has basis $1,a,b,c$ and  $a^2=a, \ c=ab, \ b^2=b+1$, related by $a=1+x+y$, $b=x$, $c=1+z$. We have  $x'=x+y+z,y'=1+x+z,z'=1+z$ we get the relations $x'^2=z',y'^2=1+x'+y+z',z'^2=x',x'y'=0,y'z'=0$ for the description in \cite{MaPac1}.
\\
\textit{No bialgebras}

0:\quad  $\F_{16}=\F_2[z]/\<z^4+z+1\>\cong \F_2[w]/\<w^4+w^3+1\>= \F_2[a]/\<a^4+a^3+a^2+a+1\>$. The second quartic version is related by $w=1+z^2+z^3$ and conversely $z=1+w+w^2$. The third is related by $a=z^3$ and conversely $z=a+a^3$, or equivalently $w=a^2+1$. There is also a different isomorphism $a=z^2+z^3$, equivalent to  $w=a+1$.   We have $x=z^2, y=1+z^3, z$  (or $z=1+x^2$ and $y=1+x+x^3$) for the description in  \cite{MaPac1}.\\
\textit{No bialgebras}


P:\quad This is $\F_2(\Z_2)\otimes \F_2(\Z_2)\cong \F_2(\Z_4)$ as algebras with $x^2=x,\ y^2=y,\ xy=yx=z$ and the  induced  relations $xz=z,yz=z,z^2=z$. The two notations indicate two canonical Hopf algebra structures, where we could also canonically identify the first version of the algebra as $\F_2(\Z_2^2)$.  Here  $x'=y+z,\ y'=x+z,\  z$ for the description in \cite{MaPac1}. \\
\textit{16 Hopf algebras and 608 further bialgebras} 
\begin{itemize}
\item  36 bialgebras with coalgebras of type C*
\item  96 bialgebras with coalgebras of type D* (includes projector bialgebra$\tens \F_2(\Z_2)$)
\item  4 Hopf algebras with coalgebras of type E* (includes $\F_2(\Z_2)\tens\F_2(\Z_2)=\F_2(\Z_2^2)$) 
\item  12 Hopf algebras with coalgebras of type G* (includes $\F_2(\Z_4)$) 
\item  24 bialgebras with coalgebras of type J*
\item  144 bialgebras with coalgebras of type K*
\item  24 bialgebras with coalgebras of type L*
\item  36 bialgebras with coalgebras of type P* (includes projector bialgebra${}^{\tens 2}$) 
\item  48 bialgebras with coalgebras of type NC*
\item 48 bialgebras with coalgebras of type ND*
\item 8 bialgebras with coalgebras of type NE*
\item  144 bialgebras with coalgebras of type NG*
\end{itemize}

In the course of the above, we have identified all 16 potential quartic algebras in dimension 4. 
Summarising them in binary $abcd$ we found only 8 distinct ones up to isomorphism, namely
\[ {\rm D}=0100,\quad {\rm G}=0000=0001,\quad  {\rm H}=0101,\quad {\rm I}=1011=1100,\]
\[{\rm J}=1000=1110,\quad {\rm L}=0010,\quad {\rm M}=1010=1101\cong 0110=0111,\quad {\rm O}=0011\cong 1001=1111 \]
where all except the two marked $\cong$ are given by the shift generator by 1 map. 

We also identified all 6 possible tensor products of the three $n=2$ algebras $\F_2\Z_2$ (the zero unital algebra for $n=2$), $\F(\Z_2), \F_4$ and 
found that only 5 of them are distinct, namely 
\[ {\rm D}=\F_2(\Z_2)\otimes\F\Z_2,\quad {\rm E}=\F_2\Z_2\otimes \F_2\Z_2, \quad {\rm H}=\F_4\otimes \F_2\Z_2,\]\[
{\rm N}=\F_4\otimes \F_4\cong\F_2(\Z_2)\tens \F_4,\quad {\rm P}=\F_2(\Z_2)\otimes \F_2(\Z_2).\]
This just leaves A,B,C,F,K of neither form. Of these, A,C,K are just an $n=3$ algebra with basis $1,x,y$ and an additional $z$ with zero products with $x,y$ and itself, while B,F are a kind of central extension with $x^2=z$. 

These remarks reassure us that our $n=4$ algebra classification  indeed turns up all the commutative algebras we might expect. If one carried out the above exercise for $n=3$ then only four of 8 possible cubics $x^3=ax^2+bx+c$ give distinct algebras with $abc$ identified as
\[ {\rm C}=010=100,\quad {\rm D}=001=110,\quad {\rm E}=000=111,\quad  {\rm F}=011=101,\]
just leaving A,B among commutative $n=3$ algebras as not cubic. There is no scope for a tensor product form  as 3 is a prime number. For $n=2$ all our algebras were quadratic, with A=00=01, B=10 and C=11. 

We then listed the maximum number of bialgebra structures for each fixed algebra, obtained by solving the equations in Section~\ref{secpre} for $C^{\mu}{}_{\nu\rho}$ and $\eps^\mu$ using Mathematica. This was done by first solving the equations for $\eps$ as these are quadratic in $\eps$ and at most linear in $\Delta$. These gave some of the variables of $\Delta$ in terms of others, and we then solved the remaining conditions for $\Delta$ on the reduced set of variables, which was then feasible in terms of computer resources. We then used R to identify the coalgebra type of each of our solutions. Note that  many of the bialgebra solutions for a fixed algebra will be isomorphic as bialgebras, which we will address in Section~\ref{secisom}. Aside from this multiplicity detail, we can, however, already see  the main structure of (\ref{n4quiver}) and Figure~\ref{n4big} from the above lists. We also identified some expected bialgebras as indicated, including all ten tensor products of the four $n=2$ bialgebras: $\F_2\Z_2$, $\F_2(\Z_2)$ (now as bialgebras), the Grassmann line and the projector bialgebra. This gives some reassurance that the coding has turned up all the $n=4$ bialgebras that we might have expected. 

\subsection{Noncommutative algebras for $n=4$ and all their coproducts}\label{secn4noncomm}
These were not considered in \cite{MaPac1} but the same method gives 9 distinct noncommutative unital algebras in dimension 4. Some of them will be cross products of commutative ones above, and for this we recall a little Hopf algebra theory. If $H$ is a bialgebra and acts on an algebra $A$ then the cross product $A\lcross H$ has $A,H$ as subalgebras and cross relations $ha=(h\o\la a) h\t$ where $\la$ is the left action (say) and $\Delta h=h\o\tens h\t$ (sum of terms) is a compact `Sweedler notation'. From our results for $n=2$, we have only $H=\F_2(\Z_2)$ and $H=\F_2\Z_2$ as Hopf algebras to consider. In fact each acts canonically on the other by `translation' $h\la a= a\o \<a\t, h\>$ where $\< \ ,\ \>$ is the duality pairing. This particular cross product is shown in \cite{Ma:book} to be the algebra of linear maps on $A$. In our case this means that
\[ M_2(\F_2)\cong \F_2(\Z_2)\lcross\F_2\Z_2\cong \F_2\Z_2\lcross\F_2(\Z_2)\]
must be one of our examples and necessarily admits no bialgebra structures. We also have an action of $\Z_2$ on $\F_4$ which will give us another example. Our 9 noncommutative algebras for $n=4$ are:

NA:\quad All products of $x,y,z$ zero except $xy=z$.  Isomorphic to its opposite algebra. \\
\textit{No bialgebras}

NB:\quad All products of $x,y,z$ zero except $x^2=z, \  xy=z,\  y^2=z$.  Isomorphic to its opposite algebra. \\
 \textit{No bialgebras}    
 
NC:\quad All products of $x,y,z$ zero except $x^2=x$ and $xy=y$. 
Opposite algebra to ND. Note that $yx=0$, so this is noncommutative. (This is the $n=3$ algebra G with $z$ adjoined with zero products.) \\
\textit{30 bialgebras and no Hopf algebras:}
\begin{itemize}
\item  2 bialgebras with coalgebras of type C*
\item  4 bialgebras with coalgebras of type D*
\item  1 bialgebras with coalgebra of type I*
\item  7 bialgebras with coalgebras of type K*
\item  2 bialgebras with coalgebras of type L*
\item  4 bialgebras with coalgebras of type P*
\item  2 bialgebras with coalgebras of type NC*
\item  2 bialgebras with coalgebras of type ND*
\item  6 bialgebras with coalgebras of type NG*
\end{itemize}

ND:\quad All products of $x,y,z$ zero except $x^2=x, \  yx=y$. 
Opposite algebra to NC. (This is  the $n=3$ algebra ${\rm G}^{\rm op}$ with $z$ adjoined with zero products.) \\
\textit{30 bialgebras and no Hopf algebras:}
\begin{itemize}
\item  2 bialgebras with coalgebras of type C*
\item  4 bialgebras with coalgebras of type D*
\item  1 bialgebra with coalgebra of type I*
\item  7 bialgebras with coalgebras of type K*
\item  2 bialgebras with coalgebras of type L*
\item  4 bialgebras with coalgebras of type P*
\item  2 bialgebras with coalgebras of type NC*
\item  2 bialgebras with coalgebras of type ND*
\item  6 bialgebras with coalgebras of type NG*
\end{itemize}

NE:\quad All products of $x,y,z$ zero except $x^2=x, \ xy=y, \  xz=z$.  Isomorphic to its opposite algebra.\\
\textit{152 bialgebras and no Hopf algebras:}
\begin{itemize}
\item 8 bialgebras with coalgebras of type C*
\item 24 bialgebras with coalgebras of type J*
\item 48 bialgebras with coalgebras of type K*
\item 24 bialgebras with coalgebras of type L*
\item 16 bialgebras with coalgebras of type M*
\item 8  bialgebras with coalgebras of type P*
\item 24 bialgebras with coalgebras of type NG*
\end{itemize}

NF:\quad All products of $x,y,z$ zero except $x^2=x, \  yx=y,\ xz=z$. Isomorphic to its opposite algebra.\\
\textit{8 Hopf algebras and no further bialgebras}:
\begin{itemize}
\item  4 Hopf algebras with coalgebras of type E*
\item  4 Hopf algebras with coalgebras of type NF*
\end{itemize}

NG:\quad All products of $x,y,z$ zero except $x^2=x,\  y^2=y, \  xz=z$. 
Isomorphic to its opposite algebra.\\ 
\textit{112 bialgebras and no Hopf algebras:}
\begin{itemize}
\item 10  bialgebras with coalgebras of type C*
\item 16 bialgebras with coalgebras of type D*
\item 2  bialgebra with coalgebras of type I*
\item 22 bialgebras with coalgebras of type K*
\item 4  bialgebras with coalgebras of type L*
\item 12  bialgebras with coalgebras of type P*
\item 6  bialgebras with coalgebras of type NC*
\item 6 bialgebras with coalgebras of type ND*
\item 2 bialgebras with coalgebras of type NE*
\item 32  bialgebras with coalgebras of type NG*
\end{itemize}

NH:\quad  $x^2=x,\ xy=0,\  yx=y,\ y^2=0,\ xz=z,\ zx=0,\  yz=1+x,\  zy=x,\ z^2=0$. Isomorphic to its opposite algebra. This is $M_2(\F_2)$ with basis 
\[ 1=\begin{pmatrix}1&0\\ 0 & 1\end{pmatrix},\quad x=\begin{pmatrix}1&0\\ 0 & 0\end{pmatrix},\quad y=\begin{pmatrix}0&0\\ 1 & 0\end{pmatrix},\quad z=\begin{pmatrix}0&1\\ 0 & 0\end{pmatrix}\]
Note that $x,y$ generate a 3-dimensional subalgebra ${\rm G}^{\rm op}$ of lower triangular matrices and $x,z$ a 3-dimensional subalgebra G of upper triangular ones, much like a finite version of $U({\rm gl}_2)$ with $1,x$ generating the Cartan subalgebra. One also has that $y+1,z+1$ generate the group algebra of $S_3$ as the group of invertible elements of $M_2(\F_2)$. For the isomorphism with the cross product, we take $\F_2(\Z_2)$ generated by $x^2=x$ and $\F_2\Z_2$ generated by $w$ say, with $w{}^2=0$ and grouplike element $s=w+1$ acting as $s\la x=1+x$. Then the cross relations are $sx=(1+x)s$ or $wx=xw+w+1$. One can check that $z=x(w+1)$ and $y=(w+1)x$ inherit the remaining relations. In matrix terms, $w$ is the matrix with all entries $1$.\\
 \textit{No bialgebras }
 
NI:\quad 
 $x^2=0, \  xy=x+z,\ yx=z,\ y^2=1+y,\  xz=0,\  zx=0,\ yz=x+z,\ zy=x,$ and $z^2=0$. 
Opposite isomorphic to itself. (This is  $\F_4\lcross\F_2\Z_2$ where $\F_2\Z_2$ with grouplike element $z=x+1$ acts on $\F_4$ generated by $y$ by the Frobenius automorphism $z\la y=y^2$. One can check that $zy=y^2z$ as equivalent to the cross relation $xy+yx=z$.) \\
 \textit{No bialgebras }

\subsection{Classification of inequivalent bialgebras for $n=4$}\label{secisom}
Unlike the $n=3$ case, most of our algebras have too many bialgebra coproducts to identify equivalence classes by hand. Therefore this was implemented by computer, using R. 
We take the full set of coalgebra solutions compatible with a fixed algebra and partition them into the different coalgebra types (as listed in the preceding subsections and obtained using Mathematica and R). We then consider all transformations $\psi\in GL(4,\F_2)$ which are coalgebra maps in the sense of \eqref{CHo} between every pair of coalgebras of the same type. Of these we ask which are algebra automorphisms (of our fixed algebra). This gives the equivalence classes of bialgebras as well as, where applicable, of Hopf algebras.

\subsubsection{Analysis for algebra {\rm G} {\bf (4 distinct Hopf algebras:  E*,P*, G*,L*)}}\label{secG}
\mbox{}\\
{\small {
\begin{tabular}{|l|l|}
\hline
Number & Hopf Algebra structure (all with $\epsilon x=\epsilon y=\epsilon z=0$) \\ \hline
\begin{tabular}{l}
G.1 \\ 
coalg. type E*\\ $\F_2[x]/\<x^4\>$
\end{tabular}
& 
\begin{tabular}{l}
$\Delta x=1\otimes x+x\otimes 1,\quad \Delta
y=1\otimes y+y\otimes 1,$ \\ 
$\Delta z=1\otimes z+x\otimes y+y\otimes x+z\otimes 1,\quad  Sx=x,\ Sy=y,\ Sz=z$
\end{tabular}
$\quad $ \\ \hline
\begin{tabular}{l}
G.2 \\ 
coalg. type G*
\end{tabular}
& 
\begin{tabular}{l}
$\Delta x=1\otimes x+x\otimes 1+y\otimes
y,\quad \Delta y=1\otimes y+y\otimes 1,$ \\ 
$\Delta z=1\otimes z+x\otimes y+y\otimes x+z\otimes 1,\quad Sx=x,\ Sy=y,\ Sz=z$
\end{tabular}
\\ \hline
\begin{tabular}{l}
G.3 \\ 
coalg. type E*
\end{tabular}
& 
\begin{tabular}{l}
$\Delta x=1\otimes x+x\otimes 1+x\otimes
y+y\otimes x+z\otimes y+y\otimes z,\quad \Delta y=1\otimes y+y\otimes 1,$ \\ 
$\Delta z=1\otimes z+x\otimes y+y\otimes x+z\otimes y+y\otimes z+z\otimes 1$,\\
$Sx=x,\ Sy=y,\ Sz=z$
\end{tabular}
\\ \hline
\begin{tabular}{l}
G.4 \\ 
coalg. type G*
\end{tabular}
& 
\begin{tabular}{l}
$\Delta x=1\otimes x+x\otimes 1+x\otimes y+y\otimes x+z\otimes y+y\otimes
z+y\otimes y$,\\
$ \Delta y=1\otimes y+y\otimes 1,\quad \Delta z=1\otimes z+x\otimes y+y\otimes x+z\otimes y+y\otimes z+z\otimes 1
$, \\ 
$Sx=x,\ Sy=y,\ Sz=z$
\end{tabular}
\\ \hline
\begin{tabular}{l}
G.5 \\ 
coalg. type P*\\ $\F_2\Z_4$
\end{tabular}
& 
\begin{tabular}{l}
$\Delta x=1\otimes x+x\otimes 1+x\otimes x,\quad \Delta y=1\otimes y+y\otimes 1+y\otimes y,$ \\ 
$\Delta z=\left( 1+x+y\right) \otimes z+z\otimes \left( 1+x+y\right)
+z\otimes z+y\otimes x+x\otimes y$, \\ 
$Sx=x+y+z,\ Sy=y,\ Sz=z$
\end{tabular}
\\ \hline
\begin{tabular}{l}
G.6 \\ 
coalg. type L*
\end{tabular}
& 
\begin{tabular}{l}
$\Delta x=1\otimes x+x\otimes 1+x\otimes x+y\otimes y+z\otimes y+y\otimes
z+z\otimes z,$ \\ 
$\Delta y=1\otimes y+y\otimes 1+y\otimes y,$ \\ 
$\Delta z=\left( 1+x+y\right) \otimes z+z\otimes \left( 1+x+y\right)
+z\otimes z+y\otimes x+x\otimes y$, \\ 
$Sx=x+y+z,\ Sy=y,\ Sz=z$
\end{tabular}
\\ \hline
\begin{tabular}{l}
G.7 \\ 
coalg. type P*
\end{tabular}
& 
\begin{tabular}{l}
$\Delta x=1\otimes x+x\otimes 1+x\otimes
x+x\otimes y+y\otimes x,\quad \Delta y=1\otimes y+y\otimes 1+y\otimes y,$ \\ 
$\Delta z=\left( 1+x\right) \otimes z+x\otimes y+y\otimes x+z\otimes
\left( 1+x\right) +z\otimes z$, \\ 
$Sx=x+y+z,\ Sy=y,\ Sz=z$
\end{tabular}
\\ \hline
\begin{tabular}{l}
G.8 \\ 
coalg. type L*
\end{tabular}
& 
\begin{tabular}{l}
$\Delta x=1\otimes x+x\otimes 1+x\otimes
x+x\otimes y+y\otimes x+y\otimes y+z\otimes y+y\otimes z+z\otimes z,$ \\ 
$\Delta y=1\otimes y+y\otimes 1+y\otimes y,$ \\ 
$\Delta z=\left( 1+x\right) \otimes z+x\otimes y+y\otimes x+z\otimes
\left( 1+x\right) +z\otimes z$, \\ 
$Sx=x+y+z,\ Sy=y,\ Sz=z$
\end{tabular}
\\ \hline
\end{tabular}
} }

The algebra G has four algebra automorphisms (including the identity) forming the group $\Z_2^2$ with order 2 generators leaving $y,z$ invariant and resulting in Hopf algebra isomorphisms
\begin{align*} &x\mapsto x+y: \quad \ {\rm G.5\cong G.7},\quad  {\rm G.6\cong G.8}\\
& x\mapsto x+z: \quad \ {\rm G.1\cong G.3}, \quad  {\rm G.2\cong G.4},\quad  {\rm G.5\cong G.7},\quad  {\rm G.6\cong G.8}.\end{align*}
Thus, the bialgebras on G are all Hopf algebras and up to isomorphism are:

(i)  The anyonic line $\F_2[x]/\<x^4\>$ Hopf algebra of type (G,E*)  with  primitive $\Delta x=x\tens 1+1\tens x$.
 
 (ii) The group Hopf algebra $\F_2\Z_4$  of type (G,P*) with $s=1+x$ grouplike, $\Delta s=s\tens s$.   
 
 (iii) A new Hopf algebra of  type (G,L*) where the coproduct G.6 and antipode in terms of the generator $x$ obeying $x^4=0$ are  
 \[ \Delta x=x\tens 1+1\tens x+x\tens x+  (x^2+x^3)\tens (x^2+x^3),\quad Sx=x+x^2+x^3.\]

 (iv) A new Hopf algebra of type (G, G*) which is self-dual. In terms of the generator $x$ obeying $x^4=0$, the coproduct G.2 and the antipode are 
 \[ \Delta x=x\tens 1+1\tens x+x^2\tens x^2,\quad Sx=x.\] 
 
 \begin{proposition}\label{propG} The self duality pairing for this Hopf algebra on the basis $1,x,y=x^2,z=x^3$ is
 \[ \<1,1\>=\<x,y\>=\<y,x\>=\<z,z\>=1\]
 and the others zero. 
 \end{proposition}
\proof Here the dual basis has $y_0=1$, the unit element of the dual,  and the other relations   $y_2^2=y_1$, $y_1y_2=y_3=y_2y_1$  and other products zero. This dual algebra to the coproduct is then isomorphic to $G$  with $y_2=x$, $y_1=y$ and $y_3=z$, giving the self-pairing stated. Finally, the coproduct on the $y_\mu$ by dualising the product of G as in Proposition~\ref{propBself} is 
\[\Delta y_1=y_1\tens 1+1\tens y_1,\  \Delta y_2=y_2\tens 1+1\tens y_2+y_1\tens y_1,\ 
 \Delta y_3=y_3\tens 1+1\tens y_3+y_1\tens y_2+y_2\tens y_1\]
 and $\eps y_1=\eps y_2=\eps y_3=0$. When written in terms of $x,y,z$, this is G.2 again. \endproof

\subsubsection{Analysis for algebra {\rm I}  {\bf (4 distinct bialgebras: NC*, ND*, 2 NG*)}}\label{secI}
\mbox{}\\
{\small {\ 
\begin{tabular}{|l|l|}
\hline
Number & Bialgebra structure (all with $\epsilon x=\epsilon y=\epsilon z=0$)\\ \hline
\begin{tabular}{l}
I.1 \\ 
coalg. type ND*
\end{tabular}
& 
\begin{tabular}{l}
$\Delta x=1\otimes x+x\otimes 1+x\otimes
x+y\otimes x,$ \\ 
$\Delta y=1\otimes y+y\otimes 1+x\otimes y+y\otimes y,$ \\ 
$\Delta z=1\otimes z+x\otimes z+y\otimes z+z\otimes 1+z\otimes x+z\otimes
y$ 
\end{tabular}
$\quad $ \\ \hline
\begin{tabular}{l}
I.2 \\ 
coalg. type NG*
\end{tabular}
& 
\begin{tabular}{l}
$\Delta x=1\otimes x+x\otimes 1+x\otimes
x+y\otimes x,$ \\ 
$\Delta y=1\otimes y+y\otimes 1+x\otimes y+y\otimes y,$ \\ 
$\Delta z=1\otimes z+x\otimes z+y\otimes z+z\otimes 1+z\otimes x+z\otimes
y+z\otimes z$ 
\end{tabular}
\\ \hline
\begin{tabular}{l}
I.3 \\ 
coalg. type NC*
\end{tabular}
& 
\begin{tabular}{l}
$\Delta x=1\otimes x+x\otimes 1+x\otimes
x+x\otimes y,$ \\ 
$\Delta y=1\otimes y+y\otimes 1+y\otimes x+y\otimes y,$ \\ 
$\Delta z=1\otimes z+x\otimes z+y\otimes z+z\otimes 1+z\otimes x+z\otimes
y$ 
\end{tabular}
\\ \hline
\begin{tabular}{l}
I.4 \\ 
coalg. type NG*
\end{tabular}
& 
\begin{tabular}{l}
$\Delta x=1\otimes x+x\otimes 1+x\otimes
x+x\otimes y,$ \\ 
$\Delta y=1\otimes y+y\otimes 1+y\otimes x+y\otimes y,$ \\ 
$\Delta z=1\otimes z+x\otimes z+y\otimes z+z\otimes 1+z\otimes x+z\otimes
y+z\otimes z$ 
\end{tabular}
\\ \hline
\end{tabular}
} }

Here ${\rm I.1}={\rm I.3}^{cop}$ and ${\rm I.2}={\rm I.4}^{cop}$ as bialgebras by inspection. 

The algebra I has only one nontrivial algebra automorphism: $x\leftrightarrow y$ which provides an isomorphism of each coalgebra with itself. Thus all four bialgebras are distinct but fall onto two pairs namely of types (I,NC*), (I,ND*) and two of type (I,NG*), where in each pair one is the co-opposite bialgebra of the other. Recall that NC is the opposite algebra to ND and NG is its own opposite algebra.

\subsubsection{Analysis for algebra {\rm J}  {\bf (5 distinct bialgebras: C*, J*, P*, 2 NE*).}}\label{secJ}
\mbox{}\\
{\small {
\begin{tabular}{|l|l|}
\hline
Number & Bialgebra structure (all with $\epsilon x=\eps y=1,\epsilon z=0)$
\\ \hline
\begin{tabular}{l}
J.1 \\ 
coalg. type P*
\end{tabular}
& 
\begin{tabular}{l}
$\Delta x=x\otimes x,\quad \Delta y=y\otimes
y,$\quad $\Delta z=x\otimes z+z\otimes x+z\otimes z$\end{tabular}
$\quad $ \\ \hline
\begin{tabular}{l}
J.2 \\ 
coalg. type P*
\end{tabular}
& 
\begin{tabular}{l}
$\Delta x=x\otimes x,\quad \Delta y=x\otimes
z+y\otimes y+y\otimes z+z\otimes x+z\otimes y,$ \\ 
$\Delta z=x\otimes z+z\otimes x+z\otimes z$\end{tabular}
\\ \hline
\begin{tabular}{l}
J.3 \\ 
coalg. type J*
\end{tabular}
& 
\begin{tabular}{l}
$\Delta x=x\otimes x+z\otimes z,$ \quad $\Delta z=x\otimes z+z\otimes x$,\\ 
$\Delta y=x\otimes x+x\otimes y+x\otimes z+y\otimes x+y\otimes z+z\otimes
x+z\otimes y$ 
\end{tabular}
\\ \hline
\begin{tabular}{l}
J.4 \\ 
coalg. type C*
\end{tabular}
& 
\begin{tabular}{l}
$\Delta x=x\otimes x+z\otimes z,$ \quad $\Delta z=x\otimes z+z\otimes x$,\\ 
$\Delta y=x\otimes x+x\otimes y+x\otimes z+y\otimes x+y\otimes z+z\otimes
x+z\otimes y+z\otimes z$
\end{tabular}
\\ \hline
\begin{tabular}{l}
J.5 \\ 
coalg. type NE*
\end{tabular}
& 
\begin{tabular}{l}
$\Delta x=x\otimes x+z\otimes 1+z\otimes x,$\quad $\Delta z=x\otimes z+z\otimes 1+z\otimes z$,
\\ 
$\Delta y=x\otimes 1+x\otimes y+y\otimes 1+z\otimes 1+z\otimes y$ \\ 
\end{tabular}
\\ \hline
\begin{tabular}{l}
J.6 \\ 
coalg. type NE*
\end{tabular}
& 
\begin{tabular}{l}
$\Delta x=1\otimes z+x\otimes x+x\otimes z,$
\quad $\Delta y=1\otimes \left( x+y+z\right) +y\otimes x+y\otimes z,$ \\ 
$\Delta z=1\otimes z+z\otimes x+z\otimes z$
\end{tabular}
\\ \hline
\end{tabular}
}}

Here ${\rm J.5}={\rm J.6}^{cop}$ as bialgebras by inspection. 

The algebra J has only one nontrivial algebra automorphism:  $y\rightarrow y+z$ (leaving $x,z$ invariant). This provides a bialgebra isomorphism  
$${\rm J.1\cong J.2}.$$ Thus up to isomorphism there is only one bialgebra of type (J,P*), one bialgebra of type (J,C*), two distinct  bialgebras of type (J,NE*)  (one is the  co-opposite bialgebra of the other), and one self-dual bialgebra of  type (J,J*). 

\subsubsection{Analysis for algebra {\rm M} {\bf (1 distinct Hopf algebra E* and 2 distinct bialgebras NE*).}}\label{secM}
\mbox{}\\
{\small {
\begin{tabular}{|l|l|}
\hline
Number & Bialgebra and Hopf algebra structure (all with $\epsilon x=0,\epsilon y=1,\epsilon z=0$) \\ \hline
\begin{tabular}{l}
M.1 \\ 
coalg. type NE*
\end{tabular}
& 
\begin{tabular}{l}
$ \Delta x=x\otimes 1+y\otimes x,\quad \Delta
y=y\otimes y,$ \quad $\Delta z=y\otimes z+z\otimes 1$
\end{tabular}
$\quad $ \\ \hline
\begin{tabular}{l}
M.2 - Hopf algebra \\ 
coalg. type E*
\end{tabular}
& 
\begin{tabular}{l}
$\Delta x=1\otimes z+x\otimes y+x\otimes z+y\otimes x+y\otimes z+z\otimes
1+z\otimes x+z\otimes y,$ \\ 
$\Delta y=1\otimes \left( 1+x+y+z\right) +\left( x+y+z\right) \otimes
1+x\otimes y+x\otimes z+y\otimes x$\\
$\qquad\qquad+y\otimes z+z\otimes x+z\otimes y$, \\ 
$\Delta z=1\otimes x+x\otimes 1+x\otimes y+x\otimes z+y\otimes x+y\otimes
z+z\otimes x+z\otimes y$, \\ 
$Sx=x,\ Sy=y,\ Sz=z$
\end{tabular}
\\ \hline
\begin{tabular}{l}
M.3 \\ 
coalg. type NE*
\end{tabular}
& 
\begin{tabular}{l}
$ \Delta x=1\otimes x+x\otimes y,\quad \Delta
y=y\otimes y,$ \quad $\Delta z=1\otimes z+z\otimes y$
\end{tabular}
\\ \hline
\end{tabular}
}}

Here ${\rm M.1}={\rm M.3}^{cop}$ as bialgebras by inspection. 

Algebra M has 3 algebra automorphisms (including the identity) forming the group $\Z_3$ generated by: $x\rightarrow z,\quad z\rightarrow 1+x+y+z$ (leaving $y$ invariant). This provides bialgebra automorphisms of M.1, M.2, M.3 separately, so these remain distinct. Thus we have two distinct bialgebras of type (M,NE*)  (one is the co-opposite of the other) and one distinct Hopf algebra of type (M,E*). The latter, if we take  M in the alternative form $\F_2[w]/\<w^4+w^2+w\>$ with $w=z^3$,  is the canonical primitive coproduct  $\Delta w=w\tens 1+1\tens w$ (as well as $\eps w=0$ and antipode $Sw=w$) which exists because all powers of $w$ in the relations are powers of 2. 
 
\subsubsection{Analysis for algebra {\rm NF}  {\bf (2 distinct Hopf algebras:  E*, NF*).}}\label{secNF}
\mbox{}\\
{\small {
\begin{tabular}{|l|l|}
\hline
Number & Hopf algebra structure (all with $\epsilon y=\epsilon z=0$)\\ \hline
\begin{tabular}{l}
NF.1 \\ 
coalg. type E*\\ $c[B_+]^*$
\end{tabular}
& 
\begin{tabular}{l}
$ \Delta x=x\otimes 1+1\otimes x,\quad \Delta y=1\otimes y+x\otimes y+x\otimes z+y\tens 1 +y\otimes x+z\otimes x,$
\\ 
$\Delta z=1\otimes z+x\otimes y+x\otimes z+y\otimes x+z\otimes 1+z\otimes
x,\quad\eps x=0$,\\
$Sx=x,\ S y=z,\ Sz=y$
\end{tabular}
$\quad $ \\ \hline
\begin{tabular}{l}
NF.2 \\ 
coalg. type NF*\\ $d_{sl_2}$
\end{tabular}
& 
\begin{tabular}{l}
$ \Delta x=1\otimes x+x\otimes 1+y\otimes
x+z\otimes x,\quad $ \\ 
$\Delta y=1\otimes y+x\otimes y+x\otimes z+y\otimes 1+y\otimes x+y\otimes
y+y\otimes z+z\otimes x,$ \\ 
$\Delta z=1\otimes z+x\otimes y+x\otimes z+y\otimes x+z\otimes 1+z\otimes
x+z\otimes y+z\otimes z,\quad \eps x=0$\\
$Sx=x+y,\ Sy=z,\ Sz=y$
\end{tabular}
\\ \hline
\begin{tabular}{l}
NF.3 \\ 
coalg. type NF*
\end{tabular}
& 
\begin{tabular}{l}
$ \Delta x=1\otimes x+x\otimes 1+x\otimes
y+x\otimes z,$ \\ 
$\Delta y=1\otimes y+x\otimes y+x\otimes z+y\otimes 1+y\otimes x +y\otimes
y+z\otimes x+z\otimes y,$ \\ 
$\Delta z=1\otimes z+x\otimes y+x\otimes z+y\otimes x+y\otimes z+z\otimes
1+z\otimes x+z\otimes z,\quad\eps x=0,$ \\ 
$Sx=x+z,\ Sy=z,\ Sz=y$
\end{tabular}
\\ \hline
\begin{tabular}{l}
NF.4 \\ 
coalg. type E*
\end{tabular}
& 
\begin{tabular}{l}
$ \Delta x=1\otimes x+x\otimes 1+xy\otimes
+x\otimes z+y\otimes x+y\otimes z+z\otimes x+z\otimes y,\quad $ \\ 
$\Delta y=1\otimes y+x\otimes y+x\otimes z+y\otimes 1+y\otimes x+y\otimes
z+zx\otimes +z\otimes y,$ \\ 
$\Delta z=1\otimes z+x\otimes y+x\otimes z+y\otimes x+y\otimes z+z\otimes
1+z\otimes x+z\otimes y,\quad\eps x=0,$ \\ 
$Sx=x+y+z,\ Sy=z,\ Sz=y$
\end{tabular}
\\ \hline
\begin{tabular}{l}
NF.5 \\ 
coalg. type E*
\end{tabular}
& 
\begin{tabular}{l}
$ \Delta x=1\otimes 1+1\otimes x+x\otimes 1,\quad $
\\ 
$\Delta y=1\otimes z+x\otimes y+x\otimes z+y\otimes x+z\otimes 1+z\otimes x,$
\\ 
$\Delta z=1\otimes y+x\otimes y+x\otimes z+y\otimes 1+y\otimes x+z\otimes
x,\quad \epsilon x=1,$ \\ 
$Sx=x,\ Sy=z,\ Sz=y$
\end{tabular}
\\ \hline
\begin{tabular}{l}
NF.6 \\ 
coalg. type NF*
\end{tabular}
& 
\begin{tabular}{l}
$ \Delta x=1\otimes 1+1\otimes x+x\otimes
1+y\otimes 1+y\otimes x+z\otimes 1+z\otimes x,\quad $ \\ 
$\Delta y=1\otimes z+x\otimes y+x\otimes z+y\otimes x+y\otimes y+y\otimes
z+z\otimes 1+z\otimes x,$ \\ 
$\Delta z=1\otimes y+x\otimes y+x\otimes z+y\otimes 1+y\otimes x+z\otimes
x+z\otimes y+z\otimes z,\quad \epsilon x=1,$ \\ 
$Sx=x+z,\quad Sy=z,\quad Sz=y$
\end{tabular}
\\ \hline
\begin{tabular}{l}
NF.7 \\ 
coalg. type NF*
\end{tabular}
& 
\begin{tabular}{l}
$ \Delta x=1\otimes 1+1\otimes x+1\otimes
y+1\otimes z+x\otimes 1+x\otimes y+x\otimes z,\quad $ \\ 
$\Delta y=1\otimes z+x\otimes y+x\otimes z+y\otimes x+y\otimes y+z\otimes
1+z\otimes x+z\otimes y,$ \\ 
$\Delta z=1\otimes y+x\otimes y+x\otimes z+y\otimes 1+y\otimes x+y\otimes
z+z\otimes x+z\otimes z,\quad \epsilon x=1,$ \\ 
$Sx=x+y,\  Sy=z,\  Sz=y$
\end{tabular}
\\ \hline
\begin{tabular}{l}
NF.8 \\ 
coalg. type E*
\end{tabular}
& 
\begin{tabular}{l}
$\Delta x=1\otimes 1+1\otimes x+1\otimes y+1\otimes z+x\otimes 1+x\otimes
y+x\otimes z+y\otimes 1$\\
$\qquad\qquad+y\otimes x+y\otimes z+z\otimes 1+z\otimes x+z\otimes
y,$ \\ 
$\Delta y=1\otimes z+x\otimes y+x\otimes z+y\otimes x+y\otimes z+z\otimes
1+z\otimes x+z\otimes y,$ \\ 
$\Delta z=1\otimes y+x\otimes y+xz\otimes +y\otimes 1+y\otimes x+y\otimes
z+z\otimes x+z\otimes y,\quad \epsilon x=1,$ \\ 
$Sx=x+y+z,\ Sy=z,\ Sz=y$
\end{tabular}
\\ \hline
\end{tabular}
}}

Here  ${\rm NF.2}={\rm NF.3}^{cop}$ and  ${\rm NF.6}={\rm NF.7}^{cop}$ as bialgebras by inspection. 

The algebra NF has 8 algebra automorphisms (including the identity) forming the group $D_4$ with order 2 and order 4 generators 
providing the following Hopf algebra isomorphisms:
\[x\rightarrow 1+x, y\leftrightarrow z:\quad (i)\quad {\rm NF.1\cong NF.5,\quad NF.4\cong NF.8,\qquad (ii)\quad NF.2\cong NF.6,\quad NF.3\cong NF.7},\] 
\[x\rightarrow 1+x+y, y\leftrightarrow z:\quad (i)\quad {\rm NF.1\cong NF.8,\quad NF.4\cong NF.5,\qquad (ii) \quad NF.2\cong NF.7,\quad NF.3\cong NF.6}.\] 
Thus all bialgebras on NF are Hopf algebras and we have up to isomorphism:

(i) A new Hopf algebra of type (NF,E*), which we will denote $c[B_+]^*$. If we use the new variable $w=y+z$ in place of $z$  and coproduct NF.1 then 
\[ x^2=x,\quad xw=y+w,\quad wx= yx=y,\quad xy=y^2=w^2=yw=wy=0\]
\[ \Delta x=x\tens 1+1\tens x,\quad \Delta w=w\tens 1+1\tens w,\quad \Delta y=y\tens 1+1\tens y+x\tens w+w\tens x \]
\[ \eps x=\eps y=\eps w=0,\quad Sx=x,\quad Sw=w,\quad Sy=y+w.\]
This has the structure of a Hopf algebra cross product \cite[Prop. 6.2.1]{Ma:book}, namely ${\rm gra}\rtimes \F_2(\Z_2)$ where the Grassmann line is generated by $w$ and $\F_2(\Z_2)$ by $x$. The latter acts by $x\la 1=0$, $x\la w=w$ so that $(1\tens x)(w\tens 1)=x\la w\tens 1+w\tens x$ or $xw=w+wx$. Here $y=w\tens x$ and the coproduct is the tensor product one.

(ii) A new Hopf algebra of type (NF,NF*). Surveying all our results, we can say this more strongly:  

\begin{proposition}\label{nfnf} There is up to isomorphism a unique noncommutative noncocommutative Hopf algebra  over $\F_2$ of dimension $n=4$, given on a basis $1,s,x,w$ by 
\[ s^2=1,\quad sx=wx=w^2=w,\quad xs=1+s+w,\quad ws=1+s+x,\quad sw=xw=x^2=x\] 
\[ \Delta s=s\tens s,\quad \Delta x=s\tens x+x\tens 1,\quad \Delta w=1\tens w+w\tens s,\quad \eps s=1,\quad \eps x=\eps w=0\]
\[ Ss=s,\quad Sx= w,\quad Sw=1+s+x.\]
This is of type {\rm (NF,NF*)} and is self-dual with duality pairing
\[ \<1,1\>=\<1,s\>=\<s,1\>=\<s,s\>=\<s,x\>=\<s,w\>= \<x,s\>=\<w,s\>=\<w,w\>=1  \]
and the others zero in this basis. We denote it $d_{sl_2}$ due to similarities with $u_q(sl_2)$. 
\end{proposition}
\proof This is the algebra NF and coproduct NF.2 in terms of new variables $w=x+y$, $s=1+y+z$ in place of $y,z$.  Letting $y_\mu$ be the dual basis to $1,s,x,w$, we first dualise the coproduct NF.2 to obtain the algebra in the dual as $y_0^2=y_0, y_2y_0=y_2, y_0y_3=y_3,y_1^2=y_1,y_1y_2=y_2,y_3y_1=y_3$ and all others zero. This is isomorphic to NF by 
\[ 1=y_0+y_1,\quad s=y_0+y_1+y_2+y_3,\quad x=y_1,\quad w=y_1+y_3.\]
This gives the self-pairing shown. Finally, we dualise the product of NF on the Hopf algebra  to a coproduct on the dual using the formula in the proof of Proposition~\ref{propBself}. Using the relations in NF to simplify and collecting terms, one arrives at
$\Delta y_\rho$. For example, $\Delta y_1=y_0\tens y_1+y_1\tens y_0+y_2\tens y_1+y_3\tens y_1$ which then reproduces the coproduct NF.2 or the stated coproduct on our identification of the $y_\rho$ with $1,s,x,w$. Hence this is a self-duality pairing of the stated Hopf algebra with itself. Note that $S^4=\id$. \endproof

\subsubsection{The $n=4$ bialgebra graph and survey of all Hopf algebras for $n\le 4$}
\mbox{}\\

\begin{figure}
\[\includegraphics[scale=0.95]{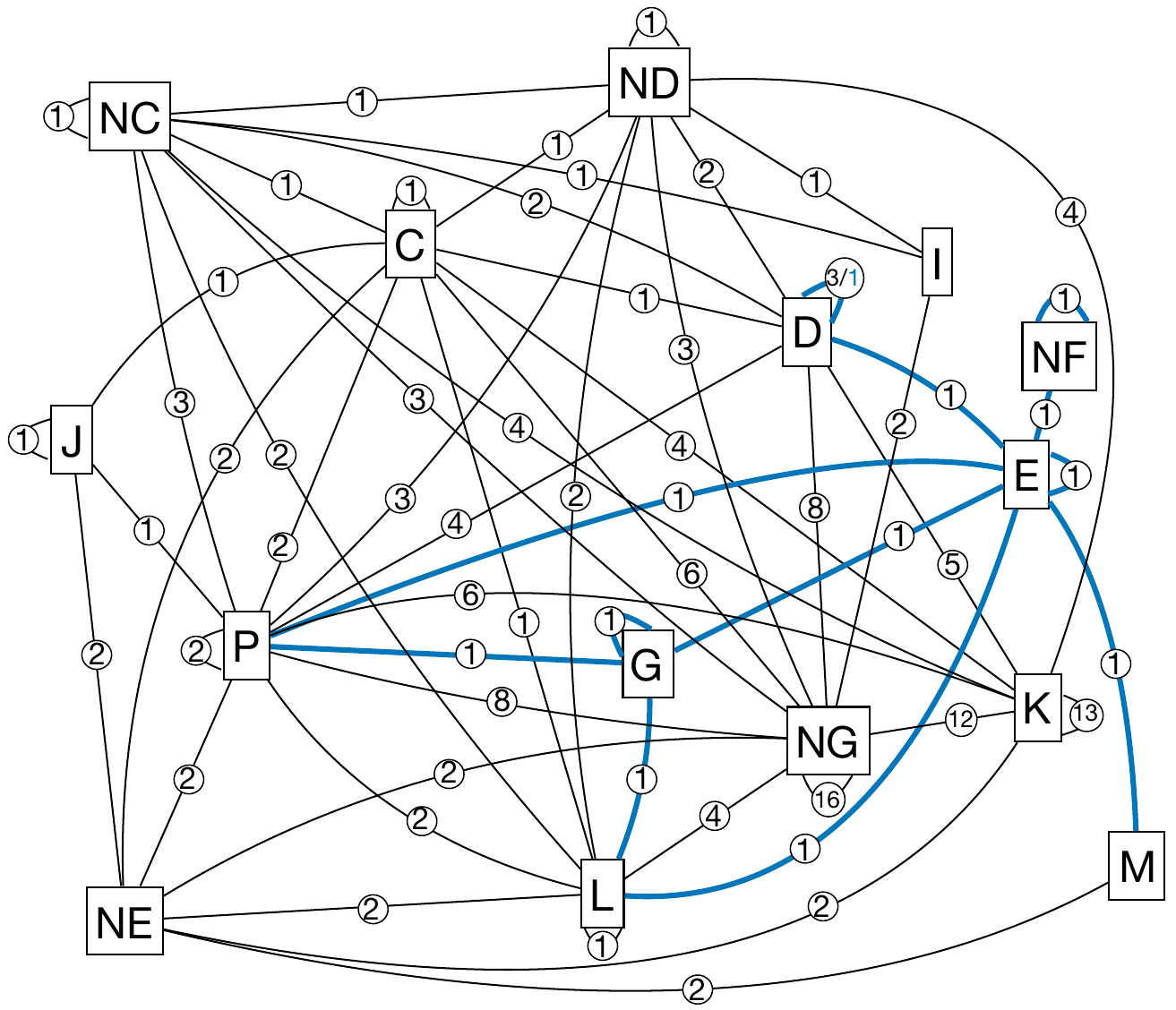}\]
\caption{Extended graph of  all $n=4$ algebras bonding with the dual of another to form bialgebras, with multiplicity (Hopf algebras in bold blue). 3/1 means 3 bialgebras of which 1 is a Hopf algebra. \label{n4big}}
\end{figure}

The remaining $n=4$ algebras admit too many bialgebras for us to list and analyse in order to find the distinct ones. Instead we summarise the resulting number of distinct bialgebras and Hopf algebras in the  weighted extended graph of  Figure~\ref{n4big}. Here  an edge means one algebra forms a bialgebra with the dual of the other algebra and the number on the edge is the number of distinct such bialgebras or Hopf algebras. More precisely, A--$i$--B means there are $i$ distinct bialgebras of type (A,B*) and $i$ of type (B, A*), while A--$i$--A means $i$ of type (A,A*). We show in bold/blue when these are actually Hopf algebras (and we show the split bialgebras/Hopf algebra multiplicity in the one case where there are some bialgebras which are not Hopf algebras). We also make the raw data available online \cite{github} from which all bialgebra coproducts can be extracted (as done in the Appendix \ref{AppA} for $n=3$) as well as lists of isomorphisms among them. 

\begin{figure}
\[\includegraphics[scale=0.77]{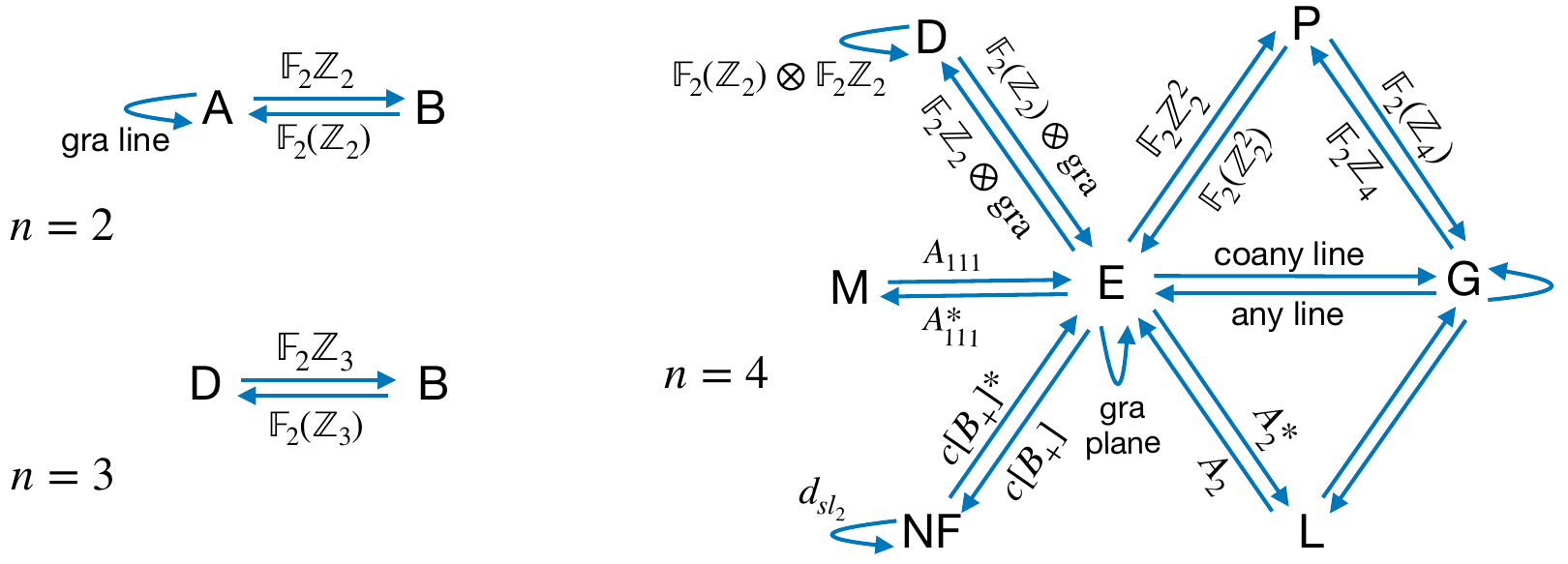}\]
\caption{Annotated quivers of all Hopf algebras for $n\le 4$.  \label{summary}}
\end{figure}

In the remainder of this section, we limit attention to the more important Hopf algebra case and make an overview of our results as shown in quiver form in Figure~\ref{summary}. The $n\le 3$ cases have already been identified as have some of the arrows in the $n=4$ diagram, notably G$\to$P is the group Hopf algebra $\F_2\Z_4$ and P$\to$G is the function Hopf algebra $\F_2(\Z_4)$ on the group $\Z_4$. In  Section~\ref{secn4comm} we took the time to identify which algebras are tensor products of our $n=2$ algebras and what should be type of coproduct when the $n=2$ algebras are given their possible bialgebra structures. Recall that the Grassmann line (i.e. one variable with $x^2=0$ and primitive coproduct $\Delta x=x\tens 1+1\tens x$) is self-dual while $\F_2\Z_2$ and $\F_2(\Z_2)$ are dual (but not isomorphic as they would be over $\C$). This accounts for 6 inequivalent tensor products as marked, where the Grassmann line tensored with itself is the two variable Grassmann algebra or Grassmann plane with $x,y$ primitive, and $\F_2\Z_2\tens\F_2\Z_2$ is the group algebra of $\Z_2^2=\Z_2\times\Z_2$. In addition, we have called G$\to$E the anyonic line:  it has one generator with $x^4=0$ and primitive coproduct. Hence E$\to$G is the dual of the anyonic line, so we denote it the coanyonic line in Figure~\ref{summary}. Taking E as for the Grassmann plane with $x^2=y^2=0$ and $xy=yx$, the different (non-primitive) coalgebra and antipode for the coanyonic line are
\begin{equation}\label{anystar} \Delta x=x\tens 1+1\tens x,\  \Delta y=y\tens 1+1\tens y+x\tens x,\quad \eps x=\eps y=0, \quad Sx=x,\ Sy=y.\end{equation}
Here the dual basis $y_\mu$ to the standard basis  of G has product dual to G.1 isomorphic to E by $y_0=1,y_1=x, y_2=y, y_3=z=xy$. The coproduct dual to the product of G was already computed in the proof of Proposition~\ref{propG}, which then gives the stated coproduct on $x,y$. 

We also expect that L$\to$E is the Hopf algebra $A_2=\F_2[x]/\<x^4+x\>$ introduced in \cite{BasMa} (it is part of a family $A_d$ defined for all $\F_p$) with its primitive coproduct $\Delta x=x\tens 1+1\tens x$, which exists because the relations only involve powers that are powers of 2. This paper also computed $A_2^*$ which would therefore be E$\to$L and is indeed built on the algebra $E$ in the alternative form $s^2=t^2=1$ with $st=ts$, justifying the identification.  

The remaining arrows in Figure~\ref{summary} are less familiar and we discuss each in turn. The arrow  M$\to E$ was described in Section~\ref{secM} as naturally built on $\F_2[w]/\<w^4+w^2+w\>$ with primitive coproduct $\Delta w=w\tens 1+1\tens w$, i.e. in the same family as $A_2$ but with different quartic relations. If we adopt a uniform notation by labels $i,j,k\in\{0,1\}$ (more compressed notation than we used before) by setting
\begin{equation}\label{Aijk} A_{ijk}=\F[w]/\<iw^4+j w^2+k w\>,\quad \Delta w=w\tens 1+1\tens w,\quad \eps w=0,\quad Sw=w\end{equation}
then $\F_2(\Z_2)=A_{011}$,  the Grassmann line is $A_{010}$, their tensor product $\F_2(\Z_2)\tens{\rm gra}$ is $A_{110}$, the anyonic line is $A_{100}$, the Hopf algebra from \cite{BasMa} is $A_2=A_{101}$ and finally, the arrow M$\to E$ is $A_{111}$, as labelled in Figure~\ref{summary}.  Its dual  
 E$\to$M  with the algebra $E$ in the above Grassmann plane form now has the non-standard (non-primitive) coproduct
\begin{equation}\label{EtoM} \begin{gathered} \Delta x=x\tens1+1\tens x+y\tens y+ xy\tens xy+ x\tens xy+ xy\tens x\\
 \Delta y=y\tens 1+1\tens y+x\tens y+y\tens x+ (x+y+xy)\tens(x+y+xy)\end{gathered}\end{equation}
along with $\eps x=\eps y=0$ and a certain antipode, as some kind of non-linear Grassmann plane. Here the dual basis $y_\mu$ to the standard basis of M  has product dual to M.2 isomorphic to E if $y_0+y_2=1$, $y_1=x+y+z$, $y_2=x$ and $y_3=y$ (say) and then the coproduct on the $y_\mu$ dual to the product of M gives the above coproducts. 

Next, the arrow NF$\to$E was described in  Section~\ref{secNF} and is noncommutative but cocommutative. This is therefore `like' a nonAbelian group algebra although not precisely. Rather, we identified it as a cross product ${\rm gra}\rtimes\F_2(\Z_2)$. Its dual E$\to$NF is therefore commutative but noncocommutative, so something like  functions on a nonAbelian group, although not precisely (as there is none of order 4). In fact it is a reduced version of an algebraic group of the Borel subgroup $B_+\subset SL_2$ and we denote it $c[B_+]$ in Figure~\ref{summary} for this reason. Using the algebra E in the mixed form with $s=1+x$ in place of $x$ (so that the relations are $s^2=1,y^2=0$ and $s,y$ commute),  the Hopf algebra structure  
\begin{equation}\label{CBplus} \Delta s=s\tens s,\quad \Delta y=y\tens 1+s\tens y,\quad \eps s=1,\quad \eps y=0,\quad Ss=s,\quad Sy=sy,\end{equation}
which is a cross coproduct of the Grassmann line and $\F_2\Z_2$. Here we take NF in the basis $1,x,y,w=y+z$ for NF.1 in Section~\ref{secNF} then its dual basis $y_\mu$ has product dual to NF.1 which can be identified with E in its standard basis by $y_0=1, y_1=x, y_2=z, y_3=y$. Then the product of NF dualises to the coproduct shown by calculations similar to those done previously. One can also think of this Hopf algebra as a quotient of the Taft algebra or the reduced quantum group $u_q(b_+)\subset u_q(sl_2)$ at $q=-1$.  

It remains to discuss the four arrows in Figure~\ref{summary} which do not appear to be part of known constructions. The self-arrow on NF was studied in detail in Proposition~\ref{nfnf} as another self-dual Hopf algebra but this time noncommutative and noncocommutative, and denoted $d_{sl_2}$.  The unmarked self-arrow on G is the self-dual Hopf algebra in Proposition~\ref{propG} built on the same algebra $\F_2[x]/\<x^4\>$ as the anyonic line but now with a certain non-standard (non-primitive) coalgebra, i.e.  non-linear version of the anyonic line. We also described in that section the arrow G$\to$L as built on this same algebra with another non-standard (non-primitive) coalgebra, i.e. another nonlinear anyonic line. Its dual L$\to G$ is built on the same algebra $\F_2[x]/\<x^4+x\>$  as $A_2$ but again with a non-standard (non-primitive) coalgebra and antipode
\begin{equation}\label{LtoG} \Delta x=x\tens 1+1\tens x+(x+x^2)\tens(x+ x^2),\quad \eps x=0,\quad Sx=x^2\end{equation}
as a nonlinear version of $A_2$. Here the dual basis $y_\mu$ to the standard basis  of G has product dual to G.6 isomorphic to L by $y_0=1, y_1+y_2=x, 1+y_1+y_2+y_3=y, y_2=z$. The coproduct dual to the product of G was already computed in the proof of Proposition~\ref{propG}, which now gives the stated coproduct on $x$.  

This completes our narrative of Figure~\ref{summary}. We also make a general observation about the figure itself. In principle, this should have been a quiver with possible multiple arrows between nodes, but in practice we see that it is in fact a directed graph extended to include self-arrows. This amounts to:

\begin{proposition} For dimension $n\le 4$, there is at most one Hopf algebra over $\F_2$ of any given algebra-coalgebra type, up to isomorphism.
\end{proposition}

\section{Integrals and Hopf algebra Fourier transforms for $n\le 4$} \label{secfou}

Integrals and Fourier transform are canonically determined up to scale for any finite-dimensional Hopf algebra, which over $\F_2$, means unique. Hence, in this section and Appendix~\ref{AppB},  we compute and list them for each of the 25  Hopf algebras for $n\le 4$ with a  fixed representative coalgebra in each isomorphism class and a each fixed each algebra. Fourier transform then  amounts to a linear map for every arrow in  Figure~\ref{summary} or quiver representation, which one can think of as a connection on the quiver \cite[Lem.~4.1]{MaTao2}. At any given algebra, we have arrows pointing to different possible coalgebra types  and if we choose one of these arrows then Fourier transform on the corresponding Hopf algebra `transports'  an element of our source algebra to an element of the target algebra. This quiver connection generically has curvature in the sense of monodromy around a loop in the graph. This is already clear for 1-step loops (on a self-dual Hopf algebra) and 2-step loops. From a cocommutative Hopf algebra to its dual and back, the composite is the antipode \cite{Ma:book,BegMa} (more generally, the canonical way back is a second `adjoint' Fourier transform or quiver representation with left and right reversed in the integrand). 

We recall that every finite-dimensional Hopf algebra $H$ has an  integral $\int:H\to k$ on $H$ characterised abstractly and then in tensor terms by 
\begin{equation}\label{rightintH}  (\int\tens\id)\Delta=1\int;\quad \int x^\mu=I^\mu,\quad   C^\mu{}_{\nu\rho}I^\nu= I^\mu \delta_{\rho,0}\end{equation}
for some $I^\mu\in k$. This is for a right-integral (one can also have left one) and over $\F_2$ it is unique as there is no scaling. We then define the Hopf algebra Fourier transform $H\to H^*$ as in \cite[Prop 1.7.7]{Ma:book}, which in our terms comes down to
\begin{equation}\label{FouH}   \CF(x^\mu)=(\int x^\nu x^\mu)y_\nu=F^{\mu\nu}y_\nu;\quad F^{\mu\nu}=V^{\nu\mu}{}_\rho I^\rho,\end{equation}
where $\{y_\mu\}$ is  the dual basis to a basis $\{x^\mu\}$ of $H$. This can be shown to be invertible \cite{BegMa}, key to which is the adjoint Fourier transform 
\begin{equation}\label{adjFouH} \CF^\#(x^\mu)=(\int x^\mu x^\nu)y_\nu=F^{\#\mu\nu}y_\nu;\quad F^{\#\mu\nu}=V^{\mu\nu}{}_\rho I^\rho\end{equation}
which differs only in the noncommutative case and obeys $\CF^\#\circ\CF\propto S$. We work over $k=\F_2$ and find $I^\mu$ and $F^{\mu\nu}$. In the spirit of earlier sections, we first do the $n\le 3$ case in detail to show the method: 

(i) For the self-dual Grassmann line, we have $x^2=0$ and basis $1,x$. Then
\[ \int 1=0,\quad \int x=1,\quad \CF_{\rm A\to A}=\begin{pmatrix}0&1\\ 1 & 0\end{pmatrix},\quad \CF_{\rm A\to A}^2=\id.\]
Here, the primitive coproduct means we need $1\int x=(\int x)1+(\int 1)x$ which fixes the integral as we don't want $\int=0$ entirely. The Fourier transform A$\to$A is then $\CF(1)=(\int x.1)y_1=x$ and $\CF(x)=(\int 1.x)y_0=1$ given the self-paring on (A,A*) in Section~\ref{secn2}. It squares correctly as the antipode is the identity map. 

(ii) For Fourier transform $\F_2\Z_2\to \F_2(\Z_2)$, we have standard form $\F_2\Z_2$ with $x^2=0$ and $\F_2(\Z_2)$ with $x^2=x$, both with denoted basis $1,x$. Then
\[ \int 1=\int x=1,\quad \CF_{\rm A\to B}=\begin{pmatrix}1&1\\ 1 & 0\end{pmatrix}.\]
Here $\Delta x=x\tens 1+1\tens x+x\tens x$ tells us that $\int 1=\int x$ so these are 1 for $\int$ on $\F_2\Z_2$. The two bases specified are already dual, $y_0=1$ and $y_1=x$ of $\F_2(\Z_2)$, so Fourier transform A$\to$ B is $\CF(1)=(\int 1.1)1+(\int x.1)x=1+x$, $\CF(x)=(\int 1.x)1=1$. Similarly, for Fourier transform $\F_2(\Z_2)\to \F_2\Z_2$: 
\[ \int 1=0,\quad \int x=1,\quad \CF_{\rm B\to A}=\begin{pmatrix}0&1\\ 1 & 1\end{pmatrix}.\]
Here the target coalgebra and hence $\int$ on $\F_2(\Z_2)$ is the same as (i) and $\CF(1)=(\int x.1)1=x$, $\CF(x)=(\int (1.x)1+(\int x.x)x=1+x$. These are mutually inverse as the antipode is the identity on these Hopf algebras. 

(iii) For Fourier transform $\F_2\Z_3\to \F_2(\Z_3)$, we have standard form $\F_2\Z_3$ with $x^2=y, y^2=x, xy=x+y=yx$ and $\F_2(\Z_3)$ with $x^2=x, y^2=y, xy=0$, both with basis $1,x,y$. Then 
\[ \int 1=\int  x=\int y=1,\quad \CF_{\rm D\to B}= \begin{pmatrix}1&1& 1\\ 1& 1& 0\\ 1 &0 & 1\end{pmatrix}. \]
Here the coproduct D.1 in Appendix~\ref{AppA} has the same form on $x$ and $y$ as in (ii) so $\int 1=\int x=\int y$, so these are all 1 for $\int$ on $\F_2\Z_3$. It is also evident from their relations in  Appendix~\ref{AppA} that the bases are dual, $y_0=1, y_1=x$, $y_2=y$, so Fourier transform D$\to$ B is $\CF(1)=(\int 1.1)1+(\int x.1)x+(\int y.1)y=1+x+y$, $\CF(x)=(\int 1.x)1+(\int x.x)x=1+x$ and $\CF(y)=1+y$ analogously. Similarly, for Fourier transform $\F_2(\Z_3)\to \F_2\Z_3$: 
\[ \int 1=\int x=\int y=1,\quad \CF_{\rm B\to D}= \begin{pmatrix}1&1& 1\\ 1& 1& 0\\ 1 &0 & 1\end{pmatrix}.\]
We use the coproduct B.4 in Appendix~\ref{AppA} where $\Delta x=x\tens 1+1\tens x+ x\tens y+y\tens x+y\tens y$ from which $\int 1=\int x=\int y$ so these are 1 for $\int$ on $\F_2(\Z_3)$. Then $\CF$ comes out with the same formula by a similar calculation. The composite is
\[ \CF_{\rm D\to B}\circ \CF_{\rm B\to D}= S,\quad \CF_{\rm B\to D}\circ \CF_{\rm D\to B}=S,\quad S=\begin{pmatrix}1&0& 0\\ 0& 0& 1\\ 0 &1 & 0\end{pmatrix}\]
being the antipode on either Hopf algebra. We also do the self-dual ones from $n=4$ by hand, as operators $\CF:H\to H$.

\begin{proposition}\label{ft} The self-dual Hopf algebras for $n=4$ have integrals and Fourier transforms: 

(i) (Grassmann plane) in the form $x^2=y^2=0$, $xy=yx$ and basis $1,x,y,xy$:
\[ \int 1=\int x=\int y=0,\quad  \int xy=1,\quad \CF=\begin{pmatrix}0&0&0&1\\ 0 & 0&1 & 0\\ 0& 1 & 0 & 0\\ 1& 0& 0& 0\end{pmatrix},\quad \CF^2=\id \]

(ii) $\F_2(\Z_2)\tens\F_2\Z_2$  in the form $x^2=x$, $y^2=0$, $xy=yx$ and basis $1,x,y,xy$:
\[ \int 1=\int y=0,\quad \int x=\int xy=1,\quad \CF=\begin{pmatrix}0&0&0&1\\ 0 & 1&0 & 0\\ 0& 0 & 1 & 0\\ 1& 0& 0& 0\end{pmatrix},\quad \CF^2=\id \]

(iii)  The self-dual version of the anyonic line in Proposition~\ref{propG} in the form $x^4=0$ and basis $1,x,x^2,x^3$:
\[  \int 1=\int x=\int x^2=0,\quad \int x^3=1,\quad  \CF=\begin{pmatrix}0&0&0&1\\ 0 & 1&0 & 0\\ 0& 0 & 1 & 0\\ 1& 0& 0& 0\end{pmatrix}, \quad \CF^2=\id \]
(iv) $d_{sl_2}$ in the form in Proposition~\ref{nfnf} with basis $1,s,x,w$: 
\[ \int 1=\int s=\int w=0,\quad \int x=1,\quad \CF=\begin{pmatrix}1&1&1&1\\ 0 & 0&1 & 1\\ 0& 1 & 0 & 1\\ 0& 0& 1& 0\end{pmatrix}, \quad \CF^3=\id.\]
\end{proposition}
\proof Here $\CF(x^\mu)=\tilde F^\mu{}_\nu x^\nu$ defines the matrices $\tilde F^\mu{}_\nu$ as displayed (they are  $F^{\mu\nu}$ composed with the chosen self-duality pairing). The order of this matrix is also shown. Note that from the quiver point of view, the matrices are $\CF_{\rm E\to E}, \CF_{\rm D\to D},$  $\CF_{\rm G\to G},  \CF_{\rm NF\to NF}$ for the respective self-arrows.  

(i) As $x,y$ are primitive, we similarly have  $\int 1=0$ from their coproducts. Then $\Delta(xy)=xy\tens 1+1\tens xy+x\tens y+y\tens x$ tells us that we also need $1\int xy=(\int xy)1+(\int x)y+(\int y)x$ forces the values shown. This is the tensor product of two Grassmann lines so the dual basis is $y_0=1, y_1=x, y_2=y, y_3=xy$. Hence the Fourier transform is $\CF(1)=(\int xy.1)y_3=xy,$ $\CF(x)=(\int y.x)y_2=y,$   $\CF(y)=(\int x.y)y_1=x,$ $\CF(xy)=(\int 1.xy)=y_0=1$.

(ii) Here $x$ is primitive so $\int 1=0$ and $\Delta y=y\tens 1+1\tens y+y\tens y$ tells us $\int y=0$. Then $\Delta(xy)= xy\tens 1+ x\tens y+xy\tens y+y\tens x+y\tens xy+1\tens xy$ tell us that $\int x=\int xy$ hence these have to be 1. One can check that the dual basis is isomorphic to D with $y_0=1,$  $y_1=y,$ $y_2=x$ and $y_3=z=xy$ (so the self-pairing is with $\<1,1\>=\<x,y\>=\<y,x\>=\<z,z\>=1$ and the others zero). Then $\CF(1)=(\int xy.1)y_3=xy$, $\CF(x)=(\int y.x)y_2=x$, $\CF(y)=(\int x.y)y_1=y$ and $\CF(xy)=(\int 1.xy)y_0=1$. 

(iii) We have $\Delta x^2=(\Delta x)^2=x^2\tens 1+1\tens x^2$ which forces $\int 1=0$. Then the coproduct in $\Delta x=x\tens 1+1\tens x+x^2\tens x^2$ needs $\int x^2=0$. We then compute $\Delta x^3=x^3\tens 1+x^2\tens x+x\tens x^2+1\tens x^3$, which tells us $\int x=0$. Then $\CF(1)=(\int x^3.1)y_3=y_3=x^3$, $\CF(x)=(\int x^2.x)y_2=y_2=x$, $\CF(x^2)=(\int x.x^2)y_1=y_1=x^2$ and $\CF(x^3)=\int(1.x^3)y_0=1$ using the identification of the dual basis in Proposition~\ref{propG}. 

(iv) We have $\int 1=\int s=0$ as $s$ is grouplike, and $\int w=0$ from the coproduct of $w$, giving the stated integral. Then $\CF(1)=(\int x.1)y_2=1+s+x+w$, $\CF(s)=(\int w.s)y_3=y_3=x+w$, $\CF(x)=(\int 1.x)y+0+(\int x.x)y_2=y_0+y_2=s+w$ and
$\CF(w)=(\int x.w)y_1=y_1=x$ where we sum over basis elements that give  product $x$ in the integrand and identify the dual bases as in Proposition~\ref{nfnf}. 
\endproof
 
The remaining 16 Hopf algebra integrals and Fourier transforms for $n=4$  are done similarly, with most of the work for the dual bases and their identification with the standard basis of the dual algebra already done in Section~\ref{secisom}. Results for these are included in the combined $n=4$ Hopf algebras tables in Appendix~\ref{AppB}. In each case the algebra is in standard form with basis $1,x,y,z$ and we choose the coproduct in each Hopf algebra isomorphism class that matches  earlier sections, notably the discussion of Figure~\ref{summary}. In particular, from the third columns there,  we have  that the composition along the cycles, 
\[ \CF_{\rm E\to P\to G\to E}=\left( 
\begin{array}{cccc}
0& 0 & 0 & 1 \\ 
0 & 0 & 1 & 0 \\ 
0 & 1 & 0 & 0 \\ 
1 & 0 & 0 & 0
\end{array}
\right),\quad  \CF_{\rm E\to G\to L\to E}=\left( 
\begin{array}{cccc}
0& 0 & 0 & 1 \\ 
0 & 1 & 1 & 1 \\ 
0 & 0 & 1 & 1 \\ 
1 & 0 & 0 & 0
\end{array}
\right),\quad \CF_{\rm E\to P\to G\to L\to E}=\left( 
\begin{array}{cccc}
1 & 0 & 0 & 0 \\ 
0 & 0 & 1 & 1 \\ 
0 & 1 & 1 & 1 \\ 
0 & 0 & 0 & 1
\end{array}
\right)\]
are all nontrivial. The matrix for the first case is the product of the matrices for $\CF_{E\to P}$, $\CF_{P\to G}$, $\CF_{G\to E}$ in that order, similarly for the others, and the matrix for the third case is the product of the first two as it must be since the antipode on the anyonic line (G,E*) is the identity (so the overlapping part cancels). The above holonomies $\rm E\to E$ have orders $2,4,3$ respectively. 

\section{Quasitriangular structures  for $n\le 4$}\label{secqt}

A Hopf algebra or bialgebra is called {\em quasitriangular} if equipped with an invertible $\CR\in H\tens H$ obeying  \cite{Dri,Ma:book}
\begin{equation}\label{Dri}(\Delta\tens\id)\CR=\CR_{13}\CR_{23},\quad (\id\tens\Delta)\CR=\CR_{13}\CR_{12},\quad \CR\Delta h= (\Delta^{cop}h)\CR\end{equation}
for all $h\in H$, where the numerical suffices indicate the position in $H^{\tens 3}$. It can be shown that $(\eps\tens\id)\CR=1=(\id\tens\eps)\CR$ and in the Hopf algebra case that  $\CR^{-1}=(S\tens\id)\CR$ and $(S\tens S)\CR=\CR$. 
In the Hopf algebra case, we can use this to define $\CR^{-1}$ as long as the $\eps$ conditions hold. It is well known that $\CR$ automatically obeys the Yang-Baxter or braid relations in the form
\[ \CR_{12}\CR_{13}\CR_{23}=\CR_{23}\CR_{13}\CR_{12}.\]

For explicit formulae in the style of Section~\ref{secpre}, one can write 
\[ \CR=R_{\mu\nu}x^\mu\tens x^\nu\]
then the axioms and $\eps$ properties are easily seen to become
\begin{equation}\label{Rpent} R_{\mu\rho}C^\mu{}_{\alpha\beta}=R_{\alpha\mu}R_{\beta\nu}V^{\mu\nu}{}_\rho,\quad R_{\mu\nu}C^\nu{}_{\alpha\beta}=R_{\rho\alpha}R_{\nu\beta}V^{\nu\rho}{}_\mu\end{equation}
\begin{equation}\label{Rcocom} R_{\mu\nu}C^\rho{}_{\alpha\beta} V^{\mu\alpha}{}_\sigma V^{\nu\beta}{}_\tau =R_{\mu\nu}C^\rho{}_{\beta\alpha}V^{\alpha\mu}{}_\sigma V^{\beta\nu}{}_\tau \end{equation}
\begin{equation}\label{Reps} \eps^\mu R_{\mu\nu}=\delta_{\nu,0}=\eps^\mu R_{\nu\mu},\end{equation}
to which one can explicitly add the existence of an inverse $\CR^{-1}=R^-_{\mu\nu}x^\mu\tens x^\nu$ in the bialgebra case, defined by its product with $\CR$.  

For each such $\CR$, there is an associated {\em quantum Killing form} $\CQ=\CR_{21}\CR$ in \cite{Dri,Ma:book} clearly given in our case by
\begin{equation}\label{Q} \CQ=Q_{\mu\nu}x^\mu\tens x^\nu;\quad Q_{\mu\nu}=R_{\alpha\beta}R_{\sigma\tau}V^{\beta\sigma}{}_\mu V^{\alpha\tau}{}_\nu.\end{equation}
Following Drinfeld, we say that  $H$ is {\em triangular}   if $\CQ=1\tens 1$, i.e.,  $Q_{\mu\nu}=\delta_{\mu,0}\delta_{\nu,0}$ for a standard basis where $1=x^0$, and {\em factorisable} if $\CQ$ is invertible as an operator $H^*\to H$, which just amounts to $Q_{\mu\nu}$ invertible.


For hand calculations, we will find it easier to work with a dual concept, a coquasitriangular structure $\CR:H\tens H\to k$ obeying \cite{Ma:book}
\begin{equation}\label{cqt1} \mathcal{R}(fg\tens h)=\mathcal{R}(f\tens h{\o})\mathcal{R}(g\tens h{\t}),\quad \mathcal{R}(f\tens gh)=\mathcal{R}(f{\o}\tens h)\mathcal{R}(f{\t}\tens g),\end{equation}
\begin{equation}\label{cqt2} g{\o}h{\o}\mathcal{R}(h{\t}\tens g{\t})=\mathcal{R}(h{\o}\tens g{\o})h{\t}g{\t}\end{equation}
for all $f,g,h\in H$,  where $\Delta h:=h\o\tens h\t$ is a standard notation (summation understood). The counit condition now appears dually as 
\begin{equation}\label{cqt3} \CR(h\tens 1)=\CR(1\tens h)=\eps h.\end{equation}
This is equivalent in the finite-dimensional case to a quasitriangular structure on $H^*$.   The quantum Killing form $\CQ:H\tens H\to k$ is then defined by 
\begin{equation}\label{cqt4}\CQ(g\tens h)=\CR(h\o\tens  g\o)\CS(g\t\tens h\t).\end{equation}

Finally, we note that for all cocommutative Hopf algebras, $\CR=1\tens 1$ is trivially a triangular structure and also that if $H$ is commutative and cocommutative then the third axiom of (\ref{Dri})  is automatic. In this case, if $H_1\to H$ is a Hopf algebra map with $H_1$ is quasitriangular then $H$ inherits a quasitriangular structure by mapping $\CR$ to $H\tens H$. This happens, in particular, if  $H$ has a triangular sub-Hopf algebra. 

We now work over $\F_2$ and by our remark, all Hopf algebras for $n\le 4$  are trivially triangular with $\CR=1\tens 1$ other than $c[B_+]$ and $d_{sl_2}$. Likewise all are trivially cotriangular with $\CR=\eps\tens\eps$ other than $c[B_+]^*$ and $d_{sl_2}$. We are only interested in nontrivial cases. As in Section~\ref{secfou},  we analyse the $n\le 3$ Hopf algebra case by hand. In fact there are very few.

\begin{proposition} \label{qt3} For $n\le 3$, only the Grassmann line  with $x^2=0$ and $x$ primitive admits a nontrivial quasitriangular structure,  $\CR=1\tens 1+ x\tens x$, which is triangular. 
\end{proposition}
\proof (i) For the Grassmann line, we have a free choice of $\CR(x\tens x)$ so this can be 0 or 1 for a coquasitriangular structure. This translates to the trivial option and the non-trivial one shown, given the self-duality.

(ii) $\F_2\Z_2$ in the form $x^2=0$ has $\Delta x=x\tens 1+1\tens x+x\tens x$, we have $0=\CR(x^2\tens x)=\CR(x\tens x)^2$ hence $\CR$ is trivial. Equivalently, $\F_2(\Z_2)$ admits on $\CR=1\tens 1$ for the quasitriangular structure.  On the dual side, for $\F_2(\Z_2)$ in the form $x^2=x$ with $x$ primitive, we have $\CR(x\tens x)=\CR(x^2\tens x)=2\CR(x\tens x)\CR(x\tens 1)=0$ so again only the trivial $\CR=1\tens 1$ on $\F_2\Z_2$. 

(iii)  $\F_2\Z_3$ in the standard form has the same coproduct on $x$ as in (ii) but now the relation $x^2=y$ tells us that $\CR(y\tens x)=\CR(x\tens x)^2=\CR(x\tens x)$ and similarly $\CR(x\tens y)=\CR(x\tens x)$ on the other side. But by symmetry in $x,y$, these are also $\CR(y\tens y)$. Then considering $xy=x+y$ we have $0=\CR(xy\tens s)=\CR(x\tens x)\CR(y\tens x)$, so these are all zero and $\F_2(\Z_3)$ has only the trivial quasitriangular structure. On the other side with $\F_2(\Z_3)$ in standard form and coproduct B.4 in the Appendix, we have $\CR(x\tens x)=\CR(x^2\tens x)=\CR(x\tens y)^2=\CR(x\tens y)$ and similarly $=\CR(y\tens x)$ on the other side. By symmetry in $x,y$ this is also $\CR(y\tens y)$. But expanding $0=\CR(xy\tens x)$ forces them all to be zero, so only the trivial quasitriangular structure for $\F_2\Z_3$.  \endproof

 We now classify quasitriangular structures on most of the $n=4$ Hopf algebras. Recall that many of these were covered by the construction $A_{1jk}$ in (\ref{Aijk}) or its dual.

\begin{proposition}\label{qt4} (i) The Hopf algebras $A_2=A_{101}$,  $A_2^*$, $A_{111}$,  $A_{111}^*$, $\F_2\Z_2^2$, $\F_2(\Z_2^2)$, $\F_2\Z_4$, $\F_2(\Z_4)$ and the ones of type G$\to$L, L$\to$ G,  admit no nontrivial quasitriangular structure. 

(ii) The anyonic line $A_{100}$, its dual $A_{100}^*$,    $\F_2(\Z_2)\tens {\rm gra}=A_{110}$ and its dual  $\F_2\Z_2\tens {\rm gra}=A_{110}^*$ each admit a unique nontrivial quasitriangular structure. This is triangular and inherited from a natural Grassmann line sub-Hopf algebra in each case.

(iii) The self-dual version of the anyonic line, $G\to G$ in Proposition~\ref{propG} with relation $x^4=0$, admits only 4 quasitriangular structures (3 nontrivial)
  \[ \CR=1\tens 1+\alpha x^2\tens x^2+ \beta(x\tens x^2+x^2\tens x+x^3\tens x^3),\quad \alpha,\beta\in\{0,1\}.\]
These are all triangular.   

(iv) $\F_2(\Z_2)\tens \F_2\Z_2$ admits only 4 quasitriangular structures  (3 nontrivial) 
\[ \CR=(1\tens 1+ \alpha y\tens x)(1\tens 1+\beta x\tens y),\quad\alpha,\beta\in\{0,1\}\]
This is triangular {\em iff} $\alpha=\beta$ and  factorisable {\em iff} $\alpha\ne\beta$ (the `quantum double' $\CR$).  

(v) The Grassmann plane Hopf algebra admits only 16 quasitriangular structures (15 nontrivial)
\[ \CR=1\tens 1+r_{ij}x^i\tens x^j+ \det(r)xy\tens xy,\quad r\in M_2(\F_2),\]
where $x^1=x, x^2=y$ are the Grassmann variables. This is triangular {\em iff} $r$ is symmetric, and is never factorisable. \end{proposition}
\proof For $A_{1jk}$, the relation $w^4=jw^2+kw$ implies that $w^5=jw^3+kw^2$, $w^6=kw^3+jw^2+jkw$. We also need to take the form $\CR=1\tens 1+\sum_{a,b=1}^{a,b=3}R_{ab}w^a\tens w^b$ to obey the counit condition. We then write out the $(\Delta\tens\id)\CR$ and $(\id\tens\Delta)\CR$ conditions for a quasitriangular structure (the 3rd condition is automatic) and reduce all powers of $w$. This results in some quadratic equations for the coefficients $R_{ab}$, which we then solve for the different $j,k$. For (ii), the Grassmann line generators of   $A_{100}$ and $A_{110}$ are $w^2$ and $w^2+w$,  respectively.

For $A_{1jk}^*$,  we look for coquasitriangular structures on $A_{1jk}$. Since $w$ is primitive, $\CR(w^2\tens w)=\CR(w\tens w)\CR(w\tens 1)+\CR(w\tens 1)\CR(w\tens w)=0$ and in a similar way one finds all $\CR(w^a\tens w^b)=0$ for either $a,b>1$. Hence the only nontrivial option is $\CR(w\tens w)=1$. For this to be defined on the quotient we need $0=\CR(w\tens w^4)=\CR(w\tens j w^3)+\CR(w\tens k w)=k\CR(w\tens w)$. So this allows nonzero $\CR(w\tens w)$ precisely when $k=0$, giving unique nontrivial quasitriangular structures on the two duals stated. It easy enough to then identify them.  

For (iii), we similarly analyse $\CR$ on the same algebra $x^4=0$ as $A_{100}$ but with $x$ not primitive. Here $\Delta x^2=x^2\tens 1+1\tens x^2$ so $\CR(x^2\tens x^2)=0$. Analysing the other powers, we find that $\CR(x\tens x)=\alpha$ and $\CR(x^2\tens x)=\CR(x\tens x^2)=\CR(x^3\tens x^3)=\beta$ are unconstrained and the others are zero. We now use $\CR=\CR(x^\mu\tens x^\nu)y_\mu\tens y_\nu$ to translate to a quasitriangular structure for the dual, which is isomorphic to the same Hopf algebra as in the proof of Proposition~\ref{propG} by $y_2=x, y_1=y=x^2, y_3=z=x^3$. This gives $\CR$ as stated. It squares to itself and is symmetric in its tensor factors, hence triangular. 

For (iv), if this is viewed as a Drinfeld double \cite{Dri,Ma:book} $D(\F_2(\Z_2))$ then $\CR=1\tens 1+x\tens y$ is the canonical factorisable $\CR$ and its own inverse, so $\CR_{21}$ is another. These are the $\alpha\ne\beta$ cases.  To show that the four stated are all, since the Hopf algebra is self-dual,  we analyse coquasitriangular structures. Here $x^2=x$, with $x$ primitive and $y^2=0$ with $\Delta y=y\tens 1+1\tens y+y\tens y$ for the two tensor factors. Then $\CR(x\tens x)=\CR(x^2\tens x)=0$ and $0=\CR(y^2\tens y)=\CR(y\tens y)\CR(y\tens y)$ which requires $\CR(y\tens y)=0$. This implies $\CR(xy\tens x^i)=\CR(x^i\tens xy)=0$ where $x^i=x,y$ while $\CR(xy\tens xy)=\CR(x\tens y)\CR(y\tens x)$. This then dualises as stated. We then compute $\CQ=1\tens 1+(\alpha+\beta)(x\tens y+y\tens x+ xy\tens xy).$ 

Similarly for (v), this is again self-dual so we analyse coquasitriangular structures. Looking at all products we can set $\CR(x\tens x)$, $\CR(x\tens y)$, $\CR(y\tens x)$, $\CR(y\tens y)$ freely, $\CR(xy\tens x^i)=\CR(x^i,xy)=0$ for $x^i=x,y$ and $\CR(xy\tens xy)=\CR(x\tens x)\CR(y\tens y)+\CR(x\tens y)\CR(y\tens x)$. This immediately dualises to the stated quasitriangular structure. Note that the stated $\CR$ is its own inverse so the `conjugate' quasitriangular structure is given by $\CR_{21}$ which has the same form but for the transposed $r$. 

Finally, for $\F_2(\Z_2^2),\F_2(\Z_4),$ G$\to$L and their duals, we again look for coquasitriangular structures and find there are none. Proofs for the first two are similar to those of Proposition~\ref{qt3} and are omitted. For the non-primitive version G$\to$L of the anyonic line with $x^4=0$, we start with $\Delta x^2=x^2\tens 1+x^2\tens 1+x^2\tens x^2$ so that $0=\CR(x^4\tens x^2)=\CR(x^2\tens x^2)^2=\CR(x^2\tens x^2)$, then proceed to all powers. For the non-primitive version L$\to$G of $A_2$, we work in the algebra $x^4=x$ and coproduct (\ref{LtoG}), we start with $\CR(x^2\tens x)=\CR(x\tens (x+x^2))^2=\CR(x\tens( x+x^2))=\CR(x\tens x^2)+\CR(x\tens x)$ and $\CR(x^2\tens x^2)=\CR(x\tens (x+x^2))$ similarly, hence $\CR(x^2\tens x)=\CR(x^2\tens x^2)=\CR(x\tens x^2)$ by symmetry, hence $\CR(x\tens x)=0$. 
 \endproof
 
 Case (iii) here includes the non-involutive strict quasitriangular structure of the Drinfeld double $D(\F_2(\Z_2))$, inducing an R-matrix obeying the braid or Yang-Baxter equations in any representation. The remaining 3 Hopf algebras are noncommutative or noncocommutative (or both) and need more care. 
 
 \begin{proposition} 
  (i)  $c[B_+]$ admits no quasitriangular structure.
  
  (ii)  $c[B_+]^*$ admits a unique nontrivial quasitriangular structure. This is triangular and inherited from a natural Grassmann line sub-Hopf algebra.
   
(iii) $d_{sl_2}$ in Proposition~\ref{nfnf} admits only 2 quasitriangular structures, both nontrivial but triangular, namely
\[ \CR=1\tens 1+u\tens w+x\tens u+ (x+w)\tens (x+w)+\alpha u\tens u;\quad u=1+s,\quad \alpha\in\{0,1\}.\]
\end{proposition}
\proof In all cases it is easier to look for coquasitriangular structures on the dual. 

For (i), working on the noncommutative $c[B_+]^*$ described in Section~\ref{secNF}(i), $\CR(x\tens x)=\CR(x^2\tens x)=0$ as $x$ is primitive and $\CR(x\tens w)=0$ similarly. Also $\CR(x\tens y)=\CR(x^2\tens y)=2\CR(x\tens x)\CR(x\tens w)=0$, $\CR(y\tens x)=\CR(yx\tens x)=0$, $\CR(y\tens w)=\CR(yx\tens w)=0$ and $\CR(y\tens y)=\CR(yx\tens y)=\CR(y\tens x)\CR(x\tens w)+\CR(y\tens w)\CR(x\tens x)=0$, and so forth to show that only $\CR(1\tens 1)=1$ on the basis $1,x,y,w$. This translates to $\CR=1\tens 1$ for $c[B_+]$, but this is not cocommutative so the last of (\ref{Dri}) does not hold.

For (ii), working on $c[B_+]$ with algebra $s^2=1,y^2=0$, $s,y$ commuting and coalgebra (\ref{CBplus}),  $1=\CR(1\tens s)=\CR(s^2\tens s)=\CR(s\tens s)^2=\CR(s\tens s)$ and $0=\CR(s\tens y^2)=\CR(s\tens y)^2=\CR(s\tens y)$ and so forth. In this way, one arrives at $\CR(s\tens\ )=\CR(\ \tens s)=0$ when the space is $y$ or $sy$ and $\CR(y\tens y)=\CR(sy\tens y)=\CR(y\tens sy)=\CR(sy\tens sy)$. This gives potentially two coquasitriangular structures, and one can check that the third axiom also holds. This implies that $c[B_+]^*$ has the trivial and one other quasitriangular structure. We can identify the latter as $\CR=1\tens 1+w\tens w$ which  is that of the Grassmann line as a sub-Hopf algebra in the description in Section~\ref{secG}(i). One can check directly that the third of (\ref{Dri}) indeed holds on $\Delta x$ and $\Delta y$.

For (iii), working on $d_{sl_2}$ with the basis and algebra in Proposition~\ref{nfnf}, we apply the same methods as above to determine a coquasitriangular structure, being careful now that this is both noncommutative and noncocommutative. Details are omitted, but the first two axioms eventually show that 
\[ \CR(s\tens s)=1,\quad \CR(x\tens x)=\CR(w\tens x)=\CR(w\tens w)=\alpha\]
\[ \CR(s\tens x)=\CR(x\tens s)=\CR(s\tens w)=\CR(w\tens s)=\gamma,\quad \CR(x\tens w)=\alpha+\gamma;\quad \alpha\gamma=\alpha\]
for  $\alpha,\gamma\in \{0,1\}$. The third (quasi-commutativity axiom) then fixes $\gamma=1$. We then use $\CR=\CR(x^\mu\tens x^\nu)y_\mu\tens y_\nu$ and identify the $y_\rho$ as in the proof of Proposition~\ref{nfnf} to obtain the result stated. One can check that $\CR_{21}\CR=1$, so this is triangular. Note for this that $(x+w)^2=x^2+w^2+wx+xw=x+w+w+x=0$ and that $u^2=u(x+w)=(x+w)u=0$. 
\endproof

\section{Concluding remarks}\label{secconc}

We succeeded in determining all inequivalent bialgebras and Hopf algebras of dimension $n\le 4$ over $\F_2$. We presented our results in the form of extended graphs (\ref{n2quiver}), (\ref{n3quiver}) and Figure~\ref{n4big} with further details of the bialgebras available on \cite{github}. For Hopf algebras alone in Figure~\ref{summary},  we identified or described all 25  of them during the paper  and Appendix~\ref{AppB}. One important lesson is that while it is common practice to refer to an algebra by its most important role, for example $\F_2\Z_2$ for the group algebra of the group $\Z_2$, and we did the same when introducing our algebras for the first time, we now see in Figure~\ref{summary} that is much better to think of these as labels of the arrows, {\em not} of the nodes. Thus $\F_2\Z_2$ is one arrow out of the $n=2$ algebra A and the Grassmann line is another entirely different arrow. This gives a much clearer view of both familiar and unfamiliar Hopf algebras by the time we come to $n=4$, with indeed seven different arrows coming out of the 2-variable Grassmann algebra E with $x^2=y^2=0$, only one of which is the Grassmann plane.  Likewise on the algebra G with $x^4=0$ we identified four arrows, one of which is the anyonic line but another was a  self-dual Hopf algebra as a non-linear version of it. Also self-dual was a unique noncommutative noncocommutative Hopf algebra, $d_{sl_2}$, as some kind of digital quantum group.

We then looked at the canonical Fourier transform on each of our Hopf algebras, viewed as linear maps attached to the arrows, or a quiver algebra representation with monodromy at least for the representative Hopf algebras chosen. We also looked at the digital version of Drinfeld's quasitriangular Hopf algebra theory with very few found for $n\le 3$ but rather more for $n=4$, although mostly triangular. The number of Hopf algebras was small enough that we could proceed algebraically, but one could also look at the full answer for all $n=4$ quasitriangular bialgebras using computer methods and the data in \cite{github}. This should be interesting, given that there are many more noncommutative and noncocommutative bialgebras according to the counting in Figure~\ref{n4big}.  

There are many further directions to go from here. Obviously, it would be interesting to know about digital quantum groups for $n\ge 5$. Certainly, many can be constructed, for example applying the quantum double construction \cite{Dri,Ma:book} to $c[B_+]$ (which is noncocommutative) and to $d_{sl_2}$, will give 16 dimensional ones. The $A_{ijk}$ construction  (\ref{Aijk})  also works similarly for  general $A_{i_1\cdots i_d}$\cite{BasMa}.  To classify higher dimensions by our methods would, however, need much more powerful computers or a different approach. For $n=4$, we already had to search among 8,184 versions of our 25 distinct algebras to identify the latter as well as to identify the dual of each coalgebra solution for Sections~\ref{secn4comm} and~\ref{secn4noncomm}. Similarly, we searched among $20,160=|GL_4(\F_2)|$ potential maps between pairs of coalgebras of a given type in order to identify isomorphic bialgebras for Section~\ref{secisom}. The same methods can also be applied for $\F_3$ and above, for small $n$. The $A_{i_1\cdots i_d}$ construction in \cite{BasMa} works similarly over any $\F_p$ (as polynomials involving only powers of $x$ that are powers of $p$) to provide some basic examples, as do group algebras, functions on finite groups and algebraic groups, but we also expect many noncommutative and noncocommutative ones.  


It would also be interesting to look at the $\F_2$-linear or `digital' monoidal category of representations of a digital Hopf algebra. In the quasitriangular case, this is braided (symmetric in the triangular case). To give a flavour, up to equivalence,  $d_{sl_2}$ has four distinct representations which appear to generate the others by direct sums (this was checked to dimension $\le 4$). We denote these according to dimension by $1,\overline{1}, 2,\overline{2}$,  where $\rho_1=\eps$ and $\rho_{\overline{1}}(x)=\rho_{\overline{1}}(w)=\rho_{\overline{1}}(s)=1$ and, as representatives, 
\[ \rho_2(x)=\begin{pmatrix}1&1\\0&0\end{pmatrix},\ \rho_2(w)=\begin{pmatrix}0&0\\1&1\end{pmatrix},\ \rho_2(s)=\begin{pmatrix}0&1\\1&0\end{pmatrix};\quad \rho_{\overline{2}}(x)=\rho_{\overline{2}}(w)=\begin{pmatrix}1&0\\1&0\end{pmatrix},\ \rho_{\overline{2}}(s)=\begin{pmatrix}0&1\\1&0\end{pmatrix}.\]
Here $2$ descends to $e=x+w+s$, $\overline{2}$ descends to $x=w$ and $2\oplus\overline{2}$ is equivalent to 
the left regular representation. Also note that $1\subset 2$ and $\bar 1\subset \bar 2$ spanned by the vector with entries 1, so these are indecomposable but not irreducible (the algebra is not semisimple). The Hopf algebra structure then adds rules for dual representations and tensor products, which up to equivalence work out  as $2,\overline{2}$ dual to each other, $1,\overline{1}$ self-dual and 
\[ \begin{array}{c|cccc}\tens &\quad 1\quad  &\quad \overline{1}\quad & \quad 2\quad  &\quad  \overline{2}\quad  \\ \hline 1 & 1 & \overline{1} & 2 & \overline{2}\\  \overline{1} & \overline{1} & 1 & \overline{2} & 2\\ 
 2 & 2 & \overline{2} & 2\oplus\overline{2} & 2\oplus \overline{2}\\  \overline{2} & \overline{2} & 2 & 2\oplus\overline{2} & 2\oplus\overline{2}.\end{array}\]
The actual representations are many more, for example there are 20 in dimension 2 and 394 in dimension 3. We also did not have room here to cover ribbon structures and the transmutation from digital quasitriangular Hopf algebras to digital braided groups \cite{Ma:book} and braided-Fourier transform \cite{LyuMa}. The latter over $\C$ is at the heart of and leads naturally into topological quantum field theories and knot and 3-manifold invariants, hence these should all have digital versions. None of this is too surprising but should be interesting to work out in detail. Likewise over $\F_p$. A variant of such TQFT's is the Kitaev model for topological quantum computing \cite{KitNet,Kit} and it would be interesting to see and possibly actually build a digital version based on $d_{sl_2}$. This will be looked at in a sequel.  


Beyond quantum groups, we could of course relax our axioms and classify small digital weak Hopf algebras, quasi-Hopf algebras and Hopf quasigroups, for example.  As mentioned in the introduction, our companion paper \cite{MaPac2} already classified digital quantum Riemannian geometries to dimension $n=3$, with some partial results for $n=4$. There are also interesting interactions between the Hopf algebra duality, the curvature and de Morgan duality \cite{Ma:dem}. The digital quantum groups in the present paper provide many more examples once equipped with bicovariant differential structures and quantum metrics, which is another specific direction for further work. 

\begin{figure}
\[\includegraphics[scale=0.75]{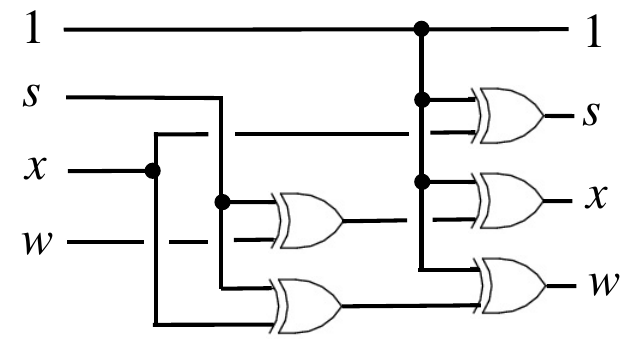}
\]
\caption{Fourier transform $\CF$ on $d_{sl_2}$ in Proposition~\ref{ft} as digital electronics. Applied three times gives the identity. \label{gates}}
\end{figure}

Finally, the constructions here may potentially be of interest to build in digital electronics. Thus both Fourier transform processes and Yang-Baxter solution R-matrices could be built in silicon as potential elements of digital quantum computers. An example is shown in Figure~\ref{gates} in terms of XOR gates. An element of $d_{sl_2}$ is input at left with a signal at wire $s$ say meaning the element includes $s$. The output similarly appears at right. Since $\CF$ is a linear operation, the effect of $1$ on the input is to add $\CF(1)=1+s+x+w$ to the output which inverts what would otherwise have been the signal output at $s,x,w$ (this is the bank of 3 XOR gates). Similarly $\CF(s)=x+w$ means that the presence of $s$ on the input inverts $x,w$, which is the bank of 2 XOR gates. Prior to that, $\CF(w)=x$ and $\CF(x)=s+w$ means we wire as shown. In this way we can translate natural operations from the quantum group theory into digital processes.  Although without the full benefits of true quantum computer gates over $\C$, we can gain experience with digital versions and meanwhile we have potential for diverse applications, possibly in signal processing and real-time encryption.

\appendix
\section{Data for $n=3$ bialgebras classification}
\label{AppA}
Here we list all possible coalgebras forming bialgebras on each of the four $n=3$ algebras for which 
these exist, obtained using Mathematica. This data is then analysed in Section~\ref{secn3} to identify the bialgebra isomorphism classes. In the following, B.1 -- B.33, for example, means we fix $V^{\mu\nu}{}_\rho$ as B in its standard form with basis $1,x,y$ and solve for all possible $C^\mu{}_{\nu\rho}$, listing the solutions. For each solution, we give the dual algebra and identify which algebra it is isomorphic to,  using R. The unit in the dual algebra is $1=\eps^\mu y_\mu$, where $\<x^\mu,y_\nu\>=\delta^\mu{}_\nu$ for the dual basis $y_\nu$.

\smallskip

\textbf{B.1. }
$\Delta x =1\otimes x+x\otimes 1+x\otimes
x+x\otimes y+y\otimes x+y\otimes y,\quad \Delta y =1\otimes
y+y\otimes 1,\quad \epsilon x =0=\epsilon y.$

Dual is commutative  algebra  C with  $1=y_0$, $y_1^2=y_1,\ y_1y_2=y_1,\ y_2^2=y_1$.

\textbf{B.2. }$\quad \Delta x =1\otimes x+x\otimes 1+x\otimes
x,\quad \Delta y =1\otimes y+y\otimes 1+x\otimes y+y\otimes x$
,\quad $\epsilon x =0=\epsilon y.$

Dual is commutative algebra C  with  $1=y_0$, 
 $y_1^2=y_1,\ y_1y_2=y_2,\ y_2^2=0$.

\textbf{B.3. }$\quad \Delta x =1\otimes x+x\otimes 1+x\otimes
x+y\otimes y,\quad \Delta y =1\otimes y+y\otimes 1+x\otimes
y+y\otimes x$,\quad $\epsilon x =0=\epsilon y.$

Dual is commutative algebra C with  $1=y_0$, $y_1^2=y_1,\ y_1y_2=y_2,\ y_2^2=y_1$.

\textbf{B.4. (Hopf algebra)}$\quad \Delta x =1\otimes x+x\otimes 1+x\otimes
y+y\otimes x+y\otimes y,$\\
$\Delta y =1\otimes y+y\otimes
1+x\otimes y+y\otimes x+x\otimes x$,\quad $\epsilon x =0=\epsilon y, \quad Sx =y,\quad Sy =x.$

Dual is commutative algebra  D  with  $1=y_0$, $y_1^2=y_2,\ y_1y_2=y_1+y_2,\ y_2^2=y_1$. 

\textbf{B.5. }$\quad \Delta x =1\otimes x+x\otimes 1+x\otimes
y+y\otimes x,\quad \Delta y =1\otimes y+y\otimes 1+y\otimes y$,\quad $\epsilon x =0=\epsilon y.$

Dual is commutative algebra C  with  $1=y_0$,  $y_1^2=0,\ y_1y_2=y_1,\ y_2^2=y_2.$

\textbf{B.6. }$\quad \Delta x =1\otimes x+x\otimes
1+x\otimes y+y\otimes x+x\otimes x,\quad \Delta y =1\otimes
y+y\otimes 1+y\otimes y$,\quad $\epsilon x =0=\epsilon y.$

Dual is commutative algebra B  with  $1=y_0$,  $y_1^2=y_1,\ y_1y_2=y_1,\ y_2^2=y_2.$

\textbf{B.7. }$\quad \Delta x =1\otimes x+x\otimes 1+x\otimes
y+y\otimes x,\quad \Delta y =1\otimes y+y\otimes 1+x\otimes
x+y\otimes y$,\quad $\epsilon x =0=\epsilon y.$

Dual is commutative algebra C   with  $1=y_0$,  $y_1^2=y_2,\ y_1y_2=y_1,\ y_2^2=y_2.$

\textbf{B.8.}$\quad \Delta x =1\otimes x+x\otimes 1+x\otimes
x+y\otimes x,\quad \Delta y =1\otimes y+y\otimes 1+x\otimes
y+y\otimes y$,\quad $\epsilon x =0=\epsilon y.$ 

Dual is noncommutative algebra G with  $1=y_0$,  $y_1^2=y_1,\ y_1y_2=y_2,\ y_2y_1=y_1$, $ 
y_2^2=y_2.$

\textbf{B.9.}$\quad \Delta x =1\otimes x+x\otimes 1+x\otimes
x+x\otimes y,\quad \Delta y =1\otimes y+y\otimes 1+y\otimes
x+y\otimes y$,\quad $\epsilon x =0=\epsilon y.$

Dual is noncommutative algebra G   with  $1=y_0$,  $y_1^2=y_1,\ y_1y_2=y_1,\ y_2y_1=y_2$, $y_2^2=y_2.$

\textbf{B.10. }$\Delta x =1\otimes
x+x\otimes 1+x\otimes x,\quad \Delta y =1\otimes y+y\otimes
1+x\otimes y+y\otimes x+y\otimes y,\quad \epsilon x =0=\epsilon y.$

Dual is commutative algebra   B  with  $1=y_0,$  $y_1^2=y_1,\ y_1y_2=y_2,\  y_2^2=y_2.$

\textbf{B.11. }$\Delta x =1\otimes x+x\otimes 1,\quad
\Delta y =1\otimes y+y\otimes 1+x\otimes y+y\otimes x+x\otimes
x+y\otimes y,\quad \epsilon x =0=\epsilon y.$ 

Dual is commutative algebra C  with  $1=y_0$,  $y_1^2=y_2,\ y_1y_2=y_2,\ y_2^2=y_2.$

\textbf{B.12. }$\quad \Delta x =x\otimes x,\quad \Delta
y =1\otimes y+y\otimes 1$,\quad $\epsilon x =1,\ \epsilon y =0.$

Dual is commutative algebra  C  with  $1=y_0+y_1$, $y_0^2=y_0,\ y_0y_1=0,\ y_0y_2=y_2$, $y_1^2=y_1,\ y_1y_2=0,\ y_2^2=0$.

\textbf{B.13. }$\quad \Delta x =x\otimes x,\quad \Delta
y =1\otimes 1+1\otimes x+x\otimes 1+1\otimes y+y\otimes
1+x\otimes x$,\quad $\epsilon x =1,\ \epsilon y =0.$

Dual is commutative algebra C  with  $1=y_0+y_1$, $y_0^2=y_0+y_2,\ y_0y_1=y_2,\ y_0y_2=y_2$, $y_1^2=y_1+y_2,\ y_1y_2=0,\ y_2^2=0.$

\textbf{B.14.}$\quad \Delta x =x\otimes x,\quad \Delta
y =x\otimes y+y\otimes 1$,\quad $\epsilon x =1,\ \epsilon y =0.$

Dual is noncommutative algebra G  with $1=y_0+y_1$, 

$y_0^2=y_0,\ y_0y_1=0=y_1y_0,\ y_0y_2=0,\ y_2y_0=y_2,\ 
y_1^2=y_1,\ y_1y_2=y_2,\ y_2y_1=0,\ y_2^2=0.$

\textbf{B.15.}$\quad \Delta x =x\otimes x,\quad \Delta
y =1\otimes y+y\otimes x$,\quad $\epsilon x =1,\ \epsilon y =0.$

Dual is noncommutative algebra G with $1=y_0+y_1$,

$y_0^2=y_0,\ y_0y_1=0=y_1y_0,\ y_0y_2=y_2,\ y_2y_0=0,\ 
y_1^2=y_1,\ y_1y_2=0,\ y_2y_1=y_2,\ y_2^2=0.$

\textbf{B.16. }$\quad \Delta x =x\otimes x,\quad \Delta
y =x\otimes y+y\otimes x$,\quad $\epsilon x =1,\ \epsilon y =0.$

Dual is commutative algebra C  with $1=y_0+y_1$, $y_0^2=y_0,\ y_0y_1=0,\ y_0y_2=0$, $
y_1^2=y_1,\ y_1y_2=y_2,\ y_2^2=0.$

\textbf{B.17. }$\quad \Delta x =x\otimes x+y\otimes y,\quad
\Delta y =x\otimes y+y\otimes x$,\quad $\epsilon x =1,\ \epsilon y =0;$

Dual is commutative algebra  C  with $1=y_0+y_1$, $y_0^2=y_0,\ y_0y_1=0,\ y_0y_2=0$, $ 
y_1^2=y_1,\ y_1y_2=y_2,\ y_2^2=y_1.$

\textbf{B.18. }$\quad \Delta x =x\otimes x,\quad \Delta
y =1\otimes 1+1\otimes x+x\otimes 1+x\otimes x+x\otimes
y+y\otimes x$,\quad $\epsilon x =1,\ \epsilon y =0.$ 

Dual is commutative  algebra C  with $1=y_0+y_1$, $y_0^2=y_0+y_2,\ y_0y_1=y_2,\ y_0y_2=0$, $
 y_1^2=y_1+y_2,\ y_1y_2=y_2,\ y_2^2=0.$

\textbf{B.19. }$\quad \Delta x =x\otimes x,\quad \Delta
y =1\otimes y+y\otimes 1+y\otimes y$,\quad $\epsilon x =1,\ \epsilon y =0.$

Dual is commutative  algebra B  with $1=y_0+y_1$, $y_0^2=y_0,\ y_0y_1=0,\ y_0y_2=y_2$, $y_1^2=y_1,\ y_1y_2=0,\ y_2^2=y_2$.

\textbf{B.20. }$ \Delta x =1\otimes 1+1\otimes x+1\otimes
y+x\otimes 1+x\otimes y+y\otimes 1+y\otimes x+y\otimes y,\\ \quad 
\Delta y =1\otimes y+y\otimes 1+y\otimes y,\quad \epsilon x =1,\ \epsilon y =0.$

Dual is commutative  algebra C  with  $1=y_0+y_1$,  $y_0^2=y_0+y_1,\ y_0y_1=y_1,\ y_0y_2=y_1+y_2$,  $y_1^2=0,\  y_1y_2=y_1,\  y_2^2=y_1+y_2$.

\textbf{B.21. (Hopf algebra)}$\quad \Delta x =1\otimes y+x\otimes x+x\otimes
y+y\otimes 1+y\otimes x$,\\
$\Delta y =1\otimes 1+1\otimes
x+1\otimes y+x\otimes 1+x\otimes x+y\otimes 1+y\otimes y,\quad \epsilon x =1,\ \epsilon y =0,\quad Sx =x,\  Sy =1+x+y.$

Dual is commutative  algebra  D  with $1=y_0+y_1$, $y_0^2=y_0+y_2,\ y_0y_1=y_2,\ y_0y_2=y_1+y_2$, $y_1^2=y_1+y_2,\ y_1y_2=y_1$, $y_2^2=y_2$.

\textbf{B.22.}$\quad \Delta x =x\otimes x,\quad \Delta
y =x\otimes y+y\otimes x+y\otimes y$,\quad $\epsilon x =1,\ \epsilon y =0.$

Dual is commutative  algebra B  with $1=y_0+y_1$, $y_0^2=y_0,\ y_0y_1=0,\ y_0y_2=0$, 
$y_1^2=y_1,\ y_1y_2=y_2$, $y_2^2=y_2$.

\textbf{B.23. }$\Delta x =1\otimes x+x\otimes 1+x\otimes
x,\quad \Delta y =1\otimes 1+1\otimes x+x\otimes 1+1\otimes
y+y\otimes 1+x\otimes y+y\otimes x+x\otimes x$,\quad $\epsilon x =0,\ \epsilon y =1.$

Dual is commutative  algebra C  with  $1=y_0+y_2$, $y_0^2=y_0+y_2,\ y_0y_1=y_1+y_2,\ y_0y_2=y_2,$  $y_1^2=y_1+y_2,\ y_1y_2=y_2$,\ $y_2^2=0$.

\textbf{B.24. }$\quad \Delta x =1\otimes x+x\otimes 1,\quad
\Delta y =y\otimes y$,\quad $\epsilon x =0,\ \epsilon y =1.$

Dual is commutative  algebra  C  with $1=y_0+y_2$, $y_0^2=y_0,\ y_0y_1=y_1,\ y_0y_2=0$, $
y_1^2=0,\ y_1y_2=0$, $y_2^2=y_2$.

\textbf{B.25. } $\quad\Delta x =1\otimes x+x\otimes 1+x\otimes x,\quad \Delta
y =y\otimes y$,\quad $\epsilon x =0,\ \epsilon y =1.$

Dual is commutative  algebra B  with  $1=y_0+y_2$, $y_0^2=y_0,\ y_0y_1=y_1,\ y_0y_2=0$, $y_1^2=y_1,\ y_1y_2=0$, $y_2^2=y_2$.

\textbf{B.26.} $\quad\Delta x =1\otimes x+x\otimes y,\quad \Delta y
=y\otimes y$,\quad $\epsilon x =0,\ \epsilon y =1.$

Dual is noncommutative algebra G  with  $1=y_0+y_2$, 

$y_0^2=y_0,\ y_0y_1=y_1,y_1y_0=0,\ y_0y_2=0=y_2y_0,\ 
y_1^2=0,\ y_1y_2=y_1,\ y_2y_1=0$, $y_2^2=y_2$.

\textbf{B.27.}$\quad \Delta x =x\otimes 1+y\otimes x,\quad
\Delta y =y\otimes y$,\quad $\epsilon x =0,\ \epsilon y =1.$

Dual is noncommutative algebra G  with  $1=y_0+y_2$, 

$y_0^2=y_0,\ y_0y_1=0,\ y_1y_0=y_1,\ y_0y_2=0=y_2y_0,\ 
 y_1^2=0,\ y_1y_2=0,\ y_2y_1=y_1$, $y_2^2=y_2.$

\textbf{B.28. }$\quad \Delta x =x\otimes y+y\otimes x,\quad
\Delta y =y\otimes y$,\quad $\epsilon x =0,\ \epsilon y =1.$

Dual is commutative  algebra  C  with  $1=y_0+y_2$, $y_0^2=y_0,\ y_0y_1=0,\ y_0y_2=0$, $ 
y_1^2=0,\ y_1y_2=y_1$, $y_2^2=y_2.$

\textbf{B.29. }$\quad \Delta x =x\otimes x+x\otimes y+y\otimes
x,\quad \Delta y =y\otimes y$,\quad $\epsilon x =0,\ \epsilon y =1.$ 

Dual is commutative  algebra B  with  $1=y_0+y_2$,  $y_0^2=y_0,\ y_0y_1=0,\ y_0y_2=0$, $ 
y_1^2=y_1,\ y_1y_2=y_1$, $y_2^2=y_2.$

\textbf{B.30. }$\Delta x =1\otimes 1+1\otimes x+1\otimes
y+x\otimes 1+y\otimes 1+y\otimes y,\quad \Delta y =y\otimes y, 
\quad \epsilon x =0,\ \epsilon y =1.$ 

Dual is commutative  algebra  C  with  $1=y_0+y_2$,  $y_0^2=y_0+y_1,\ y_0y_1=y_1,\ y_0y_2=y_1$, $y_1^2=0,\ 
 y_1y_2=0$, $y_2^2=y_1+y_2.$

\textbf{B.31. }$\quad \Delta x =1\otimes 1+1\otimes y+x\otimes
y+y\otimes 1+y\otimes x+y\otimes y,\quad \Delta y =y\otimes y,\quad 
\epsilon x =0,\ \epsilon y =1.$ 

Dual is commutative  algebra  C  with  $1=y_0+y_2$,  $y_0^2=y_0+y_1,\ y_0y_1=0,\ y_0y_2=y_1$, $y_1^2=0,\ y_1y_2=y_1$, $y_2^2=y_1+y_2.$

\textbf{B.32. }$\quad \Delta x =x\otimes y+y\otimes x,\quad
\Delta y =x\otimes x+y\otimes y$,\quad $\epsilon x =0,\ \epsilon y =1.$ 

Dual is commutative  algebra C  with  $1=y_0+y_2$, $y_0^2=y_0,\ y_0y_1=0,\ y_0y_2=0,$ 
$y_1^2=y_2,\ y_1y_2=y_1$, $y_2^2=y_2.$

\textbf{B.33. (Hopf algebra)}$\quad \Delta x =1\otimes 1+1\otimes x+1\otimes
y+x\otimes 1+y\otimes 1+x\otimes x+y\otimes y,$

$\Delta y =1\otimes x+x\otimes 1+x\otimes y+y\otimes
x+y\otimes y$, \quad $\epsilon x =0,\ \epsilon y =1, \quad
Sx =1+x+y,\quad Sy =y.$

Dual is commutative  algebra D  with  $1=y_0+y_2$,  $y_0^2=y_0+y_1,\ y_0y_1=y_1+y_2,\ y_0y_2=y_1$, $y_1^2=y_1,\  y_1y_2=y_2$,  $y_2^2=y_1+y_2$.

\smallskip

\textbf{C.1.}$\quad \Delta x =1\otimes x+x\otimes 1+x\otimes
x,\quad \Delta y =1\otimes y+y\otimes 1+x\otimes y+y\otimes x$,\quad $\epsilon x =0=\epsilon y.$ 

Dual is commutative  algebra  C  with  $1=y_0,$  $y_1^2=y_1,\ y_1y_2=y_2,\ y_2^2=0.$

\textbf{C.2.}$\quad \Delta x =1\otimes x+x\otimes 1+x\otimes
x,\quad \Delta y =1\otimes y+y\otimes 1+x\otimes y+y\otimes
x+y\otimes y$,\quad $\epsilon x =0=\epsilon y.$ 

Dual is commutative  algebra  B with  $1=y_0$,  $y_1^2=y_1,\ y_1y_2=y_2,\ y_2^2=y_2.$

\textbf{C.3.} $\quad \Delta x =x\otimes x,\quad
\Delta y =1\otimes y+y\otimes 1$,\quad $\epsilon x =1,\ \epsilon y =0.$ 

Dual is commutative  algebra  C with $1=y_0+y_1,$  $y_0^2=y_0,\ y_0y_1=0,\ y_0y_2=y_2$, $y_1^2=y_1,\ y_1y_2=0,\ y_2^2=0.$

\textbf{C.4.}$\quad \Delta x =x\otimes x,\quad \Delta
y =x\otimes y+y\otimes 1$,\quad $\epsilon x =1,\ \epsilon y =0.$ 

Dual is noncommutative algebra  G with  $1=y_0+y_1,$ 

$y_0^2=y_0,\ y_0y_1=0,\ y_0y_2=0,\ y_1y_0=0,\ 
y_1^2=y_1,\ y_1y_2=y_2,\ y_2y_0=y_2,\ y_2y_1=0,\ y_2^2=0.$

\textbf{C.5.}$\quad \Delta x =x\otimes x,\quad \Delta
y =1\otimes y+y\otimes x$,\quad $\epsilon x =1,\ \epsilon 
y =0.$

Dual is noncommutative algebra  G with  $1=y_0+y_1,$ 

$y_0^2=y_0,\ y_0y_1=0,\ y_0y_2=y_2,\ y_1y_0=0,\ 
y_1^2=y_1,\ y_1y_2=0,\ y_2y_0=0,\ y_2y_1=y_2,\ y_2^2=0.$

\textbf{C.6.}$\quad \Delta x =x\otimes x,\quad \Delta y =x\otimes y+y\otimes x$,\quad $\epsilon x =1,\ \epsilon y =0.$ 

Dual is commutative algebra  C with $1=y_0+y_1$, $y_0^2=y_0,\ y_0y_1=0,\ y_0y_2=0$, $ 
y_1^2=y_1,\ y_1y_2=y_2,\ y_2^2=0.$

\textbf{C.7.}$\quad \Delta x =x\otimes x,\quad \ \Delta
y =1\otimes y+y\otimes 1+y\otimes y$,\quad $\epsilon x =1,\ \epsilon y =0.$ 

Dual is commutative algebra  B with $1=y_0+y_1$,  $y_0^2=y_0,\ y_0y_1=0,\ y_0y_2=y_2$, $ 
y_1^2=y_1,\ y_1y_2=0,\ y_2^2=y_2.$

\textbf{C.8.}$\quad \Delta x =x\otimes x,\quad \Delta 
y =x\otimes y+y\otimes x+y\otimes y$,\quad $\epsilon x =1,\ \epsilon y =0.$ 

Dual is commutative algebra  B with $1=y_0+y_1,$  $y_0^2=y_0,\ y_0y_1=0,\ y_0y_2=0$, $
y_1^2=y_1,\ y_1y_2=y_2,\ y_2^2=y_2.$

\smallskip

\textbf{D.1.} \textbf{(Hopf algebra)}
$\Delta x =1\otimes x+x\otimes 1+x\otimes x$,
$\Delta y =1\otimes y+y\otimes 1+y\otimes y$,\\
 $\epsilon x =0=\epsilon y,\quad Sx =y,\quad Sy =x.$
 
Dual is commutative algebra  B with  $1=y_0$, $y_1^2=y_1,\ y_1y_2=0,\ y_2^2=y_2.$ 

\textbf{D.2.}\quad $\Delta x =1\otimes x+x\otimes 1+x\otimes
x+y\otimes x,\quad \Delta y =1\otimes y+y\otimes 1+x\otimes
y+y\otimes y$,\quad $\epsilon x =0=\epsilon y. $

Dual is noncommutative  algebra  G with  $1=y_0$, $y_1^2=y_1,\ y_1y_2=y_2,\ y_2y_1=y_1$, $y_2^2=y_2. 
$

\textbf{D.3.}$\quad \Delta x =1\otimes x+x\otimes 1+x\otimes
x+x\otimes y,\quad \Delta y =1\otimes y+y\otimes 1+y\otimes
x+y\otimes y$,\quad $\epsilon x =0=\epsilon y.$

Dual is noncommutative  algebra  G with  $1=y_0$, $
y_1^2=y_1,\ y_1y_2=y_1$,  $y_2y_1=y_2$, $y_2^2=y_2. 
$

\smallskip

\textbf{G.1.}$\quad \Delta x=1\otimes x+x\otimes 1+x\otimes x,\quad \Delta y =1\otimes y+y\otimes 1+x\otimes y+y\otimes x$, \quad $\epsilon x
=0=\epsilon y.$

Dual is commutative algebra   C with  $1=y_0$, $y_1^2=y_1,\  y_1y_2=y_2,\ y_2^2=0.$

\textbf{G.2}$\quad \Delta x
=1\otimes x+x\otimes 1+x\otimes x+y\otimes y,\quad 
\Delta y =1\otimes y+y\otimes 1+x\otimes y+y\otimes x$,\quad $\epsilon x
=0=\epsilon y.$

Dual is commutative algebra   C with  $1=y_0$, $y_1^2=y_1,\ y_1y_2=y_2,\ y_2^2=y_1.$

\textbf{G.3.}$\quad \Delta x
=1\otimes x+x\otimes 1+x\otimes x,\quad \Delta y =1\otimes y+y\otimes 1+x\otimes y+y\otimes x+y\otimes y$, \quad $\epsilon x
=0=\epsilon y$.

Dual is commutative algebra   B with  $1=y_0,$  $y_1^2=y_1,\ y_1y_2=y_2,\ y_2^2=y_2.$

\textbf{G.4.}$\quad \Delta x
=1\otimes x+x\otimes 1+x\otimes x+y\otimes y,\quad \Delta y =1\otimes y+y\otimes 1+x\otimes y+y\otimes x+y\otimes y$,\quad $\epsilon x
=0=\epsilon y.$

Dual is commutative algebra   D with  $1=y_0,$ 
 $y_1^2=y_1,\  y_1y_2=y_2,\ y_2^2=y_1+y_2.$ 

\textbf{G.5.}$\quad \Delta x
=x\otimes x,\quad\Delta y=x\otimes
y+y\otimes x$, \quad $\epsilon x
=1,\ \epsilon y =0.$

Dual is commutative algebra   C with $1=y_0+y_1$, $y_0^2=y_0,\ y_0y_1=0,\ y_0y_2=0$, $y_1^2=y_1,\ y_1y_2=y_2,\ y_2^2=0.$

\textbf{G.6.}$\quad \Delta x
=x\otimes x+y\otimes y,\quad \Delta y
=x\otimes y+y\otimes x$,\quad $\epsilon x
=1,\ \epsilon y =0.$

Dual is commutative algebra   C with $1=y_0+y_1$,  $y_0^2=y_0,\ y_0y_1=0,\ y_0y_2=0$,  $y_1^2=y_1,\ y_1y_2=y_2,\ y_2^2=y_1.$

\textbf{G.7.} 
$\quad \Delta x =x\otimes x,\quad \Delta y =x\otimes y+y\otimes x+y\otimes y$,\quad $\epsilon x
=1,\ \epsilon y =0.$

Dual is commutative algebra   B with $1=y_0+y_1$, $y_0^2=y_0,\ y_0y_1=0,\ y_0y_2=0$, $y_1^2=y_1,\  y_1y_2=y_2,\  y_2^2=y_2.$

\textbf{G.8.}$\quad \Delta x
=x\otimes x+y\otimes y,\quad \Delta y
=x\otimes y+y\otimes x+y\otimes y$,\quad $\epsilon x
=1,\ \epsilon y =0.$ 

Dual is commutative algebra   D with $1=y_0+y_1$, $y_0^2=y_0,\ y_0y_1=0,\  y_0y_2=0$,  
$y_1^2=y_1,\  y_1y_2=y_2,\ y_2^2=y_1+y_2. $

\section{Standard forms and Fourier transforms for $n=4$}\label{AppB}

Here we fix a standard form representative for each of the 20 distinct Hopf algebras for $n=4$, shown in a series of tables  according to algebra.  The 3rd column includes  an isomorphism between the dual Hopf algebra with dual basis $y_\mu$ and the corresponding standard  form with basis $1,x,y,z$ (written for brevity as an identification). We give the Fourier transform both in its canonical form in the 2nd column and as a `Fourier transport' map from the standard basis $1,x,y,z$ of one algebra to another. The four self-dual cases are in Proposition~\ref{ft} but included here with standard bases; the proof in general is similar, with the specific coproducts numbered according to \cite{github}. 

{\small {
\begin{tabular}{lll}
\begin{tabular}{l}
Algebra D, $\epsilon x=\epsilon y=\epsilon z=0$
\end{tabular}
& 
\begin{tabular}{l}
Dual algebra, integral, \\ 
Fourier transform
\end{tabular}
& 
\begin{tabular}{l}
Standard dual \\ 
 Fourier transport
\end{tabular}
\\ \hline
\begin{tabular}{l}
(D,D*),\quad ${\mathbb{F}}_{2}({\mathbb{Z}}_{2})\otimes {\mathbb{F}}_{2}{
\mathbb{Z}}_{2}$,\quad D.2: \\ 
\begin{tabular}{l}
$\Delta x=x\otimes 1+1\otimes x,$\newline
\\ 
$\Delta y=y\otimes 1+1\otimes y+y\otimes y,$ \\ 
$\Delta z=1\otimes z+x\otimes y+y\otimes x$ \\ 
$\qquad+z\otimes 1+z\otimes y+y\otimes z$\newline
\\ 
$Sx=x,\ Sy=y,\ Sz=z$.
\end{tabular}
\end{tabular}
& 
\begin{tabular}{l}
$
\begin{tabular}{l}
$y_{0}=1,\ y_{1}y_{2}=y_{3}=y_{2}y_{3},$\\
$ y_{2}^{2}=y_{2},\ y_{1}^{2}=y_{1}y_{3}=y_{3}^{2}=0,$ \\ 
$I=\left( 0,1,0,1\right) ,$ \\ 
${\mathcal{F}}=\left( 
\begin{array}{cccc}
$0$ & $1$ & $0$ & $1$ \\ 
$1$ & $1$ & $1$ & $1$ \\ 
$0$ & $1$ & $0$ & $0$ \\ 
$1$ & $1$ & $0$ & $0$
\end{array}
\right) $
\end{tabular}
$
\end{tabular}
& 
\begin{tabular}{l}
$y_{\mu }=(1,y,x,z),$ \\ 
${\mathcal{F}}_{\mathrm{D\rightarrow D}}=\left( 
\begin{array}{cccc}
0& 0 & 1 &1  \\ 
1& 1 & 1 &1  \\ 
0& 0 & 1 &0  \\ 
1& 0 & 1 & 0
\end{array}
\right) $
\end{tabular}
\\ \hline
\begin{tabular}{l}
(D,E*),\quad ${\mathbb{F}}_{2}({\mathbb{Z}}_{2})\otimes \mathrm{gra}$,\quad
D.1: \\ 
\begin{tabular}{l}
$\Delta x=x\otimes 1+1\otimes x$ \\ 
$\Delta y=y\otimes 1+1\otimes y,$\newline
\\ 
$\Delta z=1\otimes z+x\otimes y+y\otimes x+z\otimes 1$ \\ 
$Sx=x,\ Sy=y,\ Sz=z$.
\end{tabular}
\end{tabular}
& 
\begin{tabular}{l}
$y_{0}=1,\ y_{1}y_{2}=y_{3},$ \\ 
$y_{1}^{2}=y_{1}y_{3}=y_{2}^{2}=y_{2}y_{3}=y_{3}^{2}=0,$ \\ 
$
\begin{tabular}{l}
$I=\left( 0,0,0,1\right) ,$ \\ 
${\mathcal{F}}=\left( 
\begin{array}{cccc}
$0$ & $0$ & $0$ & $1$ \\ 
$0$ & $0$ & $1$ & $1$ \\ 
$0$ & $1$ & $0$ & $0$ \\ 
$1$ & $1$ & $0$ & $0$
\end{array}
\right) $
\end{tabular}
$
\end{tabular}
& 
\begin{tabular}{l}
$y_{\mu }=(1,y,x,z),$ \\ 
${\mathcal{F}}_{\mathrm{D\rightarrow E}}=\left( 
\begin{array}{cccc}
0& 0 & 0 & 1 \\ 
0&  1& 0 &1  \\ 
0&  0& 1 & 0 \\ 
1&  0& 1 & 0
\end{array}
\right) $
\end{tabular}
\end{tabular}
} }

{\small {
\begin{tabular}{lll}
\begin{tabular}{l}
Algebra E, $\epsilon x=\epsilon y=\epsilon z=0$%
\end{tabular}
& 
\begin{tabular}{l}
Dual algebra, integral, \\ 
Fourier transform
\end{tabular}
& 
\begin{tabular}{l}
Standard dual \\  Fourier transport
\end{tabular}
\\ \hline
\begin{tabular}{l}
(E,D*),\quad ${\mathbb{F}}_{2}{\mathbb{Z}}_{2}\otimes \mathrm{gra}$,\quad
E.2: \\ 
\begin{tabular}{l}
$\Delta x=1\otimes x+x\otimes 1$,\newline
\\ 
$\Delta y=1\otimes y+y\otimes 1+y\otimes y$,\newline
\\ 
$\Delta z=1\otimes z+x\otimes y+y\otimes x+y\otimes z$ \\ 
\qquad $+z\otimes 1+z\otimes y$, \\ 
$Sx=x,\ Sy=y,\ Sz=z$
\end{tabular}
\end{tabular}
& 
\begin{tabular}{l}
$
\begin{tabular}{l}
$y_{0}=1,$ $y_{1}y_{2}=y_{3}=y_{2}y_{3},$\\
$ y_{2}^{2}=y_{2},\  y_{1}^{2}=y_{1}y_{3}=y_{3}^{2}=0$, \\ 
$I=\left( 0,1,0,1\right) ,$ \\ 
${\mathcal{F}}=\left( 
\begin{array}{cccc}
$0$ & $1$ & $0$ & $1$ \\ 
$1$ & $0$ & $1$ & $0$ \\ 
$0$ & $1$ & $0$ & $0$ \\ 
$1$ & $0$ & $0$ & $0$
\end{array}
\right) $
\end{tabular}
$
\end{tabular}
& 
\begin{tabular}{l}
$y_{\mu }=(1,y,x,z)$, \\ 
${\mathcal{F}}_{\mathrm{E\rightarrow D}}=\left( 
\begin{array}{cccc}
$0$ & $0$ & $1$ & $1$ \\ 
$1$ & $1$ & $0$ & $0$ \\ 
$0$ & $0$ & $1$ & $0$ \\ 
$1$ & $0$ & $0$ & $0$
\end{array}
\right) $
\end{tabular}
\\ \hline
\begin{tabular}{l}
(E,E*),\quad Grass. plane, E.1: \\ 
\begin{tabular}{l}
$\Delta x=x\otimes 1+1\otimes x,$\newline
\\ 
$\Delta y=y\otimes 1+1\otimes y,$\newline
\\ 
$\Delta z=1\otimes z+x\otimes y+y\otimes x+z\otimes 1,$ \\ 
$Sx=x,\ Sy=y,\ Sz=z$\newline
\end{tabular}
\end{tabular}
& 
\begin{tabular}{l}
$y_{0}=1,$ $y_{1}y_{2}=y_{3}$,\\ 
$y_{1}^{2}=y_{1}y_{3}=y_{2}^{2}=y_{2}y_{3}=y_{3}^{2}=0$. \\ 
$
\begin{tabular}{l}
$I=\left( 0,0,0,1\right) $ \\ 
${\mathcal{F}}=\left( 
\begin{array}{cccc}
$0$ & $0$ & $0$ & $1$ \\ 
$0$ & $0$ & $1$ & $0$ \\ 
$0$ & $1$ & $0$ & $0$ \\ 
$1$ & $0$ & $0$ & $0$
\end{array}
\right) $
\end{tabular}
$
\end{tabular}
& 
\begin{tabular}{l}
$y_{\mu }=(1,x,y,z),$ \\ 
${\mathcal{F}}_{\mathrm{E\rightarrow E}}=\left( 
\begin{array}{cccc}
$0$ & $0$ & $0$ & $1$ \\ 
$0$ & $0$ & $1$ & $0$ \\ 
$0$ & $1$ & $0$ & $0$ \\ 
$1$ & $0$ & $0$ & $0$
\end{array}
\right) $
\end{tabular}
\\ \hline
\begin{tabular}{l}
(E,G*),\quad coanyonic line,\quad E.5: \\ 
\begin{tabular}{l}
$\Delta x=1\otimes x+x\otimes 1,$\newline
\\ 
$\Delta y=1\otimes y+x\otimes x+y\otimes 1,$\newline
\\ 
$\Delta z=1\otimes z+x\otimes y+y\otimes x+z\otimes 1,$ \\ 
$Sx=x,\ Sy=y,\ Sz=z$\newline
\end{tabular}
\end{tabular}
& 
\begin{tabular}{l}
$y_{0}=1,$ $y_{1}^{2}=y_{2},y_{1}y_{2}=y_{3},$\\
$y_{1}y_{3}=y_{2}^{2}=y_{2}y_{3}=y_{3}^{2}=0$, \\ 
$
\begin{tabular}{l}
$I=\left( 0,0,0,1\right) ,$ \\ 
${\mathcal{F}}=\left( 
\begin{array}{cccc}
$0$ & $0$ & $0$ & $1$ \\ 
$0$ & $0$ & $1$ & $0$ \\ 
$0$ & $1$ & $0$ & $0$ \\ 
$1$ & $0$ & $0$ & $0$
\end{array}
\right) $
\end{tabular}
$
\end{tabular}
& 
\begin{tabular}{l}
$y_{\mu }=(1,x,y,z)$, \\ 
${\mathcal{F}}_{\mathrm{E\rightarrow G}}=\left( 
\begin{array}{cccc}
$0$ & $0$ & $0$ & $1$ \\ 
$0$ & $0$ & $1$ & $0$ \\ 
$0$ & $1$ & $0$ & $0$ \\ 
$1$ & $0$ & $0$ & $0$
\end{array}
\right) $
\end{tabular}
\\ \hline
\begin{tabular}{l}
(E,L*),\quad $A_{2}^{\ast }$,\quad E.15: \\ 
\begin{tabular}{l}
$\Delta x=(1+z)\otimes x+x\otimes (1+z)+y\otimes y,$\newline
\\ 
$\Delta y=(1+z)\otimes y+y\otimes (1+z)+x\otimes x,$ \\ 
$\Delta z=1\otimes z+x\otimes y+y\otimes x+z\otimes 1+z\otimes z$\newline
, \\ 
$Sx=x,\ Sy=y,\ Sz=z$\newline
\end{tabular}
\end{tabular}
& 
\begin{tabular}{l}
$y_{0}=1,$ $y_{1}^{2}=y_{2}y_{3}=y_{2},$ \\
$y_{1}y_{2}=y_{3}^{2}=y_{3}, y_{1}y_{3}=y_{2}^{2}=y_{1},$ \\ 
$
\begin{tabular}{l}
$I=(1,0,0,1),$ \\ 
${\mathcal{F}}=\left( 
\begin{array}{cccc}
$1$ & $0$ & $0$ & $1$ \\ 
$0$ & $0$ & $1$ & $0$ \\ 
$0$ & $1$ & $0$ & $0$ \\ 
$1$ & $0$ & $0$ & $0$
\end{array}
\right) $
\end{tabular}
$
\end{tabular}
& 
\begin{tabular}{l}
$y_{\mu }=(1,x,z,1+y)$, \\ 
${\mathcal{F}}_{\mathrm{E\rightarrow L}}=\left( 
\begin{array}{cccc}
$0$& $0$ &$1$  &$0$  \\ 
$0$&$0$&$0$& $1$ \\ 
$0$& $1$ & $0$ &$0$  \\ 
$1$ & $0$&$0$  &$0$ 
\end{array}
\right) $
\end{tabular}
\\ \hline
\begin{tabular}{l}
(E,M*),\quad $A_{111}^{\ast }$,\quad E.16: \\ 
\begin{tabular}{l}
$\Delta x=(1+z)\otimes x+x\otimes (1+z)$\\
\qquad\  $+y\otimes y+z\otimes z$,\newline
\\ 
$\Delta y=(1+z)\otimes y+y\otimes (1+z)+x\otimes x$ \\ 
\qquad\ $+x\otimes z+y\otimes y+z\otimes x+z\otimes z,$\newline
\\ 
$\Delta z=1\otimes z+x\otimes y+y\otimes x+y\otimes z$ \\ 
\qquad\ $+z\otimes 1+z\otimes y+z\otimes z$,\\ 
$Sx=x,\ Sy=y,\ Sz=z$\newline
\end{tabular}
\end{tabular}
& 
\begin{tabular}{l}
$y_{0}=1,$ $y_{1}^{2}=y_{2},y_{1}y_{2}=y_{3},$ \\ 
$y_{1}y_{3}=y_{2}^{2}=y_{1}+y_{2},y_{2}y_{3}=y_{2}+y_{3},$\\
$y_{3}^{2}=y_{1}+y_{2}+y_{3}$,
\\ 
$
\begin{tabular}{l}
$I=\left( 1,1,0,1\right) ,$ \\ 
${\mathcal{F}}=\left( 
\begin{array}{cccc}
$1$ & $1$ & $0$ & $1$ \\ 
$1$ & $0$ & $1$ & $0$ \\ 
$0$ & $1$ & $0$ & $0$ \\ 
$1$ & $0$ & $0$ & $0$
\end{array}
\right) $
\end{tabular}
$
\end{tabular}
& 
\begin{tabular}{l}
$y_{\mu }=$\\
$(1,1+x+y,x+z,x)$, \\ 
${\mathcal{F}}_{\mathrm{E\rightarrow M}}=\left( 
\begin{array}{cccc}
$0$ &$0$  &$1$  &$0$  \\ 
$1$&$1$  & $0$ &$1$  \\ 
$1$& $1$ & $1$ & $0$ \\ 
$1$ & $0$ & $0$ & $0$
\end{array}\right) $
\end{tabular}
\\ \hline
\begin{tabular}{l}
(E,P*),\quad ${\mathbb{F}}_{2}{\mathbb{Z}}_{2}^{2}$,\quad E.38: \\ 
\begin{tabular}{l}
$\Delta x=1\otimes x+x\otimes 1+x\otimes x,$\newline
\\ 
$\Delta y=1\otimes y+y\otimes 1+y\otimes y,$\newline
\\ 
$\Delta z=\left( 1+x+y\right) \otimes z+z\otimes \left( 1+x+y\right) $ \\ 
$\qquad\quad +z\otimes z+y\otimes x+x\otimes y$,
\\ 
$Sx=x,\ Sy=y,\ Sz=z$\newline
\end{tabular}
\end{tabular}
& 
\begin{tabular}{l}
$y_{0}=1$, $y_{1}^{2}=y_{1}, y_{2}^{2}=y_{2},$ \\ 
$y_{1}y_{2}=y_{2}y_{3}=y_{3}^{2}=y_{1}y_{3}=y_{3}$, \\ 
$
\begin{tabular}{l}
$I=\left( 1,1,1,1\right) ,$ \\ 
${\mathcal{F}}=\left( 
\begin{array}{cccc}
$1$ & $1$ & $1$ & $1$ \\ 
$1$ & $0$ & $1$ & $0$ \\ 
$1$ & $1$ & $0$ & $0$ \\ 
$1$ & $0$ & $0$ & $0$
\end{array}
\right) $
\end{tabular}
$
\end{tabular}
& 
\begin{tabular}{l}
$y_{\mu }=(1,x,y,z)$, \\ 
${\mathcal{F}}_{\mathrm{E\rightarrow P}}=\left( 
\begin{array}{cccc}
$1$ & $1$ & $1$ & $1$ \\ 
$1$ & $0$ & $1$ & $0$ \\ 
$1$ & $1$ & $0$ & $0$ \\ 
$1$ & $0$ & $0$ & $0$
\end{array}
\right) $
\end{tabular}
\\ \hline
\begin{tabular}{l}
(E,NF*),\quad $c[B_{+}]$,\quad E.40: \\ 
\begin{tabular}{l}
$\Delta x=1\otimes x+x\otimes 1+x\otimes x,$\newline
\\ 
$\Delta y=1\otimes y+x\otimes y+y\otimes 1$, \\ 
$\Delta z= 1 \tens z+z \tens 1 + (y+z)\tens x + x \tens y$, \\  
$Sx=x,\ Sy=y+z,\ Sz=z$\newline
\end{tabular}
\end{tabular}
& 
\begin{tabular}{l}
$y_{0}=1$, $y_{1}^{2}=y_1,y_{1}y_{2}=y_{2}+y_{3},$ \\ 
$y_2y_1=y_{3}y_{1}=y_3,$\\
$y_1y_3=y_{2}^{2}=y_{2}y_{3}=y_{3}y_{2}=y_{3}^2=0$, \\ 
$
\begin{tabular}{l}
$I=\left( 0,0,1,1\right) ,$ \\ 
${\mathcal{F}}=\left( 
\begin{array}{cccc}
$0$ & $0$ & $1$ & $1$ \\ 
$0$ & $0$ & $1$ & $0$ \\ 
$1$ & $1$ & $0$ & $0$ \\ 
$1$ & $0$ & $0$ & $0$
\end{array}
\right) $
\end{tabular}
$
\end{tabular}
& 
\begin{tabular}{l}
$y_{\mu }=(1,x,y+z,y)$, \\ 
${\mathcal{F}}_{\mathrm{E\rightarrow NF}}=\left( 
\begin{array}{cccc}
$0$ & $0$ & $0$ &  $1$ \\ 
$0$ & $0$  &  $1$ &  $1$ \\ 
 $1$&  $1$ & $0$  & $0$  \\ 
 $1$& $0$ & $0$   & $0$ 
\end{array}
\right) $
\end{tabular}
\end{tabular}
} }

{\small {
\begin{tabular}{lll}
\begin{tabular}{l}
Algebra  G, $\epsilon x=\epsilon y=\epsilon z=0$
\end{tabular}
& 
\begin{tabular}{l}
Dual algebra, integral, \\ 
Fourier transform
\end{tabular}
& 
\begin{tabular}{l}
Standard dual \\  Fourier transport
\end{tabular}
\\ \hline
\begin{tabular}{l}
(G,E*),\quad {\rm anyonic\ line}\quad G.1: \\ 
\begin{tabular}{l}
$\Delta x=1\otimes x+x\otimes 1,$ \\ 
$\Delta y=1\otimes y+y\otimes 1,$ \\ 
$\Delta z=1\otimes z+x\otimes y+y\otimes x+z\otimes 1$, \\ 
$Sx=x,\ Sy=y,\ Sz=z$
\end{tabular}
\end{tabular}
& 
\begin{tabular}{l}
$
\begin{tabular}{l}
$y_{0}=1, y_{1}y_{2}=y_{3},$\\
$y_{1}^{2}=y_{1}y_{3}=y_{2}^{2}=y_{2}y_{3}=y_{3}^{2}=0,$ \\ 
$I=\left( 0,0,0,1\right) ,$ \\ 
${\mathcal{F}}=\left( 
\begin{array}{cccc}
$0$ & $0$ & $0$ & $1$ \\ 
$0$ & $0$ & $1$ & $0$ \\ 
$0$ & $1$ & $0$ & $0$ \\ 
$1$ & $0$ & $0$ & $0$
\end{array}
\right) $
\end{tabular}
$
\end{tabular}
& 
\begin{tabular}{l}
$y_{\mu }=(1,x,y,z)$, \\ 
${\mathcal{F}}_{\mathrm{G\rightarrow E}}=\left( 
\begin{array}{cccc}
$0$ & $0$ & $0$ & $1$ \\ 
$0$ & $0$ & $1$ & $0$ \\ 
$0$ & $1$ & $0$ & $0$ \\ 
$1$ & $0$ & $0$ & $0$
\end{array}
\right) $
\end{tabular}
\\ \hline
\begin{tabular}{l}
(G,G*),\quad G.2: \\ 
\begin{tabular}{l}
$\Delta x=1\otimes x+x\otimes 1+y\otimes y,$ \\ 
$\Delta y=1\otimes y+y\otimes 1,$ \\ 
$\Delta z=1\otimes z+x\otimes y+y\otimes x+z\otimes 1$, \\ 
$Sx=x,\ Sy=y,\ Sz=z$
\end{tabular}
\end{tabular}
& 
\begin{tabular}{l}
$y_{0}=1, y_{1}y_{2}=y_{3}, y_{2}^{2}=y_{1}$\\
$y_{1}^{2}=y_{1}y_{3}=y_{2}y_{3}=y_{3}^{2}=0$, \\ 
$
\begin{tabular}{l}
$I=\left( 0,0,0,1\right) ,$ \\ 
${\mathcal{F}}=\left( 
\begin{array}{cccc}
$0$ & $0$ & $0$ & $1$ \\ 
$0$ & $0$ & $1$ & $0$ \\ 
$0$ & $1$ & $0$ & $0$ \\ 
$1$ & $0$ & $0$ & $0$
\end{array}
\right) $
\end{tabular}
$
\end{tabular}
& 
\begin{tabular}{l}
$y_{\mu }=(1,y,x,z)$, \\ 
${\mathcal{F}}_{\mathrm{G\rightarrow G}}=\left( 
\begin{array}{cccc}
0& 0 &0  & 1 \\ 
0& 1 & 0 & 0 \\ 
0& 0 & 1 & 0 \\ 
1& 0 & 0 & 0
\end{array}
\right) $
\end{tabular}
\\ \hline
\begin{tabular}{l}
(G,P*),\quad $\F_2\Z_4$,\quad G.5: \\ 
\begin{tabular}{l}
$\Delta x=1\otimes x+x\otimes 1+x\otimes x,$ \\ 
$\Delta y=1\otimes y+y\otimes 1+y\otimes y,$ \\ 
$\Delta z=\left( 1+x+y\right) \otimes z+z\otimes \left( 1+x+y\right) $ \\ 
$\qquad\ +y\otimes x+x\otimes y+z\otimes z,$ \\ 
$Sx=x+y+z,\ Sy=y,\ Sz=z$
\end{tabular}
\end{tabular}
& 
\begin{tabular}{l}
$y_{0}=1,$ $y_{1}^{2}=y_{1},y_{1}y_{2}=y_{3}=y_{1}y_{3},$ \\ 
$y_{2}^{2}=y_{2},y_{2}y_{3}=y_{3}=y_{3}^{2}$, \\ 
$
\begin{tabular}{l}
$I=\left( 1,1,1,1\right) ,$ \\ 
${\mathcal{F}}=\left( 
\begin{array}{cccc}
$1$ & $1$ & $1$ & $1$ \\ 
$1$ & $1$ & $1$ & $0$ \\ 
$1$ & $1$ & $0$ & $0$ \\ 
$1$ & $0$ & $0$ & $0$
\end{array}
\right) $
\end{tabular}
$
\end{tabular}
& 
\begin{tabular}{l}
$y_{\mu }=(1,x,y,z),$ \\ 
${\mathcal{F}}_{\mathrm{G\rightarrow P}}=\left( 
\begin{array}{cccc}
$1$ & $1$ & $1$ & $1$ \\ 
$1$ & $1$ & $1$ & $0$ \\ 
$1$ & $1$ & $0$ & $0$ \\ 
$1$ & $0$ & $0$ & $0$
\end{array}
\right) $
\end{tabular}
\\ \hline
\begin{tabular}{l}
(G,L*),\quad G.6: \\ 
\begin{tabular}{l}
$\Delta x=1\otimes x+x\otimes 1+x\otimes x$ \\ 
$\qquad\quad+(y+z)\otimes (y+z),$ \\ 
$\Delta y=1\otimes y+y\otimes 1+y\otimes y,$ \\ 
$\Delta z=1\otimes z+x\otimes (y+z)+(y+z)\otimes x$ \\ 
$\qquad\ +z\otimes 1+y\otimes z+z\otimes y+z\otimes z$, \\ 
$Sx=x+y+z,\ Sy=y,\ Sz=z$
\end{tabular}
\end{tabular}
& 
\begin{tabular}{l}
$y_{0}=1,$ $y_{1}^{2}=y_{1}=y_{1}y_{3}, y_{1}y_{2}=y_{3}$ \\ 
$y_{2}^{2}=y_{1}+y_{2},y_{2}y_{3}=y_{1}+y_{3}=y_{3}^{2}$, \\ 
$
\begin{tabular}{l}
$I=(1,1,1,1),$ \\ 
${\mathcal{F}}=\left( 
\begin{array}{cccc}
$1$ & $1$ & $1$ & $1$ \\ 
$1$ & $1$ & $1$ & $0$ \\ 
$1$ & $1$ & $0$ & $0$ \\ 
$1$ & $0$ & $0$ & $0$
\end{array}
\right) $
\end{tabular}
$
\end{tabular}
& 
\begin{tabular}{l}
$y_{\mu }=$\\
$(1,x+z,z,1+x+y)$, \\ 
${\mathcal{F}}_{\mathrm{G\rightarrow L}}=\left( 
\begin{array}{cccc}
$0$ & $0$ & $1$  & $0$ \\ 
$1$ & $1$ & $0$ & $0$ \\ 
$1$ & $1$ & $0$  & $1$  \\ 
$1$ & $0$  & $0$  & $0$
\end{array}
\right) $
\end{tabular}
\end{tabular}
}}

{\small {
\begin{tabular}{lll}
\begin{tabular}{l}
Algebra L, $\epsilon x=0,\epsilon y=1,\epsilon z=0$
\end{tabular}
& 
\begin{tabular}{l}
Dual algebra, integral, \\ 
Fourier transform
\end{tabular}
& 
\begin{tabular}{l}
Standard dual \\  Fourier transport
\end{tabular}
\\ \hline
\begin{tabular}{l}
(L,E*),\quad $A_{2},$\quad L.6: \\ 
\begin{tabular}{l}
$\Delta x=1\otimes x+x\otimes 1,$ \\ 
$\Delta y=1\otimes 1+1\otimes y+x\otimes z+y\otimes 1$\\
$\qquad\quad+z\otimes x,$ \\ 
$\Delta z=1\otimes z+z\otimes 1$, \\ 
$Sx=x,\ Sy=y,\ Sz=z$
\end{tabular}
\end{tabular}
& 
\begin{tabular}{l}
$
\begin{tabular}{l}
$y_{0}^{2}=y_{0}+y_{2},y_{0}y_{1}=y_{1},$\\
$y_{0}y_{3}=y_{3}, y_{0}y_{2}=y_{1}y_{3}=y_{2},$\\
$y_{1}^{2}=y_{1}y_{2}=y_{2}^{2}=y_{2}y_{3}=y_{3}^{2}=0,$ \\ 
$I=\left( 0,0,1,0\right) ,$ \\ 
${\mathcal{F}}=\left( 
\begin{array}{cccc}
0 & 0 & 1 & 0 \\ 
0 & 0 & 0 & 1 \\ 
1 & 0 & 1 & 0 \\ 
0 & 1 & 0 & 0
\end{array}
\right) $
\end{tabular}
$
\end{tabular}
& $
\begin{tabular}{l}
$y_{\mu }=(1+z,x,z,y)$, \\ 
${\mathcal{F}}_{\mathrm{L\rightarrow E}}=\left( 
\begin{array}{cccc}
0&0  &0  & 1 \\ 
0& 0 & 1 & 0 \\ 
1 & 0 & 0 & 0 \\ 
0& 1 & 0 & 0
\end{array}
\right) $
\end{tabular}
$ \\ \hline
\begin{tabular}{l}
(L,G*),\quad L.11: \\  
\begin{tabular}{l}
$\Delta x=1\otimes x+x\otimes 1+x\otimes x+x\otimes z$\\ 
$\qquad\quad +z\otimes x+z\otimes z$, \\ 
$\Delta y=1\otimes 1+1\otimes y+x\otimes x+y\otimes 1+z\otimes z$,\\ 
$\Delta z=1\otimes z+x\otimes x+x\otimes z+z\otimes 1$\\
$\qquad\quad+z\otimes x+z\otimes z$,\\ 
$Sx=z,\ Sy=y,\ Sz=x$
\end{tabular}
\end{tabular}
& 
\begin{tabular}{l}
$y_{0}^{2}=y_{0}+y_{2},y_{0}y_{1}=y_{1},y_{0}y_{2}=y_{2},$\\
$y_{0}y_{3}=y_{3}, y_{1}^{2}=y_3^2=y_{1}+y_{2}+y_{3},$\\
$y_{1}y_{3}=y_{1}+y_{3}, y_{1}y_{2}=y_{2}^{2}=y_{2}y_{3}=0$ \\ 
$
\begin{tabular}{l}
$I=\left( 0,0,1,0\right) ,$ \\ 
${\mathcal{F}}=\left( 
\begin{array}{cccc}
0 & 0 & 1 & 0 \\ 
0 & 0 & 0 & 1 \\ 
1 & 0 & 1 & 0 \\ 
0 & 1 & 0 & 0
\end{array}
\right) $
\end{tabular}
$
\end{tabular}
& 
\begin{tabular}{l}
$y_{\mu }=$\\
$(1+z,x+z,z,x+y)$,\\ 
${\mathcal{F}}_{\mathrm{L\rightarrow G}}=\left( 
\begin{array}{cccc}
0&0  &0  &1  \\ 
0&1  &1  &0  \\ 
1& 0 &0  & 0 \\ 
0& 1 & 0 & 1
\end{array}
\right) $
\end{tabular}
\end{tabular}
}}

{\small {
\begin{tabular}{lll}
\begin{tabular}{l}
Algebra M, $\epsilon x=0,\epsilon y=1,\epsilon z=0$
\end{tabular}
& 
\begin{tabular}{l}
Dual algebra, integral, \\ 
Fourier transform
\end{tabular}
& 
\begin{tabular}{l}
Standard dual \\  Fourier transport
\end{tabular}
\\ \hline
\begin{tabular}{l}
(M,E*),\quad $A_{111}$,\quad M.2: \\ 
\begin{tabular}{l}
$\Delta x=1\otimes z+x\otimes (y+z)+(y+z)\otimes x$ \\ 
$\qquad\  +z\otimes 1+y\otimes z+z\otimes y,$ \\ 
$\Delta y=1\otimes1 +y\otimes z+z\otimes y$\\
$\qquad\ +1\tens \left(x+y+z\right) +\left( x+y+z\right) \otimes 1$
\\ 
$\qquad\ +x\otimes (y+ z)+(y+z)\otimes x$,\\
$\Delta z=1\otimes x+x\otimes 1+x\otimes (y+z)$ \\ 
$\qquad\quad +(y+z)\otimes x+y\otimes z+z\otimes y$, \\ 
$Sx=x,\ Sy=y,\ Sz=z$
\end{tabular}
\end{tabular}
& 
\begin{tabular}{l}
$
\begin{tabular}{l}
$y_{0}^{2}=y_{0}+y_{2},y_{0}y_{1}=y_{2}+y_{3},$ \\ 
$y_{0}y_{2}=y_{2},y_{0}y_{3}=y_{1}+y_{2},$ \\ 
$y_{1}y_{2}=y_{2}y_{3}=y_{1}y_{3}=y_{1}+y_{2}+y_{3},$ \\ 
$y_{1}^{2}=y_{2}^{2}=y_{3}^{2}=0,$ \\ 
$I=\left( 0,1,1,1\right) ,$ \\ 
${\mathcal{F}}=\left( 
\begin{array}{cccc}
0 & 1 & 1 & 1 \\ 
1 & 1 & 0 & 0 \\ 
1 & 0 & 1 & 0 \\ 
1 & 0 & 0 & 1
\end{array}
\right) $
\end{tabular}
$
\end{tabular}
& $
\begin{tabular}{l}
$y_{\mu }=$\\
$(1+x,x+y+z,x,y)$, \\ 
${\mathcal{F}}_{\mathrm{M\rightarrow E}}=\left( 
\begin{array}{cccc}
0& 0 & 0 & 1 \\ 
1 & 0 & 1 & 1 \\ 
1 & 0 & 0 & 0 \\ 
1 & 1 & 1 & 0
\end{array}
\right) $
\end{tabular}
$
\end{tabular}
}}

{\small {
\begin{tabular}{lll}
\begin{tabular}{l}
Algebra P, $\epsilon x=0,\epsilon y=0,\epsilon z=0$
\end{tabular}
& 
\begin{tabular}{l}
Dual algebra, integral, \\ 
Fourier transform
\end{tabular}
& 
\begin{tabular}{l}
Standard dual \\  Fourier transport
\end{tabular}
\\ \hline
\begin{tabular}{l}
(P,E*),\quad ${\mathbb{F}}_{2}({\mathbb{Z}}_{2}^{2})$,\quad P.1: \\ 
\begin{tabular}{l}
$\Delta x=1\otimes x+x\otimes 1,$ \\ 
$\Delta y=1\otimes y+y\otimes 1,$ \\ 
$\Delta z=1\otimes z+x\otimes y+y\otimes x+z\otimes 1,$ \\ 
$Sx=x,\ Sy=y,\ Sz=z$
\end{tabular}
\end{tabular}
& 
\begin{tabular}{l}
$
\begin{tabular}{l}
$y_{0}=1,y_{1}y_{2}=y_{3},$\\
$y_{1}^{2}=y_{1}y_{3}=y_{2}^{2}=y_{2}y_{3}=y_{3}^2=0$, \\ 
$I=\left( 0,0,0,1\right) ,$ \\ 
${\mathcal{F}}=\left( 
\begin{array}{cccc}
0 & 0 & 0 & 1 \\ 
0 & 0 & 1 & 1 \\ 
0 & 1 & 0 & 1 \\ 
1 & 1 & 1 & 1
\end{array}
\right) $
\end{tabular}
$
\end{tabular}
& $
\begin{tabular}{l}
$y_{\mu }=(1,x,y,z),$ \\ 
${\mathcal{F}}_{\mathrm{P\rightarrow E}}=\left( 
\begin{array}{cccc}
0 & 0 & 0 & 1 \\ 
0 & 0 & 1 & 1 \\ 
0 & 1 & 0 & 1 \\ 
1 & 1 & 1 & 1
\end{array}
\right) $
\end{tabular}
$ \\ \hline
\begin{tabular}{l}
(P,G*),\quad ${\mathbb{F}}_{2}({\mathbb{Z}}_{4})$,\quad P.3: \\ 
\begin{tabular}{l}
$\Delta x=1\otimes x+x\otimes 1,$ \\ 
$\Delta y=1\otimes y+x\otimes x+y\otimes 1,$ \\ 
$\Delta z=1\otimes z+x\otimes y+y\otimes x+z\otimes 1,$ \\ 
$Sx=x,\ Sy=y,\ Sz=z$
\end{tabular}
\end{tabular}
& 
\begin{tabular}{l}
$
\begin{tabular}{l}
$y_{0}=1,y_{1}^{2}=y_{2},y_{1}y_{2}=y_{3},$\\
$y_{1}y_{3}=y_{2}^{2}=y_{2}y_{3}=y_{3}^2=0$, \\ 
$I=\left( 0,0,0,1\right) ,$ \\ 
${\mathcal{F}}=\left( 
\begin{array}{cccc}
0 & 0 & 0 & 1 \\ 
0 & 0 & 1 & 1 \\ 
0 & 1 & 0 & 1 \\ 
1 & 1 & 1 & 1
\end{array}
\right) $
\end{tabular}
$
\end{tabular}
& $
\begin{tabular}{l}
$y_{\mu }=(1,x,y,z),$ \\ 
${\mathcal{F}}_{\mathrm{P\rightarrow G}}=\left( 
\begin{array}{cccc}
0 & 0 & 0 & 1 \\ 
0 & 0 & 1 & 1 \\ 
0 & 1 & 0 & 1 \\ 
1 & 1 & 1 & 1
\end{array}
\right) $
\end{tabular}
$
\end{tabular}
}}

{\small {
\begin{tabular}{lll}
\begin{tabular}{l}
Algebra NF, $\epsilon x=0,\epsilon y=0,\epsilon z=0$
\end{tabular}
& 
\begin{tabular}{l}
Dual algebra, integral, \\ 
Fourier transform
\end{tabular}
& 
\begin{tabular}{l}
Standard dual \\  Fourier transport
\end{tabular}
\\ \hline
\begin{tabular}{l}
(NF,E*),\quad $c[B_{+}]^{\ast },$\quad NF.1: \\ 
\begin{tabular}{l}
$\Delta x=1\otimes x+x\otimes 1,$ \\ 
$\Delta y=1\otimes y+x\otimes y+x\otimes z$ \\ 
$\qquad\quad +y\otimes 1+y\otimes x+z\otimes x,$ \\ 
$\Delta z=1\otimes z+x\otimes y+x\otimes z$ \\ 
$\qquad\quad+y\otimes x+z\otimes 1+z\otimes x,$ \\ 
$Sx=x,\ Sy=z,\ Sz=y$
\end{tabular}
\end{tabular}
& 
\begin{tabular}{l}
$
\begin{tabular}{l}
$y_{0}=1,\ y_{1}y_{2}=y_{2}+y_{3}=y_{1}y_{3}$, \\ 
$y_{1}^{2}=y_2^2=y_3^2=y_{2}y_{3}=0$,\\
$I=\left( 0,0,1,1\right) ,$ \\ 
${\mathcal{F}}=\left( 
\begin{array}{cccc}
0 & 0 & 1 & 1 \\ 
0 & 0 & 1 & 0 \\ 
1 & 0 & 0 & 0 \\ 
1& 1 & 0 & 0
\end{array}
\right) $
\end{tabular}
$
\end{tabular}
& $
\begin{tabular}{l}
$y_{\mu }=(1,x,y+z,y)$, \\ 
${\mathcal{F}}_{\mathrm{NF\rightarrow E}}=\left( 
\begin{array}{cccc}
0& 0 & 0 & 1 \\ 
0& 0 & 1 & 1 \\ 
1& 0 & 0 & 0 \\ 
1& 1 &0  & 0
\end{array}
\right) $
\end{tabular}
$ \\ \hline
\begin{tabular}{l}
(NF,NF*),\quad $d_{sl_{2}},$\quad NF.2: \\ 
\begin{tabular}{l}
$\Delta x=1\otimes x+x\otimes 1+y\otimes x+z\otimes x\ ,$ \\ 
$\Delta y=1\otimes y+x\otimes y+x\otimes z+y\otimes 1$ \\ 
$\qquad\quad+y\otimes x+y\otimes y+y\otimes z+z\otimes x\ ,$ \\ 
$\Delta z=1\otimes z+x\otimes y+x\otimes z+y\otimes x$ \\ 
$\qquad\quad+z\otimes 1+z\otimes x+z\otimes y+z\otimes z,$ \\ 
$Sx=x+y,\ Sy=z,\ Sz=y$
\end{tabular}
\end{tabular}
& 
\begin{tabular}{l}
$
\begin{tabular}{l}
$y_{0}=1,y_{1}^{2}=0$,\\
$y_{1}y_{2}=y_{2}+y_{3}=y_{1}y_{3},$\\
$y_{2}y_{1}=y_{1}+y_{2}+y_{3}= y_{3}y_{1}$,\\
$ y_{2}^{2}=y_{2}=y_{2}y_{3}, y_{3}y_{2}=y_{3}=y_{3}^{2}$, \\ 
$I=\left( 0,1,1,1\right) ,$ \\ 
${\mathcal{F}}=\left( 
\begin{array}{cccc}
0 & 1 & 1 & 1 \\ 
1 & 1 & 1 & 0 \\ 
1 & 0 & 0 & 0 \\ 
1 & 1 & 0 & 0
\end{array}
\right) $
\end{tabular}
$
\end{tabular}
& $
\begin{tabular}{l}
$y_{\mu }=$\\
$(1,y+z,x+y,x)$, \\ 
${\mathcal{F}}_{\mathrm{NF\rightarrow NF}}=\left( 
\begin{array}{cccc}
0& 0 & 0 & 1 \\ 
1& 1 & 0 & 1 \\ 
1 & 0 & 0 & 0 \\ 
1 & 0 & 1 & 1
\end{array}
\right) $
\end{tabular}
$
\end{tabular}
}}

\goodbreak

\end{document}